\documentclass[a4paper,12pt]{article}
\usepackage[english]{babel}
\usepackage[T2A]{fontenc}
\usepackage[cp1251]{inputenc}
\usepackage{amsthm}
\usepackage[tbtags]{amsmath}
\usepackage{amsfonts,amssymb}
\begin{document}

\newtheorem{theorem}{Theorem}
\newtheorem{lemma}{Lemma}
\newtheorem{proposition}{Proposition}
\newtheorem{Cor}{Corollary}

\begin{center}
{\bf Incidence Rings: Automorphisms and Derivations}
\end{center}
\begin{center}
Piotr Krylov\footnote{National Research Tomsk
State University, e-mail: krylov@math.tsu.ru .},
Askar Tuganbaev\footnote{National Research University MPEI; Lomonosov Moscow State University; e-mail: tuganbaev@gmail.com .}
\end{center}

\textbf{Abstract.} We study automorphisms, isomorphisms, and derivations of the incidence algebra $I(X,R)$, where $X$ is preordered set and $R$ is an algebra over some commutative ring $T$. 


\textbf{MSC2020 database 16W20, 16S50}

The study is supported by grants of Russian Science Foundation (=RSF), № 23-21-00375, https://rscf.ru/en/project/23-21-00375/ (P.A. Krylov) and № 22-11-00052, https://rscf.ru/en/project/22-11-00052 (A.A. Tuganbaev).

\tableofcontents

\section{Introduction}\label{section1} 

This paper is devoted to automorphisms, isomorphisms and derivations of the incidence algebra $I(X,R)$ for an arbitrary preordered set $X$ and an arbitrary algebra $R$ over some commutative ring $T$. We use and develop an approach that has proven itself well in describing automorphisms of formal matrix rings with zero trace ideals in \cite{KryN18}, \cite{KryT21}, \cite{KryT22}. Namely, an incidence algebra is written as a splitting extension of some ideal $M$ by means of some subring $L$. In other words, there is a direct decomposition of $I(X,R)=L\oplus M$. Here $L$ is the product of  square matrix rings over $R$ and $M$ is the product of rectangular matrix groups over $R$. We note that the ideal $M$ is contained in the Jacobson radical of the ring $I(X,R)$, while $M$ is a nilpotent ideal in the case of a formal matrix ring. Automorphisms and derivations of the algebra $I(X,R)$ are represented by certain $2\times 2$ matrices with respect to the direct sum $L\oplus M$. This helps a lot when studying automorphisms and derivations of the algebra $I(X,R)$.

A retrospective of studies of automorphisms, isomorphisms and derivations of the ring $I(X,R)$ is presented in the book \cite{SpiO97}. Later, these maps of the incidence rings were considered in a fairly large number of papers. At the end of each of Sections 2-4 of this work, the content of some these papers is briefly disclosed.

Other linear (and nonlinear) maps of incidence rings are also actively studied. So, antiautomorphisms and involutions were studied, see \cite{BruFS11}, \cite{BruL11}, \cite{BruFS12}, \cite{BruFS14}, \cite{ForP21b}, \cite{ForP22}, \cite{DiKL06}. In \cite{BruL11}, it is also contained an overview of the results on involutions. Jordan and Lie isomorphisms were considered (see \cite{BruFK17}, \cite{For21KS}, \cite{XiaW14}).

Much attention is paid to derivations of incidence rings, both ordinary and Lie derivations and Jordan derivations (see \cite{ForP21}, \cite{Jon95}, \cite{KayKW19}, \cite{Khr22}, \cite{Xia15}, \cite{YanW21}). 
We pay attention to a very informative introduction in  \cite{YanW21} with history of derivation studies and an extensive list of references.

For rings of generalized matrices, there are studies of their commuting, self-centralizing and other close mapings, as well as a number of related issues (see \cite{LiW12}, \cite{LiWF18}, \cite{XiaW14}). It may be interesting and advisable to study similar maps of incidence rings. Moreover, it is known how the center of an incident ring is arranged.

Incidence algebras were defined by Rota in his familiar paper \cite{Rot64} as a tool for solving some problems of combinatorics and, above all, for studying generalizations of the M\"obius inversion formula in number theory in a unified way. Over time, incidence algebras themselves, regardless of their applications in combinatorics and other areas of mathematics, have turned out to be a meaningful algebraic object. Many works are dedicated to them, including the book \cite{SpiO97}.

If $X$ is a finite set, then the ring $I(X,R)$ is isomorphic to a certain subring of the matrix ring $M(n,R)$. The ring $I(X,R)$ with finite set $X$ is often called \textsf{ring of structural matrices} (see, for example, \cite{DasW96}). It is known \cite{SmiW94} that the ring of structural matrices is representable as a ring of block triangular matrices. Thus, it is one of the types of formal matrix rings.

Formal matrix rings (they also say <<rings of generalized matrices>>) and incidence rings are closely related and represent characteristic and important algebraic structures. In addition, when studying incidence rings, it is possible to use the accumulated rich material about formal matrix rings and their linear maps. Various questions about formal matrix rings have long been considered in the framework of theory of rings and modules. Some of the results obtained are reflected in the book \cite{KryT17}.

In this paper, we consider only associative rings with a non-zero unit.

If $S$ is a ring (or an algebra), then $U(S)$ is its \textsf{group of invertible elements}, $C(S)$ is its\textsf{center} and $J(S)$ is the \textsf{Jacobson radical}, $M(n,S)$ is the usual ring of all $n\times n$ matrices with values in $S$.

Next, $\text{Aut}S$ is the \textsf{automorphism group} of the algebra $S$, $\text{In(Aut }S)$ is the \textsf{subgroup of inner automorphisms}, $\text{Out}S$ is the \textsf{group of outer automorphisms} of the algebra $S$, i.e. $\text{Out}S$ is the factor group of $\text{Aut}S/\text{In(Aut}S)$.

Using $\text{Der}S$ we denote the \textsf{group} (or the \textsf{module}) of \textsf{derivations} of the algebra $S$, $\text{In(Der}S)$ is the \textsf{subgroup of inner derivations}, $\text{Out}S$ is the \textsf{group of outer derivations} of the algebra $S$, i.e. $\text{Out}S$ is the quotient group $\text{Der}S/\text{In(Der}S)$.\\ \textbf{(The same symbols $\text{Out}S$ for groups of outer automorphisms and outer derivations are used in different places of the text, so this will not lead to confusion. The same applies to a number of other designations related to $\text{Aut }S$ and $\text{Der }S$.)}

Let $R$ be another ring, $A$ be an $R$-$S$-bimodule and let $\alpha$, $\gamma$ be two automorphisms of the rings $R$ and $S$, respectively. You can define a new bimodule structure on $A$ by setting
$$
x\circ a=\alpha(x)a,\; a\circ y=a\gamma(y)\; \text{ for all } \, x\in R,\,y\in S,\, a\in A.
$$
This bimodule is usually denoted by $_{\alpha}A_{\gamma}$.

The semidirect product of the groups $A$ and $B$ is denoted by $A\leftthreetimes B$.

\section*{\textbf{CHAPTER 1. PRELIMINARIES}}\label{chapter1}
\addtocontents{toc}{\textbf{CHAPTER 1. PRELIMINARIES}\par}

\section{On Preordered Sets}\label{section2}

We briefly outline some initial information about preordered sets (one can get acquainted with them in more detail in  \cite{SpiO97}).

Let $X$ be an arbitrary set and let $\le$ be a reflexive  transitive relation on $X$. In this case, the system $\langle X,\le\rangle$ is called a \textsf{preordered set}, and $\le$ is a \textsf{preorder} on $X$. If the relation $\le$ is also antisymmetric, then $\langle X,\le\rangle$ is a \textsf{partially ordered set}.

We further assume that $\langle X,\le\rangle$ is a preordered set. For any two elements $x,y\in X$, we denote by $[x,y]$ the set of $\{z\in X\,|\,x\le z\le y\}$. It is called an \textsf{interval} in $X$. An interval of the form $[x,x]$ is denoted by $[x]$. There are the following two useful properties of intervals.

\textbf{(a)} For any $y,z\in [x]$, we have the relation $[y,z]=[x]$.

\textbf{(b)} If $x<y$, then  $s<t$ for arbitrary elements $s\in [x]$ and $t\in [y]$.

We define a binary relation $\sim$ on $X$ by setting $x\sim y$ $\Leftrightarrow$ $x\le y$ and $y\le x$. It is clear that $\sim$ is an equivalence relation on $X$. The corresponding equivalence classes have the form $[x]$ for all possible $x\in X$. It follows from $(b)$ that the preorder relation $\le$ is consistent with the equivalence relation $\sim$. Consequently, the induced relation $\le$ appears on the quotient set $\overline{X}=X/\sim$, and $\langle\overline{X},\le\rangle$ is a partially ordered set.

A directed graph can be associated with the preordered set $X$, as well as with the partially ordered set $\overline{X}$ (we do not take into account the loops that arise in this case). It is more convenient to proceed from the set $\overline{X}$. When necessary, we consider $\overline{X}$ as a simple graph correlated to the oriented graph $\overline{X}$. At the same time, we use such standard concepts of graph theory: a \textsf{connected component}, a \textsf{semipath} and its \textsf{length}, a \textsf{spanning tree} (or a \textsf{spanning}).

We agree that all intervals in $X$ are finite. In this case, $X$ is called a \textsf{locally finite} preordered set.

For an interval in a locally finite partially ordered set, the \textsf{length} of the interval is the largest of the chain lengths in this interval.

In what follows, we agree to denote by $x$ the elements of a partially ordered set $\overline{X}$, i.e. equivalence classes of the form $[x]$. In other words, to denote the class $[x]$, we use some representative of it. This should not lead to confusion. In a particular situation, it is always clear about the elements of which set ($\overline{X}$ or $X$) we are talking about.

\section{Some Ideals and Subbimodules\\ in Incidence Algebras}\label{section3}

Starting from this section, the symbol $R$ denotes an algebra over some commutative ring $T$. However, the ring $T$ itself is almost not used.

An incidence algebra is a certain ring of functions. Let $\langle X,\le\rangle$ be an arbitrary locally finite preordered set. We set $I(X,R)=\{f\colon X\times X\to R\,|\, f(x,y)=0,\,\text{if }\, x\not\le y\}$. The functions are added point-wise. The product of the functions $f,g\in I(X,R)$ is given by the relation
$$
(fg)(x,y)=\sum_{x\le z\le y}f(x,z)\cdot g(z,y)\eqno(*)
$$
for any $x,y\in X$. Since $X$ is a locally finite set, it is possible to write $z\in X$ in $(*)$ instead of $x\le z\le y$. For any $t\in T$ and $x,y\in X$, we also assume $(tf)(x,y)=tf(x,y)$. As a result, we get a $T$-algebra $I(X,R)$, called the \textsf{incidence algebra} or an \textsf{incidence ring} of the preordered set $X$ over the ring $R$. In what follows, the specific algebra $I(X,R)$ is denoted by the symbol $K$.

We introduce some special functions from $I(X,R)$. For a given $x\in X$, we set $e_{[x]}(t,t)=1$ for all $t\in [x]$ and $e_{[x]}(z,y)=0$ for the remaining pairs $(z,y)$. The system $\{e_{[x]}\,|\,x\in X\}$ consists of pairwise orthogonal central in $L$ idempotents (the ring $L$ is defined in the next paragraph). According to the agreement, we write $e_x$ instead of $e_{[x]}$ at the end of section 2.

We define a subring $L$ and an ideal $M$ in $K$. We set $L=\{f\in K\,|\,f(x,y)=0,\, \text{if }\, x\not\sim y\}$ and $M=\{f\in K\,|\,f(x,y)=0,\, \text{ if }\, x\sim y\}$. We have a $T$-module direct sum $K=L\oplus M$, i.e. the ring $K$ is a splitting extension of the ideal $M$ with the use of the subring $L$. The ideal $M$ is naturally considered as an $L$-$L$-bimodule. In addition, $M$ is a nonunital algebra.

Let we have an arbitrary interval $[x]$. We denote by $R_{[x]}$ the set of functions $f\in K$ for which $f(z,y)=0$ if $z\not\sim x$ or $y\not\sim x$. As in the case of idempotents $e_x$, we write $R_x$ instead of $R_{[x]}$. The following relations
$$
R_x=e_xKe_x=e_xLe_x
$$
are true.
We conclude that $R_x$ is a ring with identity element $e_x$. If we go to the restrictions of functions from $K$ to $[x]\times[x]$, then in fact $R_x$ is the algebra of all functions $[x]\times[x]\to R$ with point-wise addition and the product of the convolution type as in $(*)$. We choose any numbering of the interval $[x]\colon[x]=\{x_1,\ldots,x_n\}$. After that, if the functions $f\in R_x$ match the matrix $(f(x_i,x_j))$, then we come to the isomorphism of the algebras $R_x\cong M(n,R)$. Now We take two different intervals $[x]$, $[y]$ and set
$$
M_{xy}=\{f\in K\,|\,f(s,t)=0,\, \text{if }\, s\not\sim x\, \text{ or } \, t\not\sim y\}.
$$
Then $M_{xy}=e_xKe_y$ and, therefore, $M_{xy}$ is an $R_x$-$R_y$-bimodule. The relation $M_{xy}=e_xMe_y$ is also true.

We clarify that $M_{xy}=0$ if $x\not\le y$. For $x<y$, there is a canonical isomorphism
$$
M_{xy}\cong M(n\times m,R),\quad n=|[x]|,\; m=|[y]|,
$$
with respect to the above isomorphisms $R_x\cong M(n,R)$ and $R_y\cong M(m,R)$. After identifying all algebras $R_x$ with $M(n,R)$ and bimodules $M_{xy}$ with $M(n\times m,R)$, it becomes clear that the actions of the rings $R_x$ and $R_y$ on $M_{xy}$ will be ordinary matrix multiplications. It is also clear that $M_{xy}$ is an $L$-$L$-bimodule. The action of $L$ on $M_{xy}$ is reduced to the action of $R_x$ on the left and $R_y$ on the right.

We note once again that we mean $[x]$ and $[x][y]$ in the subscripts of $R_x$ and $M_{xy}$, respectively (see Section 2).

The product $\prod_{x,y\in X}M_{xy}$ has the $L$-$L$-bimodule structure. Exactly, if $f\in L$ and $(g_{xy})\in\prod_{x,y\in X}M_{xy}$, then
$$
f(g_{xy})=(f_xg_{xy})\, \text{and }\,(g_{xy})f=(g_{xy}f_y),\eqno(1)
$$
where $f_x=e_xfe_x$ and $f_y=e_yfe_y$.

\textbf{Proposition 3.1.} There are canonical algebra isomorphisms $L\cong\prod_{x\in X}R_{x}$, as well as $L$-$L$-bimodule isomorphisms and algebra isomorphisms $M\cong\prod_{x,y\in X}M_{xy}$.

\textbf{Proof.} We define the mapping $\omega\colon L\to \prod_{x\in X}R_x$, assuming $\omega(f)=(f_x)$ for each $f\in L$, where $f_x=e_xfe_x$. Then $\omega$ is an algebra isomorphism.

The mapping $\varepsilon\colon M\to\prod_{x,y\in X}M_{xy}$, $\varepsilon(g)=(g_{xy})$, where $g_{xy}=e_xge_y$, is an isomorphism of $L$-$L$-bimodules and algebras.~$\square$

In what follows, we do not distinguish the corresponding objects with respect to the isomorphisms $\omega$ and $\varepsilon$.

We introduce some ideals and submodules contained in $M$. We pay attention to this. In view of Proposition 3.1 and our agreement, we can write the relation $M=\prod_{x,y\in X}M_{xy}$, in which all the bimodules $M_{xy}$ are non-zero.

Let $V_1=M$, then for every $k\ge 2$, let $V_k=\prod_{x,y\in X}M_{xy}$, where $x$ and $y$ are such that the length of the interval $[[x],[y]]$ of the partly ordered set $\overline{X}$ is not less than $k$. Further, for any $k\ge 1$, we set $W_k=\prod_{x,y\in X}M_{xy}$, and in the product, there are those $M_{xy}$ for which the length of the interval $[[x],[y]]$ is equal to $k$.

We write down the basic properties of subsets $V_k$ and $W_k$. They follow from the definition of multiplication in the bimodule $M$ if we represent it as a product of $\prod_{x,y\in X}M_{xy}$.

\textbf{1.} All $V_k$ are ideals in $K$ and there is a decreasing sequence 
$$
V_1\supseteq V_2\supseteq \ldots\supseteq V_k\ldots\, \text{and } \cap_{k\ge 1}V_k=0.
$$
\textbf{2.} For every $k\ge 1$, the inclusions of $MV_k\subseteq V_{k+1}$ and $V_kM\subseteq V_{k+1}$ are true.

\textbf{3.} All $W_k$ are $L$-$L$-submodules in $M$ and there are relations $M=\prod_{k\ge 1}W_k$ and $V_k=\prod_{i\ge k}W_i$.

\textbf{4.} For each $k\ge 1$, there is a direct $L$-$L$-bimodule decomposition $V_k=W_k\oplus V_{k+1}$.

\textbf{5.} For any $k,\ell\ge 1$, the inclusion $W_k\cdot W_{\ell}\subseteq W_{k+\ell}$ is true.

Like the ring of block matrices, we can consider the ring $K$ as some ring of <<block>> functions. We mean the role of the introduced rings $R_x$ and bimodules $M_{xy}$. Namely, any function $f\in K$ is mapped to a family of functions $\{e_xfe_y\,|\,x,y\in X\}$. Formulas $(1)$ and $(2)$ give the rule of multiplication of functions in block form.

If we assume that $M_{xx}$ is $R_x$, then we can write the relation $K= \prod_{x,y\in X}M_{xy}$. Then we get the equalities $K=\prod_{x\in X}e_xK=\prod_{y\in X}Ke_y$. They allow us to talk about <<rows>> and <<columns>> functions from $K$. When multiplying two functions, <<rows>> are multiplied by <<columns>>.

It is known that the algebra $K$ can be embedded in the matrix algebra $M(\mathcal{M},R)$, where $\mathcal{M}=|X|$\cite[Proposition 1.2.4]{SpiO97}. The equivalence relation $\sim$ to $X$ leads to the partitioning of matrices into certain blocks. When $K$ is embedded in $M(\mathcal{M},R)$, the functions belonging to the rings $R_x$ and the bimodules $M_{xy}$ pass to the corresponding matrix blocks from $M(\mathcal{M},R)$.

We touch on the situation when $X$ is a finite set. We can assume that $X=\{1,2,\ldots,n\}$. Let $\le$ be a preorder on $X$. Then the incidence algebra $K=I(X,R)$ is isomorphic to a certain subalgebra in $M(n,R)$. We identify $K$ with this matrix subalgebra. The latter is denoted by $M(n,\le,R)$ or $M(n,B,R)$, where $B=(b_{ij})$ is a Boolean matrix corresponding to the preorder $\le$. Usually $K$ is called the \textsf{ring of structural matrices} associated with the preorder $\le$ or the Boolean matrix $B$ (see \cite{DasW96}). When embedding the ring $K$ in $M(n,R)$, it is easy to notice that
$$
K=\{(a_{ij})\in M(n,R)\,|\,i\not\le j\,\Rightarrow a_{ij}=0\}=
$$
$$
=\{(a_{ij})\in M(n,R)\,|\,b_{ij}=0\,\Rightarrow a_{ij}=0\}.
$$
When studying rings of structural matrices, they are immediately considered as subrings in $M(n,R)$. There is a permutation $\tau$ of degree $n$ such that $\tau(B)$ is an upper block triangular matrix. There is a ring isomorphism $M(n,B,R)\cong M(n,\tau(B),R)$~\cite{SmiW94}; information about permutation actions on matrices is given in \cite[section 4.1]{KryT17}. So, rings of structural matrices are one of the types of formal matrix rings.

\textbf{Remark 3.2.} The algebra $K$ can be viewed as an algebra of functions and as an abstract ring represented as a splitting extension of $L\oplus M$. These two approaches can be called <<functional>> and <<abstract>>; they are, of course, equivalent. We use both of these approaches.

The discrete topology on $R$ induces a topology on the algebra $I(X,R)$\cite[Section 1.2]{SpiO97}. We introduce another topology using the ideals $V_i$ defined above.

\textbf{Proposition 3.3.} The algebra $K$ is a complete topological algebra in a Hausdorff topology, the base of the neighborhood of zero of which consists of the ideals $V_i$, $i\ge 1$.

\textbf{Proof.} Let $(a_i)_{i\ge 1}$ be a pure Cauchy sequence in $K$. We write $a_i=b_i+c_i$, where $b_i\in L$, $c_i\in M$ and $K=L\oplus M$. It is clear that $b_i=b_j$ for all $i$ and $j$. We denote this common element by $b$. It remains to make sure that the sequence $(c_i)_{i\ge 1}$ converges. By considering the relation
$$
V=\prod_{i\ge 1}W_i,\; V_i=W_i\oplus V_{i+1},\; i\ge 1,
$$
we conclude that all $c_i$, except $c_1$, have the same components in $W_1$, say $d_1$. Similarly, all elements $c_i$, starting from $c_3$, have the same components $d_2$ in $W_2$, and so on. We form an element $d=(d_i)_{i\ge 1}$ of $M$. The relation $\lim_{i\to\infty}c_i=d$ is true.~$\square$

It makes sense to consider linear maps of the algebra $I(X,R)$ that are continuous in the topology from Proposition 3.3. We can also talk about linear maps that leave all submodules $W_k$ in place. The inclusions are $W_k\cdot W_{\ell}\subseteq W_{k+\ell}$ from Property 5 allow us to say that the algebra $I(X,R)$ has some analogue of graduation.

\section{Invertible Elements}\label{section4}

The book \cite{SpiO97} describes the invertible elements and the Jacobson radical of the algebra $I(X,R)$ for the case of a partially ordered set $X$ and a commutative ring $R$. This description is transferred to arbitrary algebras $I(X,R)$, which is done, for example, in \cite{Vos80}. We present these facts here in a slightly different form for completeness. Without explanation, we use the decomposition $K=L\oplus M$ obtained at the beginning of Section 3.

\textbf{Proposition 4.1.} For any element $f\in M$, the element $1+f$ is invertible.

\textbf{Proof.} We make sure of the existence of an element $g\in M$ such that
$$
(1+f)(1+g)=1=(1+g)(1+f).
$$
We redefine the elements $1+f$ and $1+g$ as $f$ and $g$, respectively. Now $f$ is a function such that $f(x,x)=1$ and $f(x,y)=0$ if $x\ne y$ and $x\sim y$ for any $x,y\in X$. The function $g$ has the same properties.

The action of the function $g$ on products of the form $[x]\times[x]$ is known. You need to specify its action on the product $[x]\times [y]$ for any $x,y\in X$ with $x<y$. We do this by induction over the lengths of the intervals $[[x],[y]]$ of the partially ordered set $\overline{X}$.

Let the length of the interval $[[x],[y]]$ be equal to 1. For an arbitrary pair $(s,t)\in [x]\times [y]$, we write the relations
$$
fg(s,t)=\sum_{s\le z\le t}f(s,z)\cdot g(z,t)=
\sum_{z\in [x]}f(s,z)\cdot g(z,t)+
$$
$$
+\sum_{z\notin [x],\; z\notin [y]}f(s,z)\cdot g(z,t)+
\sum_{z\in [y]}f(s,z)\cdot g(z,t)=g(s,t)+c+f(s,t)=0.
$$
Therefore, we have
$$
g(s,t)=-f(s,t)-c,\eqno(*)
$$
and $c=0$. Hence, $g(s,t)=-f(s,t)$.

Suppose that relation $(*)$ holds for all pairs $(s,t)$ such that $(s,t)\in [x]\times[y]$ and the length of the interval $[[x],[y]]$ does exceed $k-1$, where $k\ge 2$. It is assumed that the corresponding element $c$ from $(*)$ is known.

Now let the length of the interval $[[x],[y]]$ be equal to $k$ and $(s,t)\in [x]\times [y]$. For a pair $(s,t)$, one can write equalities such as above, as well as relation $(*)$. In the sum $\sum_{z\not in [x],\; z\notin [y]}f(s,z)\cdot g(z,t)$, the length of each interval $[[z],[t]]$ is less than $k$. Therefore, the value of the function $g$ on pairs $(z,t)$ is known. Therefore, the element $c$ of $(*)$ is known.

So, the desired function $g$ is given, and it is given in such a way that $fg =1$. Similarly, we can find the function $h$ with the property $hf=1$. Therefore, $g=h=f^{-1}$.~$\square$

Using Proposition 4.1 and familiar properties of the Jacobson radical, several results can be obtained.

\textbf{Corollary 4.2.} For an incidence algebra $K$, the following assertions hold.

\textbf{1.} The relation $J(K)=\{f\in K\,|\,f(x,y)\int(R)$ holds for all $x,y\in X$ such that $x\sim y\}$.

\textbf{2.} The relations $M\subseteq J(K)$ and $J(K)=J(L)\oplus M$ hold, where $K=L\oplus M$.

Based on this corollary, it is now easy to prove the following theorem.

\textbf{Theorem 4.3.} For a function $f\in K$, the assertions written below are equivalent.

\textbf{1)} The function $f$ is invertible.

\textbf{2)} If $f=g+h$, where $g\in L$ and $h\in M$, then the function $g$ is invertible in $L$.

\textbf{3)} The restriction of the function $f$ to $[x]\times[x]$, equal to $e_xfe_x$, is an invertible function in the ring $R_x$ for all $x\in X$ (see Section 3 about the characters $e_x$ and $R_x$).

We agree that $\text{In}_1(\text{Aut }K)$ (resp., $\text{In}_0(\text{Aut }K)$) denotes the subgroup of inner automorphisms of the algebra $K$ defined by invertible elements of the form $1+d$, $d\in M$ (respectively, invertible elements of the algebra $L$). The first subgroup is normal in $\text{Aut}K$.

\textbf{Proposition 4.4.} There is a semidirect decomposition
$$
\text{In }(\text{Aut }K)=\text{In}_1(\text{Aut }K)\leftthreetimes \text{In}_0(\text{Aut }K).
$$
The proof of the proposition is contained in \cite[Section 4]{KryT21}.

\section*{\textbf{CHAPTER 2. AUTOMORPHISM GROUPS OF\\ $\qquad$INCIDENCE ALGEBRAS}}\label{chapter2}
\addtocontents{toc}{\textbf{CHAPTER 2. AUTOMORPHISM GROUPS OF\\ \mbox{\,}$\quad\;$INCIDENCE ALGEBRAS}\par}

\section{Two Conditions for the Algebra $K$}\label{section5}

Until the end of the work, we use all the notation and take into account all the conventions adopted in Chapter 1 with respect to the preordered set $X$ and the incidence algebra $I(X,R)$. We denote $I(X,R)$ by $K$.

We specifically recall that the algebra $K$ is a splitting extension. Namely, $K=L\oplus M$, where $L$ and $M$ are defined in Section 3, $L$ and $M$ are the subring and the ideal, respectively.

In the paper \cite{KryT21}, the authors consider automorphisms of an arbitrary algebra $A$, which is a splitting extension of its ideal $N$ using some algebra $S$, i.e. $A=S\oplus N$, where $\oplus$ is a group direct sum.

To an automorphism $\varphi$ of the algebra $A$, we relate a matrix $\begin{pmatrix}
\alpha&\gamma\\
\delta&\beta
\end{pmatrix}$, where
$$
\alpha\colon S\to S,\; \beta\colon N\to N,\; \gamma\colon N\to S,\; \delta\colon S\to N
$$
are $T$-module homomorphisms and
$$
\varphi(a+b)=\begin{pmatrix}
\alpha&\gamma\\ \delta&\beta\end{pmatrix}\begin{pmatrix}a\\b\end{pmatrix}=(\alpha(a)+\gamma(b))+(\delta(a)+\beta(b))
$$
for all $a\in S$ and $b\in N$. It is assumed that $\gamma=0$. In any case, all the results for this <<triangular>> case are applicable to the subgroup of automorphisms of the algebra $A$ whose matrices are the form
$\begin{pmatrix}
\alpha&0\\ \delta&\beta\end{pmatrix}$.

Let $\varphi=\begin{pmatrix}
\alpha&0\\
\delta&\beta
\end{pmatrix}$ be an automorphism of the algebra $A$. Then the mappings $\alpha$, $\delta$ and $\beta$ satisfy the equalities $(1)$ from \cite[Section 3]{KryT21}. In particular, $\alpha$ is an automorphism of the algebra $S$ and $\beta$ is an automorphism of the algebra $N$ (as a non-unital algebra).

Conversely, if $\alpha\colon S\to S$, $\beta\colon N\to N$, $\delta\colon S\to N$ are mappings satisfying the relations $(1)$ from \cite[section 3]{KryT21}, then the following is true.
\\ The transformation of the algebra $A$, given by the matrix
$\begin{pmatrix}
\alpha&0\\
\delta&\beta
\end{pmatrix}$, is an automorphism of the algebra $A$.

We note that if $N^2=0$ (i.e. the algebra $A$ is a trivial extension), then $\delta$ is a derivation of the algebra $S$ with values in the bimodule $_{\alpha}N_{\alpha}$ and $\beta$ is an $S$-$S$-bimodule isomorphism $N\to {}_{\alpha}N_{\alpha}$. In general, in \cite{KryT21}, it is assumed that $N$ is a nilpotent ideal.

In \cite{KryN18}, \cite{KryT21}, \cite{KryT22}, the above approach is applied to the study of automorphism groups of formal matrix rings with zero trace ideals (these include rings of block-triangular matrices). If $A$ is the specified matrix ring, then there is a splitting extension $A=S\oplus N$, where $S$ is a subring of diagonal matrices and $N$ is the ideal of matrices with zeros on the main diagonal. What is of great importance, $N$ is a nilpotent ideal.

We return to the incidence algebras. Still, let $K$ be the incidence algebra. According to Section 3, it can be written as a splitting extension $L\oplus M$. The ideal $M$ is not necessarily nilpotent. Nevertheless, after some modifications, the approach to computing the automorphism group developed in \cite{KryN18}, \cite{KryT21}, \cite{KryT22} is also applicable to incidence algebras. This is largely possible due to the fact that the ideal $M$ is contained in the Jacobson radical of the ring $K$ (corollary 4.2).

In what follows, we do not distinguish between an automorphism and its corresponding matrix. Sometimes, for brevity, we write a <<triangular automorphism $\varphi$>> if $\varphi=\begin{pmatrix}
\alpha&0\\\delta&\beta\end{pmatrix}$, and a <<diagonal automorphism $\varphi$>> if $\varphi=\begin{pmatrix}
\alpha&0\\ 0&\beta\end{pmatrix}$.

We define one group homomorphism and several automorphism groups. Let $f\colon\text{Aut }K\to\text{Aut}L$ be a homomorphism such that $f(\varphi)=\alpha$ for each automorphism $\varphi=\begin{pmatrix}
\alpha&0\\ \delta&\beta\end{pmatrix}$. Next, let $\Lambda$ be the subgroup of diagonal automorphisms and let
$$
C=\{\varphi=\begin{pmatrix}
\alpha&0\\ 0&\beta\end{pmatrix}\in \Lambda\,|\,\alpha R_x=R_x\, \text{for each }\,x\in X\}\,.
$$
We also denote by $\Psi$ the normal subgroup in $\text{Aut}K$ consisting of automorphisms of the form $\begin{pmatrix}
1&0\\0&\beta\end{pmatrix}$, and we denote by $\Psi_0$ the subgroup of all inner automorphisms corresponding to the invertible central elements of the ring $L$. Here $\Psi_0$ is a normal subgroup in $\Psi$. Let $\Omega$ be the image of the homomorphism $f$ and
$$
\Omega_1=\{\alpha\in \Omega\,|\,\alpha R_x=R_x\, \text{for all }\,x\in X\}\,.
$$
Finally, we denote by $\Phi$ the normal subgroup
$$
\{\varphi\in \text{Aut}K\,|\,\varphi=\begin{pmatrix}
\alpha&0\\ \delta&\beta\end{pmatrix},\alpha\in\text{In(Aut }L)\}
$$
of the group $\text{Aut }K$.

Information about the groups defined is very important for understanding the structure of the group $\text{Aut}K$.

In the problem of describing and calculating the group $\text{Aut}K$, the following research directions can be distinguished.

\textbf{1.} Calculation of the normal subgroup $\Phi$.

\textbf{2.} Calculation of the normal subgroup $\Psi$.

\textbf{3.} Calculation of the group $\Omega$.

We formulate two conditions for the algebra $K$ (also see 
\cite[Sections 8,9]{KryT21}, \cite{KryT22}).

\textbf{(I)} For any $\varphi\in\text{Out}K$, the relation $\varphi M=M$ is true, i.e. every automorphism is triangular.

Automorphisms can <<move>> the rings $R_x$ and the bimodules $M_{xy}$. We indicate one condition that eliminates such a phenomenon.

\textbf{(II)} For any $\varphi\in\text{Out}K$ and every $x\in X$, the inclusion of $\varphi(e_x)\in e_y+M$ holds for some $y\in X$.

Condition \textbf{(II)} is used in \cite{KryT21} and \cite{KryT22}.

Now we introduce several conditions on an arbitrary ring $S$.

\textbf{Definition 5.1.}\\
\textbf{1.} A ring $S$ is said to be \textsf{indecomposable} if 1 is its only non-zero central idempotent.

\textbf{2. \cite{AnhW13}, \cite{BirHKP00}.} An idempotent $e$ of the ring $S$ is said to be \textsf{semi-central} if $(1-e)Se=0$. A ring $S$ is said to be \textsf{strongly indecomposable} if $1$ is its only non-zero semi-central idempotent.

The matrix sring over an indecomposable ring is itself indecomposable. A strongly indecomposable ring is not isomorphic to any formal matrix ring $\begin{pmatrix}A&N\\0&B\end{pmatrix}$.

We put together several conditions for the ring $S$.

\textbf{(1)} $S$ is an indecomposable ring.

\textbf{(2)} $S$ is a strongly indecomposable ring.

\textbf{(3)} The quotient ring $S/J(S)$ is indecomposable.

\textbf{(4)} For any idempotent $e\in S$, the relation $(1-e)Se=0$ implies the relation $eS(1-e)=0$.

\textbf{(5)} For any two orthogonal idempotents $e,f\in S$, the relation $fSe=0$ implies the relation $eSf=0$.

There are the following relations between these conditions:
$$
(5)\Rightarrow (4),\quad (3)\Rightarrow (2)\Rightarrow(1)
$$
(see Lemmas 5.7 and 5.8).

Information about the action of isomorphisms, in particular automorphisms, on the product of indecomposable rings is very useful. Let us formulate and prove one general fact.

\textbf{Lemma 5.2.} Let the ring products $S=\prod_{j\in J}S_j$ and $R=\prod_{i\in I}R_i$ be given and let all rings $S_j$ and $R_i$ be indecomposable. If $\alpha\colon S\to R$ is an arbitrary ring isomorphism, then for each $j\in J$ there is an subscript $i\in I$ such that $\alpha S_j=R_i$.

\textbf{Proof.} Let $f_j$ and $e_i$ be the identity elements of the rings $S_j$ and $R_i$, respectively. We take some subscript $\ell\in J$. Then $\alpha(f_{\ell})=(e_p)\in R$, where $p$ runs through some subset of $P\subseteq I$ (we took into account that the rings $R_i$ are indecomposable). So $\alpha$ maps $S_{\ell}$ in $\prod_{p\in P}R_p$. If the restriction of $\alpha$ to $S_{\ell}$ is not an isomorphism, then there is an element $y\in S$ that has the zero component in $S_{\ell}$ and such that $\alpha(y)\in\prod_{p\in P}R_p$ and $\alpha(y)\ne 0$. However, then the relation $f_{\ell}y=0$ implies the relation $\alpha(f_{\ell})\alpha(y)=0$. This is a contradiction, since $\alpha(f_{\ell})$ is the identity element of the ring $\prod_{p\in P}R_p$. We conclude that $\alpha|_{S_{\ell}}$ is an isomorphism. Since the ring $S_{\ell}$ is indecomposable, this is possible only if $|P|=1$.~$\square$

The symmetric analogue of Lemma 5.2 is true for the isomorphism $\alpha^{-1}$. Therefore, there is a bijection $\tau\colon J\to I$ such that $\alpha(f_j)=e_{\tau(j)}$ and $\alpha(S_j)=R_{\tau(j)}$, $j\in J$.

\textbf{Corollary 5.3.} In the situation of Lemma 5.2, the isomorphism $\alpha$ acts on the product $\prod_{j\in J}S_j$ coordinate-wise.

\textbf{Proof.} We explain the statement of the corollary. Let $\alpha\colon S\to R$ be an isomorphism and let $\tau\colon J\to I$ be a bijection that $\alpha$ induces. We take some element $y=(y_j)\in S$ and write $\alpha(y)=(x_i)\in R$. The coordinate action means that $\alpha(y_j)=x_{\tau(j)}$ for each $j\in J$.

So let $y=(y_j)\in S$ and $\alpha(y)=(x_i)\in R$. We take an arbitrary subscript $\ell\in J$ and write the decompositions $S=S_{\ell}\oplus\prod_{j\ne\ell}S_j$ and $R=R_{\tau(\ell)}\oplus \prod_{i\ne \tau(\ell)}R_i$. Next, we write $y=y_{\ell}+a$, where $a\in \prod_{j\ne\ell}S_j$, and
$$
\alpha(y)=\alpha(y_{\ell})+\alpha(a),\,\text{where}\, \alpha(y_{\ell})\in R_{\tau(\ell)}.
$$
It is enough to make sure that $\alpha(a)\in\prod_{i\ne\tau(\ell)}R_i$. Otherwise, we have the relation $e_{\tau(\ell)}\alpha(a)\ne 0$. Where do we get $\alpha^{-1}(e_{\tau(\ell)}\alpha(a))=f_{\ell}a\ne 0$; this is a contradiction.~$\square$

Now we turn to the automorphisms of the ring $S$ equal to $\prod_{j\in J}S_j$, where all rings $S_j$ are indecomposable. It turns out that the automorphisms of this ring <<rearrange>> the rings $S_j$. Based on this, we select (however, in a non-canonical way) in the group $\text{Aut}S$ a subgroup isomorphic to some group of bijections of the set $J$. To do this, we partition the set of rings $S_j$ ($j\in J$) into equivalence classes with respect to the isomorphism of rings. We also automatically get a partition of the set $J$.

Let $P$ be some equivalence class in $J$. We fix some element in $P$ and denote it by the integer 1. For each $s\in P$, we choose some isomorphism $\varepsilon_{s1}\colon S_1\to S_s$, and $\varepsilon_{11}$ is the identity automorphism of the ring $S_1$. Next, we assume $\varepsilon_{1s}=\varepsilon_{s1}^{-1}$ and $\varepsilon_{pq}=\varepsilon_{p1}\cdot\varepsilon_{1q}\colon S_q\to S_p$ for all $p,q\in P$.

Let $\sigma$ be some bijection of the set $P$. We define the automorphism $\alpha_{\sigma}$ of the algebra $\prod_{q\in P}S_q$ assuming that $\alpha_{\sigma}$ acts on $S_q$ as the isomorphism $\varepsilon_{\sigma(q)q}$. After that, it is clear that $(\alpha_{\sigma})^{-1}=\alpha_{\sigma^{-1}}$. If $\tau$ is another bijection of the set $P$, then the relation $\alpha_{\tau\sigma}=\alpha_{\tau}\alpha_{\sigma}$ holds.
Therefore, the mapping $\sigma\to \alpha_{\sigma}$ is an isomorphic embedding of the bijection group of the set $P$ into the group $\text{Aut}(\prod_{q\in P}S_q)$. We denote the image of this attachment by $\Sigma_P$. Now we set $\sum=\prod_P\Sigma_P$, where $P$ runs through all equivalence classes of the set $J$ introduced above. The automorphism group $\sum$ of the ring $S$ is isomorphic to some bijection group of the set $J$. These are exactly the bijections that leave in place every equivalence class $P$ of the set $J$. We identify $\sum$ with this group of bijections, i.e. we assume that $\alpha_{\tau}=\tau$.

We define one more subgroup of the group $\text{Aut }S$. Let $\Gamma=\{\mu\in\text{Aut }S\,|\,\mu S_j=S_j$ for every $j\in J\}$. Here $\Gamma$ is a normal subgroup of the group $\text{Aut }S$ and $\Gamma\cap \sum=\langle 1\rangle$. We take an arbitrary automorphism $\varphi\in\text{Aut }S$. By Lemma 5.2, $\varphi$ induces a bijection $\tau$ of the subscript set $J$. We have an inclusion $\varphi\alpha_{\tau}^{-1}\in \Gamma$. We denote $\varphi\alpha_{\tau}^{-1}$ by $\mu$ and obtain $\varphi=\mu\alpha_{\tau}$, where $\mu\in\Gamma$, $\alpha_{\tau}\in\sum$.

We formalize everything recently said in the form of the following assertion.

\textbf{Corollary 5.4.} The group $\text{Aut }S$ is equal to the semidirect product $\Gamma\leftthreetimes\sum$, where $\Gamma=\prod_{j\in J}\text{Aut }S_j$, $\sum$ is some group of bijections of the set $J$. In addition, there exists an isomorphism $\text{Out }S\cong\prod_{j\in J}\text{Out }S_j\leftthreetimes\sum$.

We will repeatedly refer to one result of Azumaya. Therefore, we formulate it.

\textbf{Theorem 5.5 \cite{Azu51}.} Let $e_1,\ldots,e_n$ and $e_1^*,\ldots,e_n^*$ be two systems of pairwise orthogonal idempotents of some ring $S$ such that $e_i-e_i^*\in J(S)$ for every $i$. Then there exists an inner automorphism $\mu$ of the ring $S$ with $\mu(e_i)=e_i^*$, $i=1,\ldots,n$.

Let again $K=I(X,R)$ be an incidence algebra written as a splitting extension $L\oplus M$. In the remaining part of the section, we consider conditions \textbf{(I)} and \textbf{(II)} for the algebra $K$. They are closely related to each other.

\textbf{Lemma 5.6.} For an incidence algebra $K$, the following assertions hold.

\textbf{1.} Condition \textbf{(II)} implies Condition \textbf{(I)}.

\textbf{2.} For an indecomposable ring $R$, conditions \textbf{(I)} and \textbf{(II)} are equivalent.

\textbf{Proof.}\\ \textbf{1.} On the contrary, we assume that there exists an automorphism $\varphi$ of the algebra $K$ such that $\varphi M\not\subseteq M$. We take some element $a\in M$ with $\varphi(a)=b\notin M$ and choose an idempotent $e_y$ with $e_ybe_y\ne 0$. It follows from Condition \textbf{(II)} that $\varphi^{-1}(e_y)=e_x+d$ for some $x\in X$ and $d\in M$. It follows from Theorem 5.5 that there exists an inner automorphism $\mu$ with $\mu(e_x+d)=e_x$. We have $\mu\varphi^{-1}(e_y)=e_x$. As a result, we obtain impossible relations
$$
0\ne\mu\varphi^{-1}(e_ybe_y)=e_x\mu(a)e_x=0.
$$
Therefore, that Condition \textbf{(I)} holds.

\textbf{2.} We take an arbitrary automorphism 
$\varphi=\begin{pmatrix}
\alpha&0\\ \delta&\beta\end{pmatrix}$ of the algebra $K$ and some idempotent $e_x$. It follows from Lemma 5.2 that $\alpha(e_x)=e_y$ for some $y$. Therefore, we obtain
$$
\varphi(e_x)=\alpha(e_x)+\delta(e_x)=e_y+\delta(e_x)\in e_y+M.\quad \square
$$

We highlight some situations when one of the conditions \textbf{(I)} and \textbf{(II)} holds.

There is an obvious case. Let $R$ be a semiprimitive ring, i.e. $J(R)=0$. Then $J(L)=0$ and $J(K)=M$ (see Corollary 4.2). It is clear that $K$ satisfies Condition \textbf{(I)}.

\textbf{Lemma 5.7.} If all rings $R_x$ satisfy condition $(5)$, then Condition \textbf{(I)} is satisfied for the algebra $K$.

\textbf{Proof.} We assume the contrary. Let $\varphi$ be an automorphism of the algebra $K$ such that $\varphi M\not\subseteq M$. Then $\varphi(a)=b\notin M$, where $a\in M$. There exists an idempotent $e_y$ such that $e_ybe_y\ne 0$. We write down $\varphi^{-1}(e_y)=f+m$, where $f^2=f\in L$, $m\in M$ and $f\ne 0$. Let an inner automorphism $\mu$ transfer $f+m$ to $f$ (Theorem 5.5). Then $(\mu\varphi^{-1})(e_y)=f$. We represent $f$ in the form $f=(f_x)$, $f_x\in R_x$, as agreed in Section 3.

We set $S=fKf$. The restriction of the mapping $\mu\varphi^{-1}$ to $R_y$ is a ring isomorphism
$$
R_y=e_yKe_y\to fKf=S.
$$ 
Since $e_ybe_y\ne 0$, we have $f\mu(a)f\ne 0$. Therefore, $f_x\mu(a)f_z\ne 0$ for some $x,z$ such that $x\not\sim z$. Consequently,
$$
f_xKf_z\ne 0,\quad f_xKf_z=f_xfKff_z=f_xSf_z\ne 0,
$$
and $f_x,f_z\in S$. Since $S\cong R_y$, it follows from condition $(5)$ that $f_zSf_x\ne 0$. Therefore, $e_zKe_x\ne 0$ and $z<x$. Therefore, $e_xKe_z= 0$ and $f_xSf_z= 0$. This is a contradiction.~$\square$

In connection with the proved lemma, a natural question arises. For which rings $R$ does the matrix ring $M(n,R)$ satisfy Condition $(5)$ for all $n\ge 2$? The answer is positive in each of the following cases:

\textbf{a)} $R$ is a local ring;

\textbf{b)} $R$ is a left (or right) principal) ideal domain;

\textbf{c)} $R$ is a commutative Dedekind ring.

What unites these rings is that finitely generated modules over them are generators.

Lemmas 5.7 and 5.6 imply that if $R$ is one of the listed rings, then the conditions \textbf{(I)} and \textbf{(II)} are satisfied for the algebra $I(X,R)$.

Before passing to Condition \textbf{(II)}, we formulate one simple result.

\textbf{Lemma 5.8.} Let $S$ be an arbitrary ring. If the factor ring $S/J(S)$ is indecomposable, then $S$ is a strongly indecomposable ring. In particular, $S$ is indecomposable.

\textbf{Proof.} We assume that there exists a non-zero idempotent $e\ne 1$ with $(1-e)Se=0$. The ring $S$ can be identify with the formal matrix ring
$\begin{pmatrix}eSe&eS(1-e)\\ 0&(1-e)S(1-e)\end{pmatrix}$. This leads to a ring isomorphism
$$
S/J(S)\cong eSe/J(eSe)\oplus(1-e)S(1-e)/J((1-e)S(1-e)).
$$
This is contradiction, since the ring $S/J(S)$ is indecomposable.~$\square$

\textbf{Lemma 5.9.} If factor ring $R/J(R)$ is indecomposable, then Condition \textbf{(II)} holds for the ring $K$.

\textbf{Proof.} The factor ring $R_x/J(R_x)$ is indecomposable for every $x\in X$.

Taking into account the relation $J(K)=J(L)\oplus M$ from Corollary 4.2, we can identify the factor rings $K/J(K)$ and $L/J(L)$. We identify the last ring with the product $\prod_{x\in X}R_x/J(R_x)$.

We fix an automorphism $\varphi\in\text{Aut }K$. Let $\overline{\varphi}$ be the automorphism of the algebra $K/J(K)$ or $L/J(L)$ induced by the automorphism $\varphi$. It follows from Lemma 5.2 that for every idempotent $e_x$, there is an idempotent $e_y$ with $\overline{\varphi}(\overline{e}_x)=\overline{e}_y$, where the line above $e_x$ (resp., $e_y$) denotes the reduce class with respect to $J(R_x)$ (resp., $J(R_y)$) or, equivalently, with respect to $J(K)$.

Returning to the ring $K$, we obtain the relation $\varphi(e_x)=e_y+a+b$, where $a\in J(L)$, $b\in M$, and $e_y+a$ is an idempotent in $L$. Next we write down $a=c+d$, $c\in R_y$, $d\in \prod_{z\not\sim y}R_z$. Since $d\in J(L)$ and $d$ is an idempotent, we have $d=0$. It follows from $c\in J(R_y)$ that $e_y+c$ is an invertible idempotent in $R_y$. Consequently, $e_y+c=e_y$ and $c=0$. Thus, it is proved that $\varphi(e_x)=e_y+b$ and $b\in M$, which is required.~$\square$

At the end of the section, we answer such a natural question. What conditions on the ring $R$ and $I(X,R)$ are the most general and acceptable for calculating the group $\text{Aut}K$, where $K=I(X,R)$?

In fact, such a condition can be called Condition \textbf{(II)}. Indeed, if it holds, then Condition \textbf{(I)} holds as well (see Lemma 5.6). 
That is, all automorphisms of the algebra $K$ are triangular. Condition \textbf{(II)} allows predicting the nature of the action of automorphisms on the rings $R_x$ and the bimodule $M_{xy}$. This means the following. Let $\varphi=\begin{pmatrix}\alpha&0\\ \delta&\beta\end{pmatrix}\in\text{Aut }K$. Even when the rings $R_x$ are decomposable and Lemma 5.2 is not directly applicable to the product $\prod_{x\in X}R_x$, Condition \textbf{(II)} guarantees the following. For any $x\in X$, there exists an element $y\in X$ such that $\alpha(e_x)=e_y$ and $\alpha R_x=R_y$. This allow using automorphisms of the partially ordered set $\overline{X}$. Let also $\delta=0$, i.e. $\varphi$ is a diagonal automorphism. Then for each pair $x,y$, $x\not\sim y$, we have $\beta M_{xy}=M_{st}$, where $e_s=\alpha(e_x)$ and $e_t=\alpha(e_y)$.

In fact, such a condition can be called Condition \textbf{(II)}. Indeed, if it is true, then Condition \textbf{(I)} is also true (Lemma 5.6). That is, all automorphisms of the algebra $K$ are triangular. In addition, Condition \textbf{(II)} allows predicting the nature of the action of automorphisms on the rings $R_x$ and the bimodules $M_{xy}$. This means the following. Let $\varphi=\begin{pmatrix}\alpha&0\\ \delta&\beta\end{pmatrix}\in\text{Aut }K$. Even when the rings $R_x$ are decomposable and Lemma 5.2 is not directly applicable to the product $\prod_{x\in X}R_x$, Condition \textbf{(II)} guarantees the following. For any $x\in X$, there exists such a $y\in X$ that $\alpha(e_x)=e_y$ and $\alpha R_x=R_y$. This allows using automorphisms of the partially ordered set $\overline{X}$. Let also $\delta=0$, i.e. $\varphi$ is a diagonal automorphism. Then for each pair $x,y$ with $x\not\sim y$, we have $\beta M_{xy}=M_{st}$, where $e_s=\alpha(e_x)$ and $e_t=\alpha(e_y)$.

We indicate two main situations when Condition \textbf{(II)} holds.

\textbf{1.} The factor ring $R/J(R)$ is indecomposable (Lemma 5.9).

We recall that the rings $R$ and $R_x$ are also indecomposable by Lemma 5.8.

\textbf{2.} Condition \textbf{(I)} holds and the ring $R$ is indecomposable (Lemma 5.6).

By considering Lemmas 5.7, 5.9 and the text after the proof of Lemma 5.7, we make up the list of rings $R$ such that the algebra $I(X,R)$ satisfies Condition \textbf{(II)}:

\textbf{a)} $R$ is a local ring; 

\textbf{b)} $R$ is a left (or right) principal ideal domain; 

\textbf{c)} $R$ is a commutative Dedekind ring; 

\textbf{d)} $X$ is partially ordered set and $R$ loes not have non-zero idempotents besides 1; 

\textbf{e)} a ring $R$ such that factor ring $R/J(R)$ is indecomposable (local rings refer to \textbf{e)}). 

\textbf{Agreement.} \textbf{In Sections 6, 8 and 9, we assume that the algebra $K$, equal to $I(X,R)$, satisfies Condition \textbf{(II)}.}\\ 
For example, this is the case if $R$ is a one of the rings indicated in \textbf{a)}--\textbf{e)}.

\section{Main Isomorphisms and Decompositions for Group $\text{Aut }K$}\label{section6}

As before, let $K$ denote some incidence algebra $I(X,R)$ represented in the form $K=L\oplus M$, as at the beginning of Section 3.
We take an arbitrary automorphism $\varphi=\begin{pmatrix}\alpha&0\\ \delta&\beta\end{pmatrix}$ of the algebra $K$ (all automorphisms of $K$ are triangular by Lemma 5.6, see the agreement just adopted. Here $\alpha\in \text{Aut }L$. For every $x\in X$, there is $y$ such that $\alpha(e_x)=e_y$ (this is proved at the end of the previous section). If $s\ne x$ and $\alpha(e_s)=e_t$, then $e_ye_t=0$ and $y\ne t$ if $e_xe_s=0$. Now let $z\in X$. Taking into account that
$\varphi^{-1}=\begin{pmatrix}\alpha^{-1}&0\\ *&*\end{pmatrix}$, we obtain $\alpha(e_z)^{-1}=e_u$ for some $u\in X$. Therefore, $\alpha(e_u)=e_z$. Therefore, that $\alpha$ induces bijection $\tau\colon \overline{X}\to \overline{X}$, $x\to y$. In fact, $\tau$ is an automorphism of the partially ordered set $\overline{X}$. We verify this property.

Let $x,y\in X$ and $x<y$. Then $M_{xy}=e_xMe_y\ne 0$. We choose an inner automorphism $\mu$ such that
$$
\mu(e_{\tau(x)})+\delta(e_x))=e_{\tau(x)},\quad 
\mu(e_{\tau(y)})+\delta(e_y))=e_{\tau(y)}.
$$
(Theorem 5.5). We have relations
$$
\mu\varphi(M_{xy})=\mu\varphi(e_x)M\mu\varphi(e_y)=
e_{\tau(x)}Me_{\tau(y)}=M_{\tau(x)\tau(y)}\ne 0.
$$
Therefore, $\tau(x)<\tau(y)$.

Conversely, if it is given that $\tau(x)<\tau(y)$, then it is similarly possible to get $x<y$. So, $\tau\in\text{Aut}\overline{X}$. In this case, the mapping $\varphi\to\tau$ is a group homomorphism of $p\colon\text{Aut}K\to\text{Aut}\overline{X}$. It is convenient to write $\tau_{\alpha}$ instead of $\tau$. There is also a homomorphism $r\colon\Omega\to\text{Aut}\overline{X}$, $\alpha\to\tau_{\alpha}$ (the group $\Omega$ is defined at the beginning of Section 5). We note that $\alpha(e_x)=e_{\tau_{\alpha}(x)}$ for any $x\in X$.

So, if the automorphism $\varphi=\begin{pmatrix}\alpha&0\\\delta&\beta\end{pmatrix}$ is given, then the automorphism $\alpha$ <<rearranges>> of the ring $R_x$. In addition, if $\varphi$ is a diagonal automorphism, then we can specify exactly the action of $\beta$ on the bimodule $M_{xy}$.

\textbf{Proposition 6.1.} For an automorphism $\begin{pmatrix}\alpha&0\\ \delta&\beta\end{pmatrix}$ of the algebra $K$ there are the following assertions.

\textbf{1.} The automorphism $\alpha$ of the algebra $L$ permutes the rings $R_x$ according to some automorphism $\tau$ of the partially ordered set $\overline{X}$.

\textbf{2.} If $\delta=0$, i.e. $\varphi$ is a diagonal automorphism, then the automorphism $\beta$ of the $L$-$L$-bimodule $M$ permutes bimodules $M_{xy}$ according to the automorphism $\tau$. In addition, the restriction of $\beta$ to $M_{xy}$ is an isomorphism of bimodules $M_{xy}\to M_{\tau(x)\tau(y)}$ (with respect to ring isomorphisms $\alpha|_{R_x}\colon R_x\to R_{\tau(x)}$, $\alpha|_{R_y}\colon R_y\to R_{\tau(y)}$).

We recall on the equality from Proposition 4.4:
$$
\text{In }(\text{Aut }K)=\text{In}_1(\text{Aut }K)\leftthreetimes \text{In}_0(\text{Aut }K).
$$

We also recall that the groups from the following theorem appeared at the beginning of Section 5. In the following theorem, statement \textbf{1)} can be interpreted as the possibility to diagonalize any automorphism of the algebra $K$.

\textbf{Theorem 6.2.} Let algebra $K$ satisfies equality \textbf{II}.\\ The following relations are true.

\textbf{1.} $\text{Aut }K=\text{In}_1(\text{Aut }K)\leftthreetimes \Lambda$.

\textbf{2.} $\text{Ker }f=\text{In}_1(\text{Aut }K)\leftthreetimes \Psi$.

\textbf{3.} $\Phi=\text{In(Aut }K)\cdot\Psi=\text{In}_1(\text{Aut }K)\leftthreetimes(\text{In}_0(\text{Aut }K))\cdot\Psi$.

\textbf{Proof.}\\ \textbf{1.} We take automorphism $\varphi\in \text{Aut }K$ and let $\tau$ be the automorphism of the set $\overline{X}$ corresponding to $\varphi$ (it was given before Proposition 6.1). We form the element $v=(v_{st})$ of $K$, where
$$
v_{st}=\varphi(e_{\tau^{-1}(t)})(s,t) \,\text{ for any }\, s,t\in X.
$$
We remark that $ve_x=\varphi(e_{\tau^{-1}(x)})e_x$ for every $x\in X$. We have $v_{tt}=1$ and, therefore, the element $v$ is invertible in $K$ (Proposition 4.1 or Theorem 4.3).

So, for every $z\in X$, there is relation $ve_{\tau(z)}=\varphi(e_z)e_{\tau(z)}$. In addition, the relation $\varphi(e_z)v=\varphi(e_z)e_{\tau(z)}$ holds. Thus, we have 
$$
\varphi(e_z)v=ve_{\tau(z)},\quad v^{-1}\varphi(e_z)v=e_{\tau(z)}.
$$
Let $\mu$ be inner automorphism defined by the element $v$. Then
$$
\mu\in \text{In}_1(\text{Aut }K),\quad \mu\varphi(e_z)=e_{\tau(z)}\; \text{ for every }\, z\in X.
$$
By setting $\gamma=\mu\varphi$, we obtain $\varphi=\mu^{-1}\gamma$, where $\mu^{-1}\in\text{In}_1(\text{Aut }K)$, $\gamma\in\Lambda$.

We verify that $\text{In}_1(\text{Aut }K)\cap \Lambda=\langle 1\rangle$.
Let $\nu\in\text{In}_1(\text{Aut }K)\cap \Lambda$ and $\nu$ is defined by an element $1+d$, where $d\in M$. Then $\nu=\begin{pmatrix}1&0\\ 0&\gamma\end{pmatrix}$ for some $\gamma$. For an arbitrary element $a\in L$ we have
$$
\nu(a)=(1+d')a(1+d)=a+d'a+ad+d'ad,
$$
where $1+d'$ is the inverse element for $1+d$. Therefore, $d'a+ad+d'ad=0$. In particular,
$$
d'e_x+e_xd+d'e_xd=0,\; e_xd'e_x+e_xd+e_xd'e_xd=e_xd=0 
$$
for every $x\in X$. Therefore, $d=0$. Thus, $\nu=1$. It can be argued that the semidirect decomposition from \textbf{1} really exists.

\textbf{2.} The equality from \textbf{2} follows from \textbf{1} and inclusion $\text{In}_1(\text{Aut }K)\subseteq\text{Ker }f$.

\textbf{3.} It follows from decomposition in \textbf{1} that 
$\Phi=\text{In}_1(\text{Aut }K)\leftthreetimes(\Phi\cap\Lambda)$. It remains to verify that $\Phi\cap\Lambda=\text{In}_0(\text{Aut }K)\cdot\Psi$. It is sufficient to prove that the left part is contained in right part. Let $\varphi=\begin{pmatrix}\alpha&0\\ 0&\beta\end{pmatrix}\in \Phi\cap\Lambda$, where $\alpha\in\text{In}(\text{Aut }L)$. If $\alpha$ is defined by an element $v\in U(L)$, then let $\mu$ be the inner automorphism defined by the same element $v$. By setting $\psi=\mu^{-1}\varphi$, we obtain $\mu\psi=\varphi$, where $\mu\in \text{In}_0(\text{Aut }K)$, $\psi\in\Psi$.~$\square$

We gather several interesting useful equalities and isomorphisms.

\textbf{Proposition 6.3.} The following relations and isomorphisms hold.

\textbf{1.} $\Psi\cap\text{In(Aut }K)=\Psi\cap\text{In}_0(\text{Aut }K)=\Psi_0$.

\textbf{2.} $\Lambda/(\text{In}_0(\text{Aut }K))\cdot\Psi\cong\Omega/\text{In(Aut }L)$, $$C/\text{In}_0(\text{Aut }K)\cdot\Psi\cong\Omega_1/\text{In(Aut }L).$$

\textbf{3.} $\Phi/\text{Ker }f\cong\text{In}_0(\text{Aut }K)/\Psi_0$.

\textbf{4.} $\Phi/\text{In(Aut }K)\cong\Psi/\Psi_0$.

\textbf{5.} $\Lambda/C\cong\Omega/\Omega_1$.

\textbf{Proof.}\\
\textbf{1.} We verify the non-obvious inclusion $\Psi\cap\text{In(Aut }K)\subseteq\Psi_0$ holds. We take some automorphism $\mu=\mu_1\mu_0$, where $\mu\in\Psi$, $\mu_1\in\text{In}_1(\text{Aut }K)$, $\mu_0\in\text{In}_0(\text{Aut }K)$. Let $\mu_1$ is defined by an element $1+d$, $d\in M$, and $\mu_0$ is defined by an element $v\in U(L)$. Then $\mu$ is defined by an element $v+vd$, and $\mu^{-1}$ is defined by an element $v^{-1}+d'$, where $d'\in M$.

For any element $a\in L$ we have
$$
\mu(a)=v^{-1}av+(v^{-1}avd+d'av+d'avd).
$$
Since $\mu\in\Psi$, we have $v^{-1}av=a$ and $av=va$. Therefore, $v\in C(L)$. The expression in parentheses is equal to zero. Just as at the end of the proof of Theorem 6.2(1), it can be obtained that $d=0$. Therefore, $\mu_1=1$ and $\mu=\mu_0$. Since $v$ is an invertible central element, we have $\mu\in\Psi_0$.

\textbf{2.} First, we remark that an inclusion $\text{In(Aut }L)\subseteq\Omega$ is always true. We denote by $\pi$ the canonical epimorphism $\Omega\to\Omega/\text{In(Aut }L)$. The kernel of the homomorphism $\pi f|_{\Lambda}\colon \Lambda\to\Omega/\text{In(Aut }L)$ is equal to $\Phi\cap\Lambda$ and, consequently, is equal to $\text{In}_0(\text{Aut }K)\cdot\Psi$ (see the proof of Theorem 6.2(3)).

\textbf{3,4,5.} By considering \textbf{1} and Theorem 6.2, we obtain isomorphisms 
$$
\Phi/\text{Ker }f\cong(\text{In}_0(\text{Aut }K)\cdot\Psi)/\Psi\cong\text{In}_0(\text{Aut }K)/\Psi_0\cong\text{In(Aut }L).
$$
Next we have
$$
\Phi/\text{In(Aut }K)=(\text{In}(\text{Aut }K)\cdot\Psi)/\text{In}(\text{Aut }K) \cong\Psi/\Psi_0.
$$
It is also not difficult to prove isomorphism from \textbf{5}.~$\square$

We formulate the main result about the structure of the group $\text{Aut }K$.

\textbf{Theorem 6.4.} Let algebra $K$ satisfies equality \textbf{II}.
 
\textbf{a.} The following isomorphisms are true.

$\quad$\textbf{1.} $\text{Aut }K/\text{Ker }f\cong\Omega\cong\Lambda/\Psi$ and $C/\Psi\cong\Omega_1$.

$\quad$\textbf{2.} $\text{Aut }K/\Phi\cong\Omega/\text{In}(\text{Aut }L)$.

\textbf{b.} $\text{Out }K$ has a normal subgroup $H$ such that
$$
H\cong\Psi/\Psi_0,\quad \text{Out }K/H\cong\Omega/\text{In}(\text{Aut }L).
$$
\textbf{c.} If the relation $\Psi=\Psi_0$ holds, then
$$
\Phi=\text{In(Aut }K),\quad \text{Out }K\cong\Omega/\text{In}(\text{Aut }L).
$$
\textbf{Proof.}\\
\textbf{a.} The isomorphisms from \textbf{1} and \textbf{2} directly follow from Theorem 6.2 and Proposition 6.3.

\textbf{b.} We take $\Phi/\text{In(Aut }K)$ as a subgroup of $H$. Everything we need, follows from Proposition 6.3 and the  isomorphism from \textbf{2}.

\textbf{c.} If $\Psi=\Psi_0$, then it follows from Theorem 6.2 that $\Phi=\text{In(Aut }K)$.~$\square$

The following conclusion can be drawn: if it is possible to find the structure of the groups$\Psi$ and $\Omega$, then the structure of the groups $\text{Aut}K$ and $\text{Out}K$ will be known in a certain sense.

In papers \cite{KryT21} and \cite{KryT22}, for formal matrix algebra, various information about groups $\Psi$ and $\Omega$ is obtained (the definitions of the groups $\Psi$ and $\Omega$ for formal matrix rings is similar to the definitions in the given paper). 
Sections 8 and 9 of the given paper are also devoted to the groups $\Omega$ and $\Psi$.

At the end of the section, we present various information about the action of automorphisms on the bimodule $M$.

\textbf{Proposition 6.5.}\\ 
\textbf{1.} Diagonal automorphisms act on the product $M=\prod_{x,y\in X}M_{xy}$ coordinate-wise. At the same time the relation $\varphi W_k=W_k$ holds for every diagonal automorphism $\varphi$ and $k\ge 1$.

\textbf{2.} Let $\varphi\in \text{Aut }K$. Then $\varphi V_k=V_k$ for every $k\ge 1$.

\textbf{Proof.} The bimodules $M_{xy}$, $W_k$ and the ideals $V_k$ are defined in Section 3.

\textbf{1.} Let $\varphi=\begin{pmatrix}\alpha&0\\ 0&\beta\end{pmatrix}$ be a diagonal automorphism. For any $x,y\in X$, $x<y$, and each $m=(m_{xy})\in M$ we have relations
$$
\varphi(m_{xy})=\beta(m_{xy})=\beta(e_xme_y)=\alpha(e_x)\beta(m)\alpha(e_y)=
$$
$$
=e_{\tau(x)}\beta(m)e_{\tau(y)}\in M_{\tau(x)\tau(y)},
$$
where $\tau$ is the automorphism partially of the ordered of the set $\overline{X}$ induced by $\alpha$ (see Proposition 6.1 and the text before it). On the other hand, we write down
$$
\varphi(m)=(n_{xy})\in \prod_{x,y\in X}M_{xy}. 
$$
We have to verify that $\varphi(m_{xy})=\beta(m_{xy})=n_{\tau(x)\tau(y)}$ for any non-equivalent $x,y\in X$ (cf. Corollary 5.3).

We not write out all the details of further reasoning. 
They are based on the presence of direct decompositions of the bimodule $M$ for any non-equivalent $x,y\in X$:
$$
M=M_{xy}\oplus N_{xy},\, \text{ where }\, N_{xy}=
\prod_{s\not\sim x\,\text{ or }\,t\not\sim y}M_{st}
$$
and the second similar decomposition $M=M_{\tau(x)\tau(y)}\oplus N_{\tau(x)\tau(y)}$. When checking the equality of $\varphi W_k=W_k$, we should also take into account the coincidence of the lengths of the intervals $[[x],[y]]$ and $[[\tau(x)],[\tau(y)]]$ in $\overline{X}$.

\textbf{2.} We write down an arbitrary automorphism $\varphi\in \text{Aut }K$ in the form $\varphi=\mu\psi$, where $\mu\in \text{In}_1(\text{Aut }K)$ and $\psi$ is a diagonal automorphism (Theorem 6.2). Let $\mu$ is defined by an element $1+d$, where $d\in M$. Then we have
$$
\mu(m)=(1+d')m(1+d)=m+(d'm+md+d'md)
$$
for every $m\in M$, where $1+d'$ is the inverse element for $1+d$. If $m\in V_k\setminus V_{k+1}$, then $d'm+md+d'md\in V_{k+1}$ (property 2 in Section 3). Therefore, $\mu(m)\in V_k$. So, $\mu V_k=V_k$. The relation $\psi V_k=V_k$ follows from \textbf{1}. Since $\varphi=\mu\psi$, we also have $\varphi V_k=V_k$.~$\square$

\section{On Group $\text{Aut }X$}\label{section7}

We present a number of facts about the automorphism group $\text{Aut}X$ of the preordered set $\langle X,\le\rangle$. Some notation and terminology concerning the set of $X$ are in Section 2. In particular, $\sim$ is an equivalence relation on $X$, $\overline{X}$ is the factor set $X/\sim$, which is a partially ordered set.

Every automorphism $\tau\in \text{Aut }X$ induces automorphism $\overline{\tau}\in \text{Aut }\overline{X}$. This leads to the group homomorphism $k\colon \text{Aut }X\to\text{Aut }\overline{X}$. An automorphism $\tau\in \text{Aut }X$ is said to be \textsf{inner} if it maps every interval of the form $[x]$, 
$x\in X$, into itself. All inner automorphisms form a normal subgroup of $\text{In(Aut }X)$ in the group $\text{Aut }X$. The factor group $\text{Aut }X/\text{In(Aut }X)$ is called the \textsf{group of outer automorphisms} of the set $X$; it is denoted by $\text{Out }X$.

The kernel of the homomorphism $k$ coincides with $\text{In(Aut }X)$, and the image of $k$ is contained in the subgroup $G$ equal to $\{p\in \text{Aut }\overline{X}\,|\,|[x]|=|p([x])|$ for any $x\in X\}$. We define (not canonically) a homomorphism $j\colon G\to \text{Aut }X$. 
To do this, we divide the set of all intervals of the form $[x]$ into classes of intervals of the same cardinality.
Let $J$ be one from such classes. In $J$, we fix some interval $[z]$. Next for every interval $[x]\in J$ we choose bijection $\varepsilon_{[x],[z]}\colon [z]\to[x]$, and $\varepsilon_{[z],[z]}$ is the identity mapping. Next we set 
$$
\varepsilon_{[z],[x]}=\varepsilon_{[x],[z]}^{-1},\quad
\varepsilon_{[x],[y]}=\varepsilon_{[x],[z]}\varepsilon_{[z],[y]}.
$$
We do the same with each class of sets of the same cardinality. 
We take an arbitrary automorphism $q\in G$. 
Let $\tau_q$ be a mapping $X\to X$ such that on the interval $[x]$ from some class $J$, the mapping $\tau_q$ coincides with $\varepsilon_{q([x]),[x]}$. Then $\tau_q\in\text{Aut}X$ and the mapping $q\to\tau_q$, $q\in G$, is a homomorphism of $j\colon G\to\text{Aut}X$. At the same time, the equality $k_j=1_G$ holds. This implies that the groups $\text{Im }k$ and $G$ coincide. Therefore, we can identify these groups with $\text{Out }X$. Next we can identify these groups with $\text{Im }j$ which is a  defined (but non-canonically) subgroup in $\text{Aut }X$. Thus the following fact is true.

\textbf{Proposition 7.1.} There is a (non-canonical) semidirect decomposition of groups $\text{Aut }X=\text{In}(\text{Aut }X)\leftthreetimes\text{Out }X$.

The above result and the reasoning before it are ideologically close to the considerations after the proof of Corollary 5.3.

Based on the connected components of the preordered set $X$, another decomposition of the group $\text{Aut}X$ can be obtained. Let $X_i$ $(i\in I)$ be all pairwise distinct connected components of the set $X$. We divide the set of all connected components into equivalence classes with respect to an isomorphism of preordered sets $X_i$. We also get the corresponding partition of the set $I$. In essentially the same way as in the situation with the ring product $\prod_{j\in J}S_j$ after the proof of Corollary 5.3, it is possible to define (in a non-canonical way) some automorphism group $\Sigma_X$ of the set $X$. It is isomorphic to the group of those bijections of the set $I$ that leave each of its equivalence classes in place. Now we define another normal subgroup $\Gamma_X=\{\sigma\in\text{Aut }X\,|\,\sigma X_i=X_i$ for every $i\in I\}$. The following assertion holds.

\textbf{Proposition 7.2.} There ar the following semidirect decompositions:
$$
\text{Aut }X=\Gamma_X\leftthreetimes\Sigma_X,\,\text{ where }\,
\Gamma_X=\prod_{i\in I}\text{Aut }X_i,\; \text{Out }X=\prod_{i\in I}\text{Out }X_i\leftthreetimes\Sigma_X.
$$

How does partitioning the set $X$ into connected components affect the structure of the ring $I(X,R)$? The following result will allow us to apply Corollary 5.3 to the group $\text{Aut}K$.

\textbf{Proposition 7.3.} A ring $I(X,R)$ is indecomposable if and only if the ring $R$ is indecomposable and the set $X$ is connected.

\textbf{Proof.} It is not difficult to verify \textbf{necessity} conditions of the proposition. 

\textbf{Sufficiency.} We assume that the ring $I(X,R)$ is decomposable into a non-trivial direct sum of ideals. Let $g,h$ be central idempotents corresponding too this decomposition. By considering Theorem 5.5, we can assume that $g,h\in L$ (we recall about decomposition $K=L\oplus M$). Similar to Section 3 (see Proposition 3.1), we write down $g=(g_x)$ and $h=(h_x)$, where $g_x,h_x\in R_x$. Since the ring $R$ is indecomposable, we have either $g_x=1$ or $g_x=0$ for any $x$; the same is true for $h_x$. We set $Y=\{x\,|\,g_x=1\}$ and $Z=\{x\,|\,h_x=1\}$. Then $Y\cap Z=\varnothing$ and $Y\cup Z=\overline{X}$. It follows from $gMh=0=hMg$ that $e_yMe_z=0=e_zMe_y$ for any $y\in Y$, $z\in Z$. Consequently, $y$ and $z$ are not comparable in $\overline{X}$. But the existence such subsets $Y$ and $Z$ contradicts to the property that $X$ is connected.~$\square$

The partitioning the set $X$ into connected components $X_i$ leads to equality $K=I(X,R)=\prod_{i\in I}I(X_i,R)$ after identification.
From Corollary 5.4 and Proposition 7, we obtain the following assertions.

\textbf{Proposition 7.4.} Let $R$ be an indecomposable ring. Then we have a semidirect decomposition 
$$
\text{Aut }K=\Gamma_K\leftthreetimes\Sigma_K,\,\text{ where }\,
\Gamma_K=\prod_{i\in I}\text{Aut }(I(X_i,R)).
$$
We also have an isomorphism
$$
\text{Out }K\cong\prod_{i\in I}\text{Out }K_i\leftthreetimes\Sigma_K,\,\text{ where }\, K_i=I(X_i,R).
$$

So, when studying the incidence rings $I(X,R)$, to some extent we can limit ourselves to the case of a connected set $X$.

At the beginning of the previous section, the homomorphisms $p\colon\text{Aut}K\to\text{Aut}\overline{X}$ and $r\colon\Omega\to\text{Aut}\overline{X}$ were introduced. More importantly, the group $\text{Aut}X$ can be naturally embedded in the group $\text{Aut}K$. Let $\tau\in\text{Aut }X$. We define a mapping $\xi_{\tau}\colon K\to K$, assuming $\xi_{\tau}(f)(x,y)=f(\tau(x),\tau(y))$ for any $f\in K$ and $x,y\in X$. Then $\xi_{\tau}$ is an automorphism and the mapping $\tau\to\xi_{\tau}$ is a group embedding $t\colon\text{Aut}X\to\text{Aut}K$. We identify $\tau$ with $\xi_{\tau}$. Automorphisms of the form $\xi_{\tau}$ are called \textsf{ordinal}. The equality $pt=k$ is also true.

We can also specify how the automorphism $\tau\in\text{Aut}X$ acts on elements from $K$ if they are represented in the vector form as in Section 3. In general we can say that automorphism $\tau$ permutes components of the elements. As a result, the group $\Sigma_X$ can be embedded in the group $\Sigma_K$, and $\Gamma_X$ can be embedded in the $\Gamma_K$ (see Propositions 7.2 and 7.4).

Now we consider the decomposition $\text{Aut }X=\text{In}(\text{Aut }X)\leftthreetimes\text{Out }X$ from Proposition 7.1. We clarify again that we identify the automorphism $\tau$ with the automorphism $\xi_{\tau}$ of the ring $K$ (it is introduced above). An arbitrary automorphism $\tau\in \text{In}(\text{Aut }X)$ acts on $K$ as the conjugation by a matrix $P$, where $P=\prod_{x\in X}P_x$, and $P_x$ is the matrix of some permutation of degree $|[x]|$. Thus, $\tau\in \text{In}_0(\text{Aut }K)$. As for the group $\text{Out }X$, we have an inclusion $\text{Out }X\subseteq\Lambda$, whence $\text{Out }X\cap C=\langle 1 \rangle$.

The restriction of the homomorphism $f$ to $\text{Out }X$ is an embedding $\text{Out }X\to \Omega$ and $\text{Out }X\cap \Omega_1=\langle 1 \rangle$ (of the group $\Lambda$, $C$, $\Omega$ and $\Omega_1$ appeared in Section 5).

It cannot be argued that the subgroup $\text{Out}X$ will be a semi-direct summand in $\Lambda$ or $\Omega$ (and this situation is very characteristic). The obstacle is contained in the fact that there may be such isomorphic rings $R_x$ and $R_y$ that the matrices of these rings have different orders. In Section 9, we will return to this circumstance and indicate how this obstacle can be eliminated under one assumption.

\section{Groups $\Omega$ and $\Omega_1$}\label{section8}

At the end of Section 5, the condition \textbf{II} was fixed on the ring $K$, where $K=I(X,R)$. \textbf{We assume that $K$ satisfies this condition.}

The group $\Omega$ and its normal subgroup $\Omega_1$ were introduced in Section 5. After proving Theorem 6.4, the expediency of calculating the groups $\Omega$ and $\Psi$ for understanding the structure of the group $\text{Aut}K$ is emphasized. Section 10 is dedicated to the group $\Psi$.

The group $\Omega_1$ is the kernel of the homomorphism $r\colon\Omega\to \text{Aut }\overline{X}$; see the beginning of Section 6. Therefore, the factor group $\Omega/\Omega_1$ can be embedded in the group $\text{Aut }\overline{X}$. In the previous section, we have agreed to consider the group $\text{Out }X$ as a subgroup of $\text{Aut }K$ and $\Omega$ such that $\Omega_1\cap \text{Out }X=\langle 1\rangle$. It is important to find out when $\Omega=\Omega_1\leftthreetimes \text{Out }X$. This question will be considered in Section 9. If we take $L$ as the ring $S$ in Corollary 5.4, then we obtain decomposition $\text{Aut }L=\Gamma\leftthreetimes \Sigma$. It is clear that $\Omega_1\subseteq\Gamma$. It is the subgroup $\Omega_1$ that we pay attention to here.

We write down the following questions about the structure of the group $\Omega_1$.

\textbf{1.} Which automorphisms of $\Gamma$ belong to $\Omega_1$?

\textbf{2.} What is the structure of the groups $\Omega_1$ and $\Omega_1/\text{In}(\text{Out}L)$? In connection with the question \textbf{2}, it is useful to refer to Theorem 6.4.

The answer to the first question and, in some cases, to the second question will be given.

We prove one general fact. It generalizes the following result (see \cite[Chapter 2, Proposition 5.2]{Fai73}).

Let $H$ be some algebra and let $\alpha$ and $\gamma$ be two automorphisms of $H$. There exists an $H$-$H$-bimodule isomorphism $H\to {}_{\alpha}H_{\gamma}$ if and only if $\alpha^{-1}\gamma$ is an inner automorphism. The definition of the bimodule ${}_{\alpha}H_{\gamma}$ is given in Section 1.

For positive integers $k$ and $\ell$, we set
$$
P=M(k,R),\; Q=M(\ell,R),\; H=M(c,R),\, \text{ where }\, c=\text{LCM }(k,\ell),
$$
and, finally, $V=M(k\times\ell,R)$. Let $\ell'=c/k$, $k'=c/\ell$.

The ring $H$ can be represented as a ring of block matrices in two ways: as a ring of block matrices over $P$ of order $\ell'$ and as a ring of block matrices over $Q$ of order $k'$. It is also a $P$-$Q$-bimodule of block matrices over $V$ of size $\ell'\times k'$.

Let $\alpha$ and $\gamma$ be automorphisms of the algebras $P$ and $Q$, respectively. They induce automorphisms $\overline{\alpha}$ and $\overline{\gamma}$ of $H$, respectively. Namely, $\overline{\alpha}(A)=(\alpha(a_{ij}))$ for any matrix $A=(a_{ij})\in H$, $a_{ij}\in P$, and similarly for $\overline{\gamma}$. The automorphisms $\overline{\alpha}$ and $\overline{\gamma}$ are called \textsf{ring (block) automorphisms}. We denote them $\alpha$ and $\gamma$. This agreement is already in effect in the following proposition.

\textbf{Proposition 8.1.} If $\alpha\in \text{Aut }P$ and $\gamma\in\text{Aut }Q$, then an isomorphism of $P$-$Q$-bimodules $V\to{}_{\alpha}V_{\gamma}$ exists if and only if $\alpha^{-1}\gamma$ is an inner automorphism of the algebra $H$.

\textbf{Proof.} Let we have a $P$-$Q$-bimodule isomorphism $\beta\colon V\to {}_{\alpha}V_{\gamma}$. The isomorphism $\beta$ induces the $H$-$H$-bimodule isomorphism $\overline{\beta}\colon H\to {}_{\alpha}H_{\gamma}$, $\overline{\beta}(A)=(\beta(A_{ij}))$ for every matrix $A=(A_{ij})$ from $H$. 
It is assumed that the matrix $A$ is represented in the above block form, i.e. $A_{ij}$ are $\ell'\times k'$ blocks.
Consequently, $\alpha^{-1}\gamma$ is an inner automorphism of the algebra $H$.

Now we assume that $\alpha^{-1}\gamma$ is an inner automorphism of the algebra $H$. Therefore, there exists an  $H$-$H$-bimodule isomorphism $\beta\colon H\to {}_{\alpha}H_{\gamma}$.

Let $e_1,\ldots,e_{\ell'}$ and $f_1,\ldots,f_{k'}$ be diagonal matrix units corresponding to two block partitions of matrices from $H$. The following relations
$$
\alpha(e_i)=e_i\; (i=1,\ldots,\ell'),\,\text{ and }\,\gamma(f_j)=f_j\; (j=1,\ldots,k')
$$
are true. Therefore, we can deduce that $\beta$ induces an $A$-$B$-bimodule isomorphism $e_iHf_j\to $ ${}_A(e_iHf_j)_B$, where $A=e_iHe_i$, $B=f_jHf_j$. In other words, $\beta$ actually induces a $P$-$Q$-bimodule isomorphism $V\to{}_{\alpha}V_{\gamma}$.~$\square$

We pay attention to such a circumstance. An automorphism $\alpha$ of $\Gamma$ is contained in $\Omega_1$ exactly when there is a transformation $\beta$ of the algebra $M$, which is both its automorphism and an $L$-$L$-bimodule isomorphism $M\to{}_{\alpha}M_{\alpha}$. This is equivalent to the fulfillment of equality
$$
\beta(cd)=\beta(c)\beta(d),\eqno(1)
$$
where $c\in M_{xz}$, $d\in M_{zy}$, $x<z<y$, and equalities
$$
\beta(ac)=\alpha(a)\beta(c),\; \beta(db)=\beta(d)\alpha(b),\eqno(2)
$$
where $a\in R_x$, $c,d\in M_{xy}$, $b\in R_y$ (it is necessary to take into account the formulas $(1)$ and $(2)$ in Section 3 and Proposition 6.5(1) on coordinate-wise action of automorphisms).

Let $n_x$ be the order of matrices from the ring $R_x$, $x\in X$. We set $c_{xy}=\text{LCM}(n_x,n_y)$ for all $x,y\in X$ such that $x<y$. We denote by $H_{xy}$ the matrix ring $M(c_{xy},R)$. The ring is a block matrix ring over the rings $R_x$ and $R_y$. It also is an $R_x$-$R_y$-bimodule of block matrices over $M_{xy}$. Automorphisms of the rings $R_x$ and $R_y$ are assumed to be ring (block) automorphisms of the rings $H_{xy}$.

\textbf{Theorem 8.2.} An automorphism $\alpha=(\alpha_x)$ of the algebra $L= \prod_{x\in X}R_x$ belongs to the group $\Omega_1$ if and only if $\alpha_x^{-1}\alpha_y$ is an inner automorphism of the algebra $H_{xy}$ for all $x,y\in X$ such that $x<y$.

\textbf{Proof.}\\
\textbf{Necessity.} Let $\alpha\in \Omega_1$. Therefore, there exists an automorphism $\begin{pmatrix}\alpha&0\\ 0&\beta\end{pmatrix}$ of the algebra $K$, where $\beta$ is an $L$-$L$-bimodule isomorphism $M\to{}_{\alpha}M_{\alpha}$. The restriction of $\beta$ to $M_{xy}$ is an $R_x$-$R_y$-bimodule isomorphism $M_{xy}\to{}_{\alpha_x}(M_{xy})_{\alpha_y}$. By Proposition 8.1
$$
\alpha_x^{-1}\alpha_y\in \text{In}(\text{Aut }H_{xy}).
$$

\textbf{Sufficiency.} By Proposition 8.1, we have an $R_x$-$R_y$-bimodule isomorphism $\beta_{xy}\colon M_{xy}\to {}_{\alpha_x}(M_{xy})_{\alpha_y}$, where $x<y$. We define a mapping
$$
\beta=(\beta_{xy})\colon M\to M,\; (d_{xy})\to (\beta_{xy}(d_{xy})),\;
\text{ where }\;(d_{xy})\in M=\prod_{x,y\in X}M_{xy}.
$$
Here $\beta$ is an $L$-$L$-bimodule isomorphism $M\to{}_{\alpha}M_{\alpha}$ (this can be verified with the use of equalities $(2)$). In order for the transformation $\begin{pmatrix}\alpha&0\\ 0&\beta\end{pmatrix}$ of the algebra $K=L\oplus M$ to be its automorphism, it remains to verify that $\beta$ is an automorphism of the algebra $M$. To do this this, it is sufficient to verify that relation $(1)$ holds.

We fix three elements $x,y,z$ such that $x<z<y$. We begin to use one more matrix ring. We set $H=M(d,R)$, where $d=\text{LCM}(n_x,n_z,n_y)$. The ring $H$ is a block matrix ring over each of the rings $H_{xz}$, $H_{zy}$, $H_{xy}$, $R_x$, $R_z$, $R_y$. It can be also considered as a bimodule of block matrices over each of the bimodules $M_{xz}$, $M_{zy}$, $M_{xy}$.

We assume that the automorphisms $\alpha_x$, $\alpha_z$ and $\alpha_y$ are ring (block) automorphisms of the algebra $H$. We consider the isomorphisms $\beta_{xz}$, $\beta_{zy}$ and $\beta_{xy}$ as bimodule (block) isomorphisms $H\to{}_{\alpha_x}H_{\alpha_z}$, $H\to{}_{\alpha_z}H_{\alpha_y}$ and $H\to{}_{\alpha_x}H_{\alpha_y}$ (i.e. we also preserve the previous notation). The products $\alpha_x^{-1}\alpha_z$, $\alpha_z^{-1}\alpha_y$ and $\alpha_x^{-1}\alpha_y$ are inner automorphisms of the algebra $H$. 
They induce (in the sense of Proposition 8.1) the same bimodule isomorphisms $\beta_{xz}$, $\beta_{zy}$ and $\beta_{xy}$ as written above.

These bimodule isomorphisms act as follows (see \cite[Section 2]{KryT21}). There are invertible matrices $u,v,w\in H$, corresponding to the above inner automorphisms, such that there are relations
$$
\beta_{xz}(a)=\alpha_x(a)\alpha_x(u)= \alpha_x(u)\alpha_z(a),
$$
$$
\beta_{zy}(b)= \alpha_z(b)\alpha_z(v)= \alpha_z(v)\alpha_y(b), 
$$
$$
\beta_{xy}(c)= \alpha_x(c)\alpha_z(w)= \alpha_x(w)\alpha_y(c)
$$
for all $a,b,c\in H$.

We verify that $H$ satisfies the relation
$$
\beta_{xy}(cd)=\beta_{xz}(c)\beta_{zy}(d) \eqno(3) 
$$
for any $c,d\in H$. We pay attention to the relation $w=vu$ (we still have to remember that elements $u,v,w$ are defined up to invertible central elements). We have relations
$$
\beta_{xy}(cd)=\alpha_x(cd)\alpha_x(w)=\alpha_x(c)\alpha_x(d)\alpha_x(v)\alpha_x(u).
$$
We also have relations
$$
\beta_{xz}(c)\beta_{zy}(d)=\alpha_x(c)\alpha_x(u)\alpha_z(d)\alpha_z(v)=
$$
$$
=\alpha_x(c)\alpha_x(d)\alpha_x(u)\alpha_z(v)=
$$
$$
=\alpha_x(c)\alpha_x(d)\alpha_x(v)\alpha_x(u).
$$
The equality $(3)$ is proved. In $(3)$, the multiplication of matrices is performed in block form, from which the equality $(1)$ follows.~$\square$

The rest of the section is devoted to the structure problem for the group $\Omega_1$ (see Question 2 at the beginning of the section). This problem seems to be quite complicated. At the same time, information about the structure of the group $C$ is obtained. The latter is defined in Section 5. The group $C$ consists of those diagonal automorphisms $\begin{pmatrix}\alpha&0\\0&\beta\end{pmatrix}$ of the algebra $K$ such that $\alpha_xR_x=R_x$ for all $x\in X$.

We recall that $n_x$ is the order of matrices over the ring $R_x$ (the symbol $n_x$ appeared before Theorem 8.2). Let the natural number $k$ be such that every $n_x$ is a multiple of $k$. For example, $k$ can be the largest common divisor of all integers $n_x$ or it is equal to one. We note that in further considerations we could limit ourselves to the case $k=1$. To do this, you need to represent the ring $I(X,R)$ in the form of some ring of <<block>> functions (see Section 3 about this).

We denote by $P$ the matrix ring $M(k,R)$. Any automorphism $\varepsilon$ of the ring $P$ provides a ring (block) automorphism $\varepsilon_x$ of the $R_x$. Next we obtain an automorphism $\varepsilon_L=(\varepsilon_x)$ of the $L$. An automorphism $\varepsilon$ also induces block automorphism $\varepsilon_M$ of the algebra $M$ and the $L$-$L$-bimodule $M$. As a result, we obtain an automorphism $\overline{\varepsilon}=\begin{pmatrix}\varepsilon_L&0\\ 0&\varepsilon_M\end{pmatrix}$ of the algebra $K$. We call the automorphisms $\varepsilon_L$, $\varepsilon_M$ and $\overline{\varepsilon}$ \textsf{ring (block) automorphisms induced by the automorphism} $\varepsilon$.

We denote by $D$ the subgroup of all ring automorphisms of the algebra $L$. The subgroup of all ring automorphisms of the algebra $K$ is denoted by the same letter $D$. The both these subgroups are canonically isomorphic to the group $\text{Aut }P$. At the same time, $D\subseteq \Omega_1$ in $\text{Aut }L$, $D\subseteq C$ and $D\cap \Psi=\langle 1 \rangle$ in $\text{Aut }K$. Next let $D_0$ be a subgroup of $D$ which is isomorphic to $\text{In}(\text{Aut }P)$. Then $D/D_0\cong\text{Out }P$ and the following lemma is true.

\textbf{Lemma 8.3.} There are relations
$$
D_0=D\cap \text{In}(\text{Aut }L)\, \text{ in }\, \text{Aut }L\, \text{ and}
$$
$$
D_0=D\cap \text{In}_0(\text{Aut }K)\cdot\Psi=
 D\cap \text{In}_0(\text{Aut }K)\, \text{ in } \text{Aut }K.
$$
In addition, the relation $\Psi_0D_0=\Psi D\cap \text{In}_0(\text{Aut }K)$ holds.

We write down the following condition $(*)$ for the algebra $K$:
\begin{enumerate}
\item[]
For every $\alpha=(\alpha_x)\in\Omega_1$, any automorphism $\alpha_x$ is the product in $\text{Aut }R_x$ of an inner automorphism and a ring (block) automorphism.
\end{enumerate}

Condition $(*)$ holds if $R$ is a local ring or a left (or right) principal ideal domain \cite[Corollary 10.5]{KryT22}.

It is convenient to denote the factor group $C/\text{In}_0(\text{Aut }K)$ by $\text{Out}_0K$ (this is because the group $\text{Out}K$ is actually equal to $\Lambda/\text{In}_0(\text{Aut }K)$; also see Proposition 6.3).

We recall that in Section 2 we agreed to consider the preordered set $X$ as a directed graph.

\textbf{Theorem 8.4.} Let $X$ be a connected set and let the algebra $K$ satisfy Conditions \textbf{(II)} and $(*)$. There are the following relations and isomorphisms.

\textbf{1.}
 
\textbf{a.} $\Omega_1=\text{In }(\text{Aut }L)\cdot D$;

\textbf{b.} $\Omega_1/\text{In }(\text{Aut }L)\cong \text{Aut }P$;

\textbf{c.} $C=\text{In}_0(\text{Aut }K)\cdot \Psi\cdot D$.

\textbf{2.} The group $\text{Out}_0K$ contains a normal subgroup of $H$ which is isomorphic to $\Psi/\Psi_0$, and the factor group $\text{Out}_0K/H$ is isomorphic to $\Omega_1/\text{In }(\text{Aut }L)$, i.e. it is isomorphic to the group $\text{Out }P$.

\textbf{Proof.}\\
\textbf{1.} \textbf{a.} We take an arbitrary automorphism $\alpha=(\alpha_x)\in\Omega_1$. By assumption, for every $x\in X$, we have $\alpha_x=\mu_x\rho_x$, where $\mu_x$ is an inner automorphism and $\rho_x$ is a ring (block) automorphism of the ring of matrices $R_x$. We form automorphisms $\mu=(\mu_x)$ and $\rho=(\rho_x)$ of the ring $L$, and we obtain the relation $\alpha=\mu\rho$ in the group $\Gamma$. Since $\alpha,\mu\in \Omega_1$, we have $\rho\in\Omega_1$. Let $\rho_x$ be induced by an automorphism $\varepsilon_x\in \text{Aut }P$, $x\in X$. 
Now we show that it is actually possible to make it so that $\varepsilon_x=\varepsilon_y$ for all $x,y\in X$. In this regard, we note that any automorphism $\rho_x$ can obviously be considered as a ring automorphism for any ring $R_y$. In the sense that $\varepsilon_x$ induces some ring automorphism of the ring $R_y$.

We fix some element $t\in X$. Let $x\in X$. In $X$, we choose a semipath from $t$ to $x$. Applying Theorem 2.2 several times or using the induction on the semipath length, we can make sure that $\rho_x=\nu_x\rho_t$, where $\nu_x$ is some inner ring automorphism of the ring $R_x$.

Here are some small explanations. The equalities of the form $\rho_p=\nu_p\rho_q$ obtained by the above methods of reasoning arise in the ring $H_{pq}$ (see Theorem 8.2). If we set $c=\text{LCM}(n_t,n_{z_1},\ldots,n_{z_m},n_x)$, where $z_1,z_2,\ldots,z_m$ are the vertices of the selected semipath from $t$ to $x$, then we can consider all the appearing ring automorphisms as automorphisms of the algebra $M(c,R)$. Therefore, we can consider them as automorphisms of each ring $R_x$. As a corollary, we obyain $\nu_x\in\text{Aut }R_x$.

Returning to the beginning of the proof, we obtain $\alpha_x=\mu_x\rho_x=\mu_x\nu_x\rho_t$ for every $x\in X$. We set $\gamma=(\mu_x\nu_x)$ and $\rho=(\rho_t)$. Then $\alpha=\gamma\rho$, where $\gamma\in\text{In}(\text{Aut }L)$, $\rho\in D$.

\textbf{b.} By considering Lemma 8.3, we have the following relations and isomorphisms.
$$
\Omega_1/\text{In}(\text{Aut }L)=(\text{In}(\text{Aut }L)\cdot D)/\text{In}(\text{Aut }L)\cong
$$
$$
\cong D/(D\cap \text{In}(\text{Aut }L))=D/D_0\cong\text{Aut }P/\text{In}(\text{Aut }P)=\text{Out }P.
$$
\textbf{c.} Let $\varphi=\begin{pmatrix}\alpha&0\\ 0&\beta\end{pmatrix}\in C$, where $\alpha\in\Omega_1$. Similar to \textbf{a}, we write down $\alpha=\mu\rho$, where $\mu$ is an inner automorphism and $\rho$ is a ring automorphism of the ring $L$. They induce an inner automorphism and a ring automorphism of $K$, respectively. We leave the same designations $\mu$ and $\rho$ for them. Now we set $\psi=(\mu\rho)^{-1}\varphi$ and obtain relations 
$$
\varphi=\mu\rho\psi=\mu\psi'\rho,\;\text{ where }\; \mu\in \text{In}_0(\text{Aut }K), \;\psi,\psi'\in\Psi,\; \rho\in D.
$$

\textbf{2.} As the subgroup $H$, we take the factor group $\text{In}_0(\text{Aut }K)\cdot\Psi/\text{In}_0(\text{Aut }K)$ which is isomorphic to $\Psi/\Psi_0$ by Proposition 6.3. After this, we obtain relations 
$$
\text{Out}_0K/H\cong \text{In}_0(\text{Aut }K)\cdot\Psi\cdot D/\text{In}_0(\text{Aut }K)\cdot\Psi\cong
$$
$$
\cong D/(D\cap \text{In}_0(\text{Aut }K)\cdot\Psi)=D/D_0\cong\text{Out }P
$$
(we also have take into account Lemma 8.3).~$\square$

We gather some partial cases of Theorem 8.4.

\textbf{Corollary 8.5.}
 
\textbf{1.} If we set $k=1$ in Theorem 8.4, i.e. $P=R$, then $D\cong\text{Aut }R$.

\textbf{2.} We assume that the algebra $K$ is such that there exists an integer $n_s$ which divides every $n_x$. Then we can take $R_s$ as the ring $P$ and condition $(*)$ will be satisfied. In particular,
$D\cong\text{Aut }R_s$.

\textbf{3.} Let $X$ be a partially ordered set. Then the assumptions of \textbf{1} and \textbf{2} hold. If additionally the ring $R$ is commutative and only automorphisms of the $R$-algebra $K$ are taken into account, then $\Omega_1=\langle 1 \rangle$ and $D=\langle 1 \rangle$.

\textbf{Proof.} It is not necessary to prove \textbf{1} and \textbf{3}.

\textbf{2.} We assume that all the rings $R_x$ are block matrix rings over the ring $R_s$. Let $\alpha=(\alpha_x)\in \Omega_1$. We take an arbitrary automorphism $\alpha_x$ and fix a semipath from $s$ to $x$. Similar to how the equality from Theorem 8.4(1(a)) was verified, we can obtain the relation $\alpha_x=\mu_x\alpha_s$ in some ring $M(c,R)$, where $\mu_x$ is an inner automorphism, and $\alpha_x,\alpha_s$ are ring (block) automorphisms this the rings. Since $n_s$ divides all $n_z$, $z\in X$, we have that $\alpha_s$ is a ring automorphism in $R_x$. Therefore, $\mu_x\in\text{Aut }R_x$. So, $\alpha_x=\mu_x\alpha_s$, where $\mu_x$ is an inner automorphism, and $\alpha_s$ is a ring automorphism of the ring $R_x$. We have proved that condition $(*)$ holds for the ring $P$ equal to $R_s$.~$\square$

Based on Proposition 7.4, we extend Theorem 8.4 and Corollary 8.5 to incidence algebras $I(X,R)$, where $X$ is an arbitrary preordered set.

As before, let $K=I(X,R)$, $X_i$ ($i\in I$) be all connected components of the set $X$, and let $K_i=I(X_i,R)$ (the algebras $K_i$ appeared in Section 7).

According to Proposition 7.4, we have a semi-direct decomposition
$$
\text{Aut }K=\Gamma_K\leftthreetimes\Sigma_K,\, \text{where }\,
\Gamma_K=\prod_{i\in I}\text{Aut}K_i
$$
and $\Sigma_K$ is some group of bijections of the set $I$. We also have an isomorphism
$$
\text{Out }K\cong\prod_{i\in I}\text{Out }K_i\leftthreetimes\Sigma_K.
$$
We not write out all the details of the argument. We just pay attention to some key points. First, we have the equality $C=\prod_{i\in I}C_i$, where $C_i$ is the subgroup of diagonal automorphisms of the ring $K_i$ that leave the corresponding rings $R_x$ stationary. We have the inclusion $C\subseteq\Gamma_K$.

The analogue of the group $D$ for an arbitrary set $X$ should be defined as follows. 
For each $i\in I$, let the integer $k_i$ be chosen so that every $n_x$ is a multiple of $k_i$, where $x\in X_i$. Next let $P_i=M(k_i,R)$ and let $D$ be a subgroup of all ring (block) automorphisms of the ring $K_i$. Now we set $D=\prod_{i\in I}D_i$, where $D_i\cong \text{Aut }P_i$. 
There is also the equality $\Omega_1=\prod_{i\in I}(\Omega_1)_i$, in which the group $(\Omega_1)_i$ has a clear meaning.

The analogue of the group $D$ for an arbitrary set $X$ should be defined as follows. For each $i\in I$, let the number $k_i$ be chosen such that every $n_x$ is a multiple of $k_i$, where $x\in X_i$. Then let $P_i=M(k_i,R)$ and let $D$ be the subgroup of all ring (block) automorphisms of the ring $K_i$. Now we set $D=\prod_{i\in I}D_i$, where $D_i\cong\text{Aut}P_i$. There is also an equality $\Omega_1=\prod_{i\in I}(\Omega_1)_i$, in which the group $(\Omega_1)_i$ has an understandable meaning.

The condition $(*)$ in the more general situation under consideration is formulated in the following form $(**)$:
\begin{enumerate}
\item[]
For every $i\in I$ and every $\alpha=(\alpha_x)\in(\Omega_1)_i$, any automorphism $\alpha_x$ is the product in $\text{Aut}R_x$ of an inner automorphism and a ring (block) automorphism.\end{enumerate}

\textbf{Theorem 8.6.} Let the algebra $K=I(X,R)$ satisfy conditions \textbf{(II)} and $(**)$. There are the following relations and isomorphisms in that case.
\textbf{1.}
 
\textbf{a.} $\Omega_1=\text{In }(\text{Aut }L)\cdot D$;

\textbf{b.} $\Omega_1/\text{In }(\text{Aut }L)\cong \prod_{i\in I}\text{Out }P_i$;

\textbf{c.} $C=\text{In}_0(\text{Aut }K)\cdot \Psi\cdot D$.

\textbf{2.} The group $\text{Out}_0K$ contains a normal subgroup of $H$ which is isomorphic to $\Psi/\Psi_0$, and the factor group $\text{Out}_0K/H$ is isomorphic to $\Omega_1/\text{In }(\text{Aut }L)$, i.e. it is isomorphic to group $\prod_{i\in I}\text{Out }P_i$. In addition, $\text{Out }K/\text{Out}_0K\cong \Omega/\Omega_1$.

Next we consider only Corollary 8.5(2) (at the same time, not in the most general form).

\textbf{Corollary 8.7.} Let $K$ be an incidence algebra $I(X,R)$ such that there exists a number $n_s$ dividing every $n_x$. For every $i\in I$, we can take $R_s$ as the ring $P_i$; then condition $(**)$ is satisfied. In particular, we obtain an isomorphism $D\cong\prod_{|I|}\text{Aut }R_s$.

There is one simple construction that allows you to embed an arbitrary algebra $K$ into some incidence algebra $K_0$, satisfying the condition from Corollary 8.7. In this case, the group $\text{Aut}K$ is embedded in the group $\text{Auth}K_0$. In a certain sense, this construction is the opposite of the construction implied in the remark after the proof of Theorem 8.2. In one case, you need to start by partitioning each interval $[x]$ into some parts, and in the other case, on the contrary, extend each such interval. Having obtained the necessary preordered set, in both cases we arrive at the desired incidence algebra $K_0$.

There is one simple construction that allows you to embed an arbitrary algebra $K$ into some incidence algebra $K_0$ satisfying the condition from Corollary 8.7. In this case, the group $\text{Aut}K$ is embedded in the group $\text{Auth}K_0$. In a certain sense, this construction is the opposite of the construction implied in the remark after the proof of Theorem 8.2. In one case, you need to start by splitting each interval $[x]$ into some parts, and in the other case, on the contrary, expand each such interval. Having obtained the necessary preordered set, in both cases we arrive at the desired incidence algebra $K_0$.

Using this construction and assuming that the ring $R$ is commutative, we can find interesting information about the group $\Omega_1/\text{In}(\text{Aut}L)$.

We recall that an Abelian group $B$ is said to be \textsf{bounded} if there is an upper bound for the orders of its elements. A bounded group is a direct sum of cyclic groups.

It is known that for commutative rings $R$, the group of outer automorphisms of the $R$-algebra $M(n,R)$ is a bounded Abelian group such that the orders of its elements divide $n$ \cite{Isa80}.

In two remaining corollaries, the ring $R$ is commutative, and we consider automorphisms of the algebra $R$-$K$ equal to $I(X,R)$.

\textbf{Corollary 8.8.} Let $X_i$, $i\in I$, be all connected components of the set $X$. Then the group $\Omega_1/\text{In}(\text{Aut }L)$ is isomorphic to the product $\prod_{i\in I}B_i$, where every $B_i$ is a bounded Abelian group.

The group $B_i$ corresponds to the algebra $I(X_i,R)$. Corollary 8.9 gives more complete information about $I(X_i,R)$. In Corollary 8.9, the set $X$ is connected.

We denote by $d$ the greatest common divisor of all integers $n_x$, $x\in X$. Based on the construction mentioned above, we can find an embedding $\Omega_1/\text{In}(\text{Aut }L)\to \text{Out }R_s$, for every $s\in X$. Therefore, we obtain the following result.

\textbf{Corollary 8.9.} Let $X$ be a connected set. Then orders of elements of the Abelian group $\Omega_1/\text{In}(\text{Aut }L)$ divide the integer $d$.

\section{Structure of Groups $\text{Aut }K$ and $\text{Out }K$}\label{section9}

Recall that we have an agreement in force, adopted at the end of Section 5.

We systematize some of the studies conducted earlier. In Section 7, the group of outer automorphisms $\text{Out}X$ of the preordered set $X$ was defined. In this case, there is an isomorphism
$$
\text{Out }X\cong G=\{\sigma\in\text{Aut }\overline{X}\,|\,|[x]|=|[\sigma(x)]|\, \text{ for every } x\in X\}.
$$
Then we embedded the group $\text{Out}X$ in $\text{Aut }$ in such a way that the equality
$$
\text{Aut }X=\text{In}(\text{Aut }X)\leftthreetimes\text{Out }X
$$
is true. In Section 7, it was noted that there is always an embedding
$$
t\colon \text{Aut }X\to\text{Aut }K,\; \tau\to\xi_{\tau},\; \text{ where }\;\xi_{\tau}(f)(x,y)= f(\tau(x),\tau(y)), \; f\in K,\; x,y\in X.
$$
We identify $\tau$ and $\text{Aut }X$ with $\xi_{\tau}$ and the image of the embedding $t$, respectively. The relations
$$
\text{Out }X\subseteq \Lambda\; \text{ and }\; \text{Out }X\cap C=\langle 1 \rangle
$$
 are true.
The mapping $ft\colon \text{Out }X\to\Lambda$ is a monomorphism. 
Again, after identification, we get the relations
$$
\text{Out }X\subseteq \Omega\; \text{ and }\; \text{Out }X\cap \Omega_1=\langle 1 \rangle
$$ 
(the designations used here were introduced in Section 5).

On the other hand, the following homomorphisms appeared in Section 6:
$$
 p\colon \text{Aut }K\to \text{Aut }\overline{X}\; \text{ and }\;
r\colon\Omega\to \text{Aut }\overline{X}. 
$$
However, these homomorphisms are of little use. Here's the thing. Let $\begin{pmatrix}\alpha&0\\ \delta&\beta\end{pmatrix}$ be an automorphism and let $\alpha(e_x)=e_y$, as at the beginning of Section 6. 
Then $R_x\cong R_y$, but not necessarily $n_x=n_y$ ($n_x$ is the order of matrices in the ring $R_x$; the symbol $n_x$ is introduced before Theorem 8.2). For this reason, the subgroup $\text{Im}p$ is not always contained in the group $G$, i.e. in $\text{Out}X$ (we take into account the identifications adopted in Section 7).

Now we formulate one condition under which the inclusion of $\text{Im} p\subseteq\text{Out}X$ is true. This clarifies the structure of the group $\text{Aut }K$. The condition is as follows.
$$
\text{If }\;R_x\cong R_y,\; \text{then }\; n_x=n_y\; \text{for any }\;x,y\in X.\eqno (1)
$$
It is clear that this condition holds if $n_x=n_y$ for all $x,y\in X$. For example, if $X$ is a partially ordered set or the following more useful condition holds which is called the \textsf{$(n,m)$-condition}.
$$
\text{For any }\, n,m\in\mathbb{N} \text{ the isomorphism }\, M(n,R)\cong M(m,R) \text{ implies } n=m.\eqno(2)
$$
For example, the rings $R$ from the following list satisfy the $(n,m)$-condition:

commutative rings, local rings, principal left (or right) ideal domains; $(n,m)$-condition appears in \cite{Khr10} and \cite{Vos80}.

Fulfillment of one of Conditions $(1)$ and $(2)$ guarantees the presence of the homomorphism $p\colon\text{Aut}K\to\text{Out}K$.
The composition $p t$ acts identically on $\text{Out }X$. Consequently, $p$ splits and we have the semidirect product $\text{Aut }K=\text{Ker }p\leftthreetimes \text{Out }X$, where 
$$
\text{Ker }p=\{\varphi=\begin{pmatrix}\alpha&0\\ \delta&\beta\end{pmatrix}\in\text{Aut }K\,|\,\alpha R_x=R_x\;\text{ for all }\,x\in X\}.
$$
There are also homomorphisms
$$
ft\colon \text{Out }X\to \Omega,\; r\colon \Omega\to\text{Out }X
$$
and a semidirect decomposition 
$$
\Omega=\text{Ker }r\leftthreetimes \text{Out }X=\Omega_1\leftthreetimes \text{Out }X.
$$
Indeed, the composition $r(ft)$ acts identically on $\text{Out }X$.

Now we can write Theorem 9.1. We note that before Theorem 8.4, the factor group $C/\text{In}_0(\text{Author}K)$ was denoted by $\text{Out}_0K$.

\textbf{Theorem 9.1.} Let algebra $K$ satisfy condition \textbf{(II)} and holds one of Conditions $(1)$ and $(2)$. Then we can formulate the following assertions.

\textbf{a.} There are semidirect decompositions
$$
\Omega=\Omega_1\leftthreetimes \text{Out }X,\; \Lambda=C\leftthreetimes \text{Out }X,
$$
$$
\text{Aut }K=\text{In}_1(\text{Aut }K)\leftthreetimes C\leftthreetimes \text{Out }X.
$$
\textbf{b.} The following assertions hold.

\textbf{1.} There is an isomorphism $\text{Out }K\cong\text{Out }_0K \leftthreetimes \text{Out }X$, the group $\text{Out }_0K$ contains a normal subgroup of $F$ such that $F\cong\Psi/\Psi_0$ and $(\text{Out }_0K)/F\cong \Omega_1/\text{In}(\text{Aut }L)$.

\textbf{2.} If $\Psi=\Psi_0$, then $\text{Out }K\cong \Omega_1/\text{In}(\text{Aut }L)\leftthreetimes \text{Out }X$.

\textbf{3.} If $\Omega_1=\text{In}(\text{Aut }L)$, then $\text{Out }K\cong\Psi/\Psi_0\leftthreetimes \text{Out }X$.

\textbf{Remarks.}

\textbf{1.} It is useful to compare Theorem 9.1 with Theorems 6.4 and 8.4, and also \cite[Theorems 3 and 4]{Khr10}.

\textbf{2.} Section 10 contains conditions which imply the equality $\Psi=\Psi_0$. This section also contains information about the structure of the factor group $\Psi/\Psi_0$.

\textbf{3.} The equality $\Omega_1=\text{In}(\text{Aut }L)$ is  satisfied if every automorphism of the algebra $R_x$ is inner for any $x\in X$ (this follows from Theorem 8.4).

\textbf{Proof.}\\
\textbf{a.} Recently we obtained the relation $\text{Aut }K=\text{Ker }p\leftthreetimes\text{Out }X$. By Theorem 6.2, we have the relation
$\text{Aut }K=\text{In}_1(\text{Aut }K)\leftthreetimes \Lambda$. It follows from the inclusion $\text{In}_1(\text{Aut }K)\subseteq\text{Ker }p$ that
$$
\text{Ker}p=\text{In}_1(\text{Aut }K)\leftthreetimes (\text{Ker }p\cap\Lambda)=\text{In}_1(\text{Aut }K)\leftthreetimes C ,
$$
$$
\Omega=\Omega_1\leftthreetimes \text{Out }X,\; \Lambda=C\leftthreetimes \text{Out }X,
$$
$$
\text{Aut }K=\text{In}_1(\text{Aut }K)\leftthreetimes C\leftthreetimes \text{Out }X.
$$
In particular, the relation $\Lambda= C\leftthreetimes \text{Out }X$ holds.

\textbf{b.} The isomorphism in \textbf{1} follows from Proposition 4.4 and \textbf{a}. 
The remaining statements can be deduced from Theorem 6.4(b), and it is possible, as in Theorem 6.4, to verify their validity directly.
To do this, take $(\text{In}_0(\text{Aut}K)\cdot\Psi)/\text{In}_0(\text{Aut}K)$ as a subgroup of $F$ (in this case, we must take into account Proposition 6.3(1)). Then one can write isomorphisms
$$
\text{Out}_0K/F\cong C/(\text{In}_0(\text{Aut }K)\cdot\Psi)\cong\Omega_1/\text{In }(\text{Aut }L)
$$ 
(see Proposition 6.3(2)).

\textbf{2} and \textbf{3}. The assertions follow from  \textbf{1}.~$\square$

We return to the situation before Theorem 8.4. We assume that all rings $R$ are block matrix rings over the ring $P$, where $P=M(k,R)$ for some $k\ge 1$.

\textbf{Corollary 9.2.} We assume that $X$ is a connected set, Condition $(*)$ from Section 8 holds, and one of Conditions $(1)$, $(2)$ holds. Then the following relations and isomorphisms hold.
$$
\text{Aut }K=\text{In}_1(\text{Aut }K)\leftthreetimes (\text{In}_0(\text{Aut }K)\cdot\Psi\cdot D)\leftthreetimes \text{Out }X, \leqno \textbf{1.}
$$
where $D\cong \text{Aut }P$ or, more briefly, 
$$
\text{Aut }K=(\text{In }(\text{Aut }K)\leftthreetimes \Psi\cdot D)
\leftthreetimes \text{Out }X;
$$
$$
\text{Out }K\cong\Psi/\Psi_0\leftthreetimes\text{Out }P\leftthreetimes \text{Out }X. \leqno \textbf{2.}
$$
\textbf{Proof.}

\textbf{1.} The required relations follow from Theorems 8.4 and 9.1.

\textbf{2.} It follows from Theorem 9.1 that it is sufficient to prove the isomorphism $\text{Out}_0K\cong\Psi/\Psi_0\leftthreetimes D/D_0$ (it is necessary to take into account that $D/D_0\cong\text{Out }P$; see the paragraph before Lemma 8.3). Taking into account Lemma 8.3 and Theorem 8.4, we write down several equalities and isomorphisms
$$
\text{Out}_0K=C/\text{In}_0(\text{Aut }K)=
(\text{In}_0(\text{Aut }K)\cdot\Psi\cdot D)/\text{In}_0(\text{Aut }K)\cong
$$
$$
\cong \Psi D/(\Psi D\cap\text{In}_0(\text{Aut }K))=\Psi D/\Psi_0D_0=
$$
$$
=(\Psi\leftthreetimes D)/(\Psi_0\leftthreetimes D_0)\cong
\Psi/\Psi_0\leftthreetimes D/D_0.\quad\square
$$

Theorem 8.6 transfers Theorem 8.4 to the case of an arbitrary (i.e. not necessarily connected) set $X$. Now we will do the same with Corollary 9.2. Before Theorem 8.6, it is explained how, with such transfers, the groups $D$ and $D_0$ should be defined and what is meant by the group $P_i$.

\textbf{Corollary 9.3.} Let condition $(**)$ from Section 8 and one of conditions $(1)$ and $(2)$ satisfy. Then the following relations and isomorphism hold.
$$
\text{Aut }K=\text{In}_1(\text{Aut }K)\leftthreetimes (\text{In}_0(\text{Aut }K)\cdot\Psi\cdot D)\leftthreetimes \text{Out }X, \leqno \textbf{1.}
$$
or, in short, 
$$
\text{Aut }K=(\text{In }(\text{Aut }K)\cdot\Psi\cdot D)
\leftthreetimes \text{Out }X;
$$
$$
\text{Out }K\cong\Psi/\Psi_0\leftthreetimes\prod_{i\in I}\text{Out }P_i\leftthreetimes \text{Out }X. \leqno \textbf{2.}
$$

\textbf{Proof.}

\textbf{1.} The assertion follows from equalities $\text{Aut }K=\text{In}_1(\text{Aut }K)\leftthreetimes C\leftthreetimes \text{Out }X$, $C=\text{In}_0(\text{Aut }K)\cdot\Psi\cdot D$ and Theorems 9.1, 8.6.

\textbf{2.} The proof of Corollary 9.2(2) remains true in this more general situation. We only need to take into account that we now have the isomorphism $D/D_0\cong\prod_{i\in I}\text{Out}P_i$ (here $I$ is the subscript set for the connected components $X_i$ of the set $X$; see the text before Theorem 8.6).~$\square$

We assume that the ring $R$ is commutative. Then conditions $(1)$ and $(2)$ hold. It is natural to consider automorphisms of the $R$-algebra $I(X,R)$. If we set $k_i=1$ in Section 8 for all $i\in I$ (the integers $k_i$ are defined before Theorem 8.6), then the algebra $I(X,R)$ automatically satisfy Condition $(**)$. For the commutative indecomposable ring $R$, Corollary 9.3  takes the following form.

\textbf{Corollary 9.4.} Let $R$ be a commutative indecomposable ring.

\textbf{a.} The following relations are true. 
$$
\text{Aut }K=\text{In}_1(\text{Aut }K)\leftthreetimes (\text{In}_0(\text{Aut }K)\cdot\Psi)\leftthreetimes \text{Out }X=(\text{In}(\text{Aut }K)\cdot\Psi)\leftthreetimes \text{Out }X, \leqno \textbf{1.}
$$
$$
\text{Out }K\cong\Psi/\Psi_0\leftthreetimes\text{Out }X. \leqno \textbf{2.}
$$
\textbf{b.} If $X$ is a partially ordered set, then there are the following relations:
$$
\text{Aut }K=\text{In}_1(\text{Aut }K)\leftthreetimes\Psi\leftthreetimes \text{Aut }X, \leqno \textbf{1.~\cite{SpiO97}}
$$
$$
\text{Out }K\cong\Psi/\Psi_0\leftthreetimes\text{Aut }X. \leqno \textbf{2.}
$$

\section{Subgroups $\Psi$ and $\Psi_0$}\label{section10}

Below, we consider the ring $L$ and the bimodules $M$, $M_{xy}$ defined in Section 3. It is also useful to remember the decomposition of $K=I(X,R)=L\oplus M$ of the incidence algebra $I(X,R)$.

In this section, we do not actually impose any conditions on the preordered set $X$, the ring $R$, and the algebra $K=I(X,R)$.

The normal subgroup $\Psi$ consists of automorphisms of the form $\begin{pmatrix}1&0\\ 0&\beta\end{pmatrix}$. The subgroup $\Psi_0$ contains all inner automorphisms defined by invertible central elements of the ring $L$. These subgroups appeared in Section 5. The relationship between the subgroups $\Psi$ and $\Psi_0$ is visible from the equalities
$$
\Psi\cap\text{In }(\text{Aut }K)=\Psi\cap\text{In}_0(\text{Aut }K)=\Psi_0
$$
(Proposition 6.3). 
Automorphisms from $\Psi$ are called \textsf{multiplicative}, and automorphisms from $\Psi_0$ are called \textsf{fractional} (there is more text on this topic at the end of the section).

This section provides various information about the groups $\Psi$ and $\Psi_0$. If $X$ is a finite connected set, then the subgroup $\Psi_0$ can be precisely calculated (Corollary 10.6) and it is a direct factor of the group $\Psi$ (Theorem 10.5).

The following fact is standard and essentially known (see also \cite[Proposition 13.2]{KryT21}).

\textbf{Proposition 10.1.} Let $P=M(n,R)$, $Q=M(m,R)$ and $V=M(n\times m,R)$. Endomorphisms (resp., automorphisms) of $P$-$Q$-bimodule $V$ coincide with multiplications by the central elements (resp., invertible central elements) of the ring $R$.

Let $\begin{pmatrix}1&0\\ 0&\beta\end{pmatrix}$ be a multiplicative automorphism. Then $\beta$ is an automorphism of the algebra $M$ and an automorphism of the $L$-$L$-bimodule $M$. Conversely, if a mapping $\beta\colon M\to M$ has indicated properties, then $\begin{pmatrix}1&0\\ 0&\beta\end{pmatrix}$ is a multiplicative automorphism (see the beginning of Section 5 and Proposition 6.1). 

We recall that the subbimodules $M_{xy}$ are defined in Section 3 and they were significantly used in Section 8. Let we have a multiplicative automorphism $\begin{pmatrix}1&0\\ 0&\beta\end{pmatrix}$. Similar to Section 8, we denote the restriction of $\beta$ to $M_{xy}$ by $\beta_{xy}$. Since $\beta$ is an automorphism of the algebra $M$, we have that for any $x,z,y\in X$ with $x<z<y$ and each of $a\in M_{xz}$, $b\in M_{zy}$, the following equality must be fulfilled:
$$
\beta_{xy}(ab)=\beta_{xz}(a)\beta_{zy}(b).\eqno (1)
$$
Due to the fact that $\beta$ is simultaneously an automorphism of an $L$-$L$-bimodule $M$, then $\beta_{xy}$ is an automorphism of the $R_x$-$R_y$-bimodule $M_{xy}$ for any $x,y$ such that $x<y$ (see Proposition 6.1). Next, it follows from Proposition 10.1 that there is element $c_{xy}\in C(U(R))$ with $\beta_{xy}(g)=c_{xy}g$, $g\in M_{xy}$. Taking into account the relation $M_{xz}M_{zy}=M_{xy}$, we use $(1)$ to deduce the relation
$$
c_{xy}=c_{xz}c_{zy},\eqno (2)
$$
where $x<z<y$. Thus, to the automorphism $\psi\in\Psi$, we can match a system of invertible central elements $c_{xy}$ ($x,y\in X$, $x<y$) of the ring $R$. Conversely, every system of elements 
$\{c_{xy}\in C(U(R))\,|\,x<y\}$, satisfying equality $(2)$, provides a multiplicative automorphism $\psi$. Namely, for the element $g=(g_{xy})\in M$, we set $\psi(g)=(c_{xy}g_{xy})$. At the same time, we have to take into account the multiplication rule for elements of $M$ (see Section 3).

We formalize everything just stated above as follows.

\textbf{Proposition 10.2.}\\
\textbf{1.} The is one-to-one correspondence between multiplicative automorphisms and systems of the elements $\{c_{xy}\in C(U(R))\,|\,x<y\}$ such that $c_{xy}=c_{xz}c_{zy}$ for all $x,z,y\in X$ with $x<z<y$.

\textbf{2.} The group $\Psi$ can be embedded in the product $\mathcal{M}$ of copies of the group $C(U(R))$, where $\mathcal{M}=|\{(x,y)\,|\,x<y\}|$.

A more formal point of view on the situation under consideration is as follows. For every pair $(x,y)$ of elements of $X$ with $x<y$, we take a copy $U_{xy}$ of the group $C(U(R))$. Then the group $\Psi$ is isomorphic to a subgroup of the product $\prod_{x<y}U_{xy}$ consisting of vectors $(u_{xy})$, elements of which satisfy equality $u_{xy}=u_{xz}u_{zy}$ if $x<z<y$.

What can be said about fractional automorphisms, i.e.  automorphisms from the subgroup $\Psi_0$? Let a fractional automorphism $\psi=\begin{pmatrix}1&0\\ 0&\beta\end{pmatrix}$ be defined by an invertible central element $v=(v_x)\in L$, where $v_x\in C(U(R))$. For any $x,y$ ($x<y$) and every $g\in M_{xy}$, we have relations
$$
\beta(g)=v^{-1}gv=v_x^{-1}v_yg.
$$
Consequently, the system of the elements $\{c_{xy}\,|\,x<y\}$ from Proposition 10.2 corresponding to the automorphism $\psi$, consists of all elements of the form $v_x^{-1}v_y$, where $x<y$.

We have the following result.

\textbf{Proposition 10.3.} Let $\psi$ be a multiplicative automorphism $\psi$ and let
an element system $c_{xy}$ $(x,y\in X,\; x<y)$, satisfying $(2)$, correspond to $\psi$. The following assertions are equivalent:

\textbf{1)} $\psi$ is a fractional automorphism; 

\textbf{2)} $\psi$ is an inner automorphism;

\textbf{3)} there is an element $v=(v_x)\in L$, where $v_x\in C(U(R))$ such that the relation $c_{xy}=v_x^{-1}v_y$ holds for all $x,y$ $(x<y)$.

\textbf{In the remaining part of the section, the set $X$ is finite and connected. For technical convenience, we also assume that $X$ is partially ordered.}

Let we have an arbitrary set $\{c_{xy}\,|\,x,y\in X,\; x<y\}$ of central invertible elements of the ring $R$. Based on this set, we assign weithts to edges and semipaths of this graph $X$ (see \cite{BruFS15}). If $(x,y)$ is an edge, then we set $w(x,y)=c_{xy}$ for $x<y$ and $w(x,y)=c_{yx}^{-1}$ for $x>y$. Next for semipath $P=z_1z_2\ldots z_{k+1}$ in $X$, we set
$$
w(P)=\prod_{i=1}^kw(z_i,z_{i+1}).
$$
If $\psi$ is some multiplicative automorphism, then we obtain weithts of edges and semipaths of the graph $X$. They are induced by the system $\{c_{xy}\,|\, x<y\}$ corresponding to the automorphism $\psi$ in the sense of Proposition 10.2. Let
$\psi$ be a fractional automorphism (i.e. an inner multiplicative automorphism; see Proposition 10.3) which is defined by the  invertible element $v=(v_x)$. It is easy to verify the following. For any elements $x,y\in X$, the weight of any semipath from $x$ to $y$ is equal to $v_x^{-1}v_y$ (therefore, it does not depend on  the specific semipath).

We fix some spanning tree $T$ of the graph $X$. For every edge $(x,y)$ in $T$ with $x<y$, we choose an element $c_{xy}\in C(U(R))$. Similar to what was done above, we use the elements $c_{xy}$ to assign the weights $w$ to edges and semipaths in $T$. Let $(x,y)$ be an edge in $X$, not contained in $T$, and let $x<y$. We set $c_{xy}=w(P)$, where $P$ is the (unique) semipath in $T$ from $x$ to $y$. As a result, we obtain a system of the elements $\{c_{xy}\in C(U(R))\,|\, x,y\in X,\, x<y\}$ and also weithts of edges and semipaths in $X$.

\textbf{Lemma 10.4.} The constructed system of elements $\{c_{xy}\,|\, x<y\}$ satisfies relations $(2)$. Consequently, it induces a multiplicative automorphism $\psi$ in the sense of Proposition 10.2.

\textbf{Proof.} Let $x,z,y$ be elements with $x<z<y$. We denote by $P_{xz}$, $P_{zy}$ and $P_{xy}$ semipaths in $T$ from $x$ to $z$, from $z$ to $y$ and from $x$ to $y$, respectively. Since every semipath in $T$ is unique, we have the relation $P_{xy}=P_{xz}P_{zy}$. Therefore, $w(P_{xy})=w(P_{xz})\cdot w(P_{zy})$ and, therefore, $c_{xy}=c_{xz}c_{zy}$. By Proposition 10.2, a quite definite automorphism $\psi\in \Psi$ corresponds to the system $\{c_{xy}\,|\, x<y\}$.~$\square$

Since $X$ is partially ordered set, the bimodules $M_{xy}$ are copies of the ring $R$ in our situation. We also can define analogues of matrix units. Let $x,y\in X$ with $x<y$. We denote by $e_{xy}$ a function $X\times X\to R$ such that $e_{xy}(s,t)=1$ provided $s=x$, $t=y$ and $e_{xy}(s,t)=0$ for all other pairs $(s,t)$. The functions $e_{xy}$ satisfy the following property: if $x<z<y$, then $e_{xz}e_{zy}=e_{xy}$.

We pay attention to that we have chosen the spanning tree $T$ of the graph $X$. We denote by $\Psi_1$ the subgroup in $\Psi$ consisting of automorphisms $\psi$ such that $\psi(e_{xy})=e_{xy}$ for all edges $(x,y)\in T$ with $x<y$. Equivalently, $c_{xy}=1$ for indicated edges $(x,y)$. Here $\{c_{xy}\in C(U(R))\,|\, x<y\}$ is the system from Proposition 10.2. In addition, it follows from this proposition that the group $\Psi$ is commutative.

\textbf{Theorem 10.5.} There is a group direct decomposition $\Psi=\Psi_1\times \Psi_0$.

\textbf{Proof.} Let $\psi\in\Psi$ and let $\{c_{xy}\in C(U(R))\,|\, x<y\}$ be the corresponding system of the elements from Proposition 10.2. For every vertices from $X$, we define some central invertible element.

In the tree $T$, we choose some root $x_0$. Next, let $y$ be some vertex in $T$ which is adjacent to $x_0$. If $x_0<y$, then we set $v_{x_0}=1$ and $v_y=v_{x_0}c_{x_0y}$. If  $x_0>y$, then we set $v_y=v_{x_0}c_{yx_0}^{-1}$. Next we take some vertex $z\in T$ adjacent to $y$ and act similarly. Namely, let $v_z=v_yc_{yz}$ for $y<z$ and let $v_z=v_{y}c_{zy}^{-1}$ for $y>z$. We act in this way until we get to some hanging vertex in $T$. We obtain semipath in $T$ from $x_0$ to this hanging vertex.

We choose some root $x_0$ in the tree $T$. Next, let $y$ be some vertex adjacent to $x_0$ in $T$. If $x_0<y$, then we assume $v_{x_0}=1$ and $v_y=v_{x_0}c_{x_0y}$. If $x_0>y$, then we assume $v_y=v_{x_0}c_{yx_0}^{-1}$. Then we take some vertex $z$ adjacent to $y$ in $T$ and do the same. Namely, for $y<z$, let $v_z=v_yc_{yz}$ and for $y>z$, let $v_z=v_{y}c_{zy}^{-1}$. We act in this way until we reach some hanging vertex in $T$. We get a semipath in $T$ from $x_0$ to this hanging vertex. 
If there are vertices in $T$ that do not belong to this semipath, then let $s$ be the first vertex in the semipath with degree exceeding two.

Starting from the vertex $s$, we similarly form another semipath, and so on. As a result, each vertex $x$ in $T$ is mapped to an element $v_x\in C(U(R))$.

We take an invertible central element $v=(v_x)$ in $L$. Let $\psi_0$ be the inner automorphism induced by $v$. Then $\psi_0\in \Psi_0$, since $v\in C(U(L))$. We set $\psi_1=\psi\psi_0^{-1}$. If $(x,y)\in T$ and $x<y$, then $\psi(e_{xy})= c_{xy}e_{xy}$ by Proposition 10.1 and $\psi_0(e_{xy})=v_x^{-1}v_ye_{xy}$. By construction of the elements $v_x$, the relation $c_{xy}=v_x^{-1}v_y$ holds. Consequently, $\psi$ and $\psi_0$ act the same on elements $e_{xy}$ such that $(x,y)\in T$. Therefore, the inclusion $\psi_1\in \Psi_1$ is true. So, $\psi=\psi_1\psi_0$, where $\psi_1\in \Psi_1$, $\psi_0\in \Psi_0$.

It remains to verify that $\Psi_1\cap \Psi_0=\langle 1 \rangle$. Let $\xi\in \Psi_1\cap \Psi_0$. Then $\xi(e_{xy})=e_{xy}$ for every edge $(x,y)\in T$, where $x<y$. At the same time, the automorphism $\xi$ is defined as an inner automorphism by some central invertible element $v=(v_x)$ of the ring $L$. Therefore, we also have relation $\xi(e_{xy})=v_x^{-1}v_ye_{xy}$. Consequently, $v_x^{-1}v_y=1$ and $v_x=v_y$. Now let edge $(x,y)$ be not contained in $T$. There is a semipath in $T$ from $x$ to $y$. 
Using this semipath, we deduce the equality $v_x=v_y$ from what has just been proved
Therefore, we obtain the relation $\xi=1$. Thus, $\Psi_1\cap \Psi_0=\langle 1 \rangle$ which completes the proof.~$\square$

We denote by $m(X)$ and $\lambda(X)$ the number of edges of the graph $X$ and the cyclomatic number of $X$, respectively.

\textbf{Corollary 10.6.} There is a group isomorphism
$$
\Psi_0\cong \prod_{m(X)-\lambda(X)}C(U(R)).
$$

\textbf{Proof.} We denote by $G$ the product in the right part of the isomorphism. We define mapping $g\colon \Psi\to G$ as follows. Let $\psi\in \Psi$ and let $\{c_{xy}\,|\, x<y\}$ be the system of the elements from Proposition 10.2 corresponding to $\psi$. The mapping $g$ maps the automorphism $\psi$ to the vector $(c_{xy})$, where $(x,y)\in T$, $x<y$. It is clear that $g$ is a group homomorphism. It follows from Lemma 10.4 that $g$ is surjective. It is also clear that $\text{Ker }g=\Psi_1$. We have an isomorphism $\Psi/\Psi_1\cong G$. On the other hand, By Theorem 10.5 we have isomorphism $\Psi/\Psi_1\cong \Psi_0$.~$\square$

\textbf{Corollary 10.7.} The equality $\Psi=\Psi_0$ is true if and only if $\psi=1$ for any $\psi\in \Psi$ such that $\psi(e_{xy})=e_{xy}$ for any edge $(x,y)\in T$.

\textbf{Corollary 10.8 \cite{BruFS15}.} The multiplicative automorphism $\psi$ is inner if and only if for every elements $x,y\in X$, weithts of any two semipaths from $x$ to $y$ coincide.

\textbf{Proof.} If the multiplicative inner automorphism $\psi$ is defined by an element $v=(v_x)$, then the weight of any semipath from $x$ to $y$ is equal to $v_x^{-1}v_y$. This has already been noted above.

We assume that holds condition for weithts of semipaths for automorphism $\psi$. As in Theorem 10.5, we write down $\psi=\psi_1\psi_0$. If $(x,y)$ is an edge in $T$, then $\psi(e_{xy})=\psi_0(e_{xy})$. Now let $(x,y)$ be an edge in $X$ not contained in $T$. Then $\psi(e_{xy})=c_{xy}e_{xy}$ and $w(x,y)=c_{xy}$. If $P$ is a semipath in $T$ from $x$ to $y$, then $w(P)=v_x^{-1}v_y$, where the invertible element $v=(v_x)$ defines an inner automorphism $\psi_0$. It follows from the assumption that $c_{xy}=v_x^{-1}v_y$. Thus, $\psi$ acts on all elements $e_{xy}$ as $\psi_0$ and, therefore $\psi=\psi_0$.~$\square$

Let $C$ be some simple cycle of the graph $X$ and let $\psi$ be some multiplicative automorphism. We define the weight $w(C)$ of the cycle $C$ with respect to $\psi$ as the weight of one of two semipaths from any vertex of the cycle to itself. When choosing a different direction of the cycle pass, the weight will be equal to $w(C)^{-1}$. Thus, the weight of the cycle $C$ this one from elements $w(C)$ and $w(C)^{-1}$.

Let $T$ be a spanning tree of the graph $X$ and let $C_1,\ldots,C_k$ be the fundamental system of cycles, associated with $T$, where $k=\lambda(X)$.

\textbf{Corollary 10.9.} The multiplicative automorphism $\psi$ is inner if and only if $w(C_1)=\ldots=w(C_k)=1$. 

\textbf{Proof.}Necessity of the assumption follows from the first paragraph of the proof of Corollary 10.8.

Now we assume that $w(C_1)=\ldots=w(C_k)=1$. We represent the multiplicative automorphism $\psi$ in the form $\psi=\psi_1\psi_0$ (see the proof of Corollary 10.8). Let an edge $(x,y)$, $x<y$, is not contained in $T$. Then $(x,y)\in C_i$ for some $i$. If $P$ is a semipath in $C_i$ from $x$ to $y$, then there are relations
$$
1=w(C_i)=w(P)w(x,y)^{-1},\quad w(x,y)=w(P).
$$
By repeating the argument from the proof of Corollary 10.8, we obtain that $c_{xy}=v_x^{-1}v_y$ and $\psi=\psi_0$.~$\square$

In the book \cite{SpiO97}, multiplicative and fractional automorphisms are defined with the use of some special functions in $I(X,R)$, where $R$ is a commutative ring. These definitions can be extended to the case of an arbitrary ring $R$.

A function $m\in I(X,R)$ is said to be \textsf{multiplicative} if for every elements $x,z,y\in X$ such that $x\le z\le y$, we have 
$$
m(x,y)\in C(U(R)),\quad m(x,y)=m(x,z)\cdot m(z,y).
$$
Let $m$ be some multiplicative function in $I(X,R)$. Then the mapping
$$
\psi_m\colon K\to K,\; \psi_m(f)=m*f,\; f\in K,
$$
is an automorphism of the algebra $K$ (here $*$ is the product of Hadamard functions, i.e. the point-wise product). It is clear that $\psi_m\in \Psi$. Conversely, every automorphism from $\Psi$ can be obtained in this way, with the use of some multiplicative function.

Let $q\colon X\to C(U(R))$ be some mapping. Then the function $m_q\in K$ such that $m_q(x,y)=q^{-1}(x)\cdot q(y)$ for $x\le y$ and $m_q(x,y)=0$ in the remaining cases, is multiplicative. In addition, the corresponding automorphism $\psi_{m_q}$ is fractional. All fractional automorphisms are obtained in this way.

\section{Reduced Incidence Algebras}\label{section11}

The approaches developed in the previous sections and the results obtained can be transferred to certain subalgebras of the algebra $I(X,R)$. We are talking about reduced incidence algebras. They are given quite a lot of attention in the book \cite{SpiO97}.

As in \cite{SpiO97}, we further assume that $X$ is a partially ordered set, $R$ is a commutative ring, and consider the ring $I(X,R)$ as an $R$-algebra. Although much can be done in the more general situation that is discussed in this chapter. We not accept any additional conditions, such as those written in the agreement at the end of Section 5.

Here is some material from \cite{SpiO97}. A reduced incidence algebra is defined using some equivalence relation on the set of all intervals of a partially ordered set $X$. Any such equivalence relation is denoted by $\sim$ or by the letter $E$.

\textbf{Definition 11.1.} Let $\sim$ (or $E$) be such an equivalence relation on the set of all intervals of the set $X$ that, whenever $[x,y]\sim[s,t]$, there exists a bijection $\varepsilon\colon [x,y]\to [s,t]$ such that the relations
$$
[[x,z]\sim [s,\varepsilon(z)]\; \text{ and }\; [z,y]\sim [\varepsilon(z),t]
$$
are fulfilled for each $z\in [x,y]$.
In this case, it is said that $\sim$ is a \textsf{order-compatible relation}.

Let $\sim$ be a relation compatible with the order. We denote by $I(X_E,R)$ the set of all functions $f$ of $I(X,R)$ $f(x,y)=f(s,t)$ whenever $[x,y]\sim[s,t]$. In other words, $f$ is a constant on the equivalence classes $\sim$.

\textbf{Proposition 11.2.} The following assertions are true.

\textbf{1.} $I(X_E,R)$ is a subalgebra of the algebra $I(X,R)$.

\textbf{2.} We have the relation
$$
U(I(X_E,R))=I(X_E,R)\cap U(I(X,R)).
$$

The subalgebra $I(X_E,R)$ is called a \textsf{reduced incidence algebra}. Let $\sim$ be an equivalence relation such that $[x,y]\sim[s,t]$ exactly when the intervals $[x,y]$ and $[s,t]$ are isomorphic as partially ordered sets. It is clear that this relation $\sim$ is compatible with the order. The corresponding algebra $I(X_E,R)$ is called a \textsf{standard reduced incidence algebra}. Such algebras represent the most natural objects among reduced incidence algebras. There are quite a lot of them. So, if we take as $X$ the set of all non-negative integers with the usual order, then the corresponding standard reduced incidence algebra is isomorphic to the algebra of formal power series (\cite[Example 1.3.10]{SpiO97}).

Let $I(X_E,R)$ be a reduced incidence algebra defined by some equivalence relation $\sim$ (or $E$) compatible with the order. We also denote it by the symbol $K_E$. The algebra $K$ itself, which is equal to $I(X,R)$, is a splitting extension of $K=L\oplus M$ (Section 3). Our study of the group  $\text{Aut}K$ is based on this circumstance. There is also a similar extension for the algebra $K_E$: $K_E=L_E\oplus M_E$ (see details below). Moreover, the subring $L_E$ and the ideal $M_E$ are also products of copies of the ring $R$. In order to achieve a more complete analogy with the extension $K=L\oplus M$, we have to transform the subring $L_E$ and the ideal $M_E$. 

We set $L_E=\{f\in L\,|\,f(x,x)=f(s,s),\;\text{if }\,[x,x]\sim [s,s]\}$ and $M_E=\{f\in M\,|\,f(x,y)=f(s,t),\;\text{if }\,[x,y]\sim [s,t],\;\text{where }\, x<y,\,s<t\}$. We have the relations $L_E=L\cap K_E$ and $M_E=M\cap K_E$. If $f\in K_E$ and $f=g+h$, where $g\in L$, $h\in M$, then it is true that $g\in L_E$, $h\in M_E$. We can conclude that the equality $K_E=L_E\oplus M_E$ is true, where $L_E$ is a subring and $M_E$ is an ideal in $K_E$.
So, we have presented the algebra $K_E$ in the form of a splitting extension.

Following \cite{SpiO97}, we call \textsf{types} the equivalence classes of the set of intervals with respect to the relation $\sim$. We denote by $T$ the set of all such types. We can write $T=T_0\bigcup T_1$, where $T_0$ is the set of types consisting of intervals of the form $[x,x]$ and $T_1$ is the set of all types of the set of intervals of the form $[x,y]$ with $x<y$. Then we obtain the relation
$$
L_E=\prod_{t\in T_0}R_t,\qquad M_E=\prod_{t\in T_1}M_t, \eqno(1)
$$
where all $R_t$ and $M_t$ are equal to $R$.

So, if $t\in T$, then any function $f\in K_E$ is a constant on all intervals of $t$. If so, then we can consider $f$ as a function defined on the set $T$. To describe the composition of such functions or, equivalently, the multiplication of elements of $K_E$, a suitable set of numbers is introduced.

\textbf{Definition 11.3.} Let $r,s,t$ be three types and let $[x,y]$ be an interval of type $t$. The symbol $\left[\begin{smallmatrix}&t&\\ r&&s\end{smallmatrix}\right]$ denotes the number of elements of the (finite) set $\{z\in[x,y]\,|\,[x,z]\in r,\;[z,y]\in s\}$. This number is called the \textsf{incidence coefficient} of the algebra $K_E$.

We take two arbitrary elements $f,g$ of the algebra $K_E$. With the use of decomposition $K_E=L_E\oplus M_E$ and relations $(1)$, we represent $f$ and $g$ in the form $f=(f_t)$, $g=(g_t)$. Let $f\cdot g=h=(h_t)$. Then there is relation (see \cite[Proposition 3.1.4]{SpiO97}):
$$
h_t=\sum_{r,s\in T}\left[\begin{smallmatrix}&t&\\ r&&s\end{smallmatrix}\right]f_r\cdot g_s. \eqno(2)
$$

Based on the splitting extension $K_E=L_E\oplus M_E$, we can study the group $\text{Aut}K_E$. But there is one significant difference from the situation with the extension $K=L\oplus M$. In the case of the algebra $K$, the idempotents $e_x$ are a useful technical tool (defined in Section 3). It is important that $e_xMe_x=0$.

We denote by $e_t$ the identity element of the ring $R_t$ for each $t\in T_0$. It may happen that $e_tM_Ee_t\ne 0$. To fix this, we move on to another splitting extension.

We introduce one notation. Let $p,q\in T_0$. We set
$$
p\times q=\{[x,y]\,|\,[x,x]\in p,\, [y,y]\in q,\, x<y\}.
$$
For every $p\in T_0$, let $S_p$ be the ring $e_pK_Ee_p$. We have a splitting extension
$$
S_p=R_p\oplus N_p,\, \text{ where }\, N_p=e_pM_Ee_p=\prod_{t\subseteq p\times p}M_t.\eqno(3)
$$
Next we set
$$
P_E=\prod_{p\in T_0}S_p.\eqno(4)
$$

For an arbitrary type $t$ in $T_1$, there are types $p,q\in T_0$ that have the following property:
\begin{enumerate}
\item[]
For any interval $[x,y]\in t$, inclusions $[x,x]\in p$ and $[y,y]\in q$ hold, the types $p$ and $q$ are uniquely determined, and we can write $t\subseteq p\times q$.
\end{enumerate}

For the types $p,q\in T_0$, $p\ne q$, we assume
$$
N_{pq}=e_pM_Ee_q,\quad N_E=\prod_{p,q\in T_0}N_{pq}.\eqno(5)
$$
Here $N_{pq}=\prod_{t\subseteq p\times q}M_t$. In addition to the decomposition $K_E=L_E\oplus M_E$, we obtain another decomposition:
$$
K_E=P_E\oplus N_E.\eqno(6)
$$
For this new decomposition, we have $e_pN_Ee_p=0$ for each $p\in T_0$.

There is a formula $(2)$ for multiplying elements of the algebra $K_E$ written with respect to the products specified in $(1)$. How do these elements multiply if they are represented relative to the products in $(4)$ and $(5)$? We deal with this. This is useful if the decomposition of $(6)$ is required in the process of studying the automorphisms of the algebra $K_E$.

We write down the main equalities we have:
$$
K_E=P_E\oplus N_E,\,\text{where }\, P_E=\prod_{p\in T_0}S_p,\;N_E=\prod_{p,q\in T_0}N_{pq},\; N_{pq}=\prod_{t\subseteq p\times q}M_t.
$$
It is not difficult to clarify how the elements of $P_E$ are multiplied among themselves and by the elements of $N_E$. We focus on multiplying elements of $N_E$ with each other.

We take two arbitrary elements $f,g$ of $N_E$. Let $f\cdot g=h$. We have
$$
f=(f_t)_{t\in T_1}=(f_{pq})_{p,q\in T_0},\, \text{ where }\,
f_{pq}=(f_t)_{t\subseteq p\times q},
$$
$$
g=(g_t)_{t\in T_1}=(g_{pq})_{p,q\in T_0},\, \text{ where }\,
g_{pq}=(g_t)_{t\subseteq p\times q}.
$$
We also have similar relations for $h$:
$$
h=(h_t)_{t\in T_1}=(h_{pq})_{p,q\in T_0},\, \text{ where }\,
h_{pq}=(h_t)_{t\subseteq p\times q}.
$$

We assert that for any distinct $p,q\in T_0$, the relation
$$
h_{pq}=(f_{pr}\cdot g_{rq})_{r\in T_0} \eqno(7)
$$
is true. Informally speaking, when multiplying vectors $(f_t)$ and $(g_t)$ according to the relation $(2)$, the vectors $(f_{pq})$ and $(g_{pq})$ are multiplied according to the relation $(7)$.

First, we note that the various products of $f_{pr}\cdot g_{rq}$ from $(7)$ do not have common non-zero components of $h_t$.

Then, given $(2)$, we can write the equalities
$$
h_{pq}=(h_t)_{t\subseteq p\times q}=\left(\sum_{u,v\in T_1}\left[\begin{smallmatrix}&t&\\ u&&v\end{smallmatrix}\right]f_u\cdot g_v\right)_{t\subseteq p\times q}.\eqno(8)
$$

It is clear that any component of $h_t$ in the product of $f_{pr}\cdot g_{rq}$ is a component for $h_{pq}$ as well. On the other hand, the element $f_u$ is contained in some bimodule $N_{pr}$, and $g_v$ is contained in some bimodule $N_{rq}$. Therefore, on the contrary, every component in parentheses on the right in $(8)$ is also present in the sum of a finite number of products of $f_{pr}\cdot g_{rq}$ in $(7)$. So, the equality $(7)$ is indeed true. This is a new formula for multiplying elements from $N_E$. It is useful in the process of studying the group $\text{Aut}K_E$, if we proceed from the decomposition of $K_E=P_E\oplus N_E$. At the same time, we have to know the structure of the group $\text{Aut}S_p$. To a certain extent, the structure is understandable.

It may happen that the algebra $K_E$ consists of a single ring of the form $S_p$. All standard reduced algebras belong to such algebras. In general, automorphisms of algebras $K_E$ require a separate systematic study.

\textbf{Remarks.} A variety of information about papers in which automorphisms (as well as isomorphisms and derivations) are studied is contained in the book \cite{SpiO97}, which was mentioned in Section 1.

Among papers published after the publication of this book, we highlight the following. The paper \cite{DroK07} is devoted to the calculation of the group of outer automorphisms of the algebra $I(X,R)$ for a finite set $X$. At the same time, an approach based on the cohomology interpretation of the automorphism group is developed.

The paper \cite{Khr10} concerns automorphisms of finitary incidence algebras. These are somewhat more general objects than just incident algebras (\cite{KhrN09}, \cite{Khr10b}). A finite incidence algebra also has a decomposition of the form $K=L\oplus M$. To it, apparently, it is possible to apply the approaches and techniques used in this work.

In \cite{BruFS15} we consider multiplicative automorphisms of the algebra $I(X,R)$ for a finite partially ordered set $X$ and a field $R$.

The review of \cite{BruL11} is devoted to automorphisms, antiautomorphisms and involutions (i.e. antiautomorphisms of order 2) of incidence algebras. There are also new results in it. In general, the involutions of incidence rings and triangular matrix rings have been actively studied recently; for example, \cite{BruFS11}, \cite{BruL11}, \cite{BruFS12}, \cite{BruFS14}, \cite{ForP21b}, \cite{ForP22}, \cite{DiKL06}.

In earlier articles \cite{AnhW11}, \cite{Kez90}, automorphisms of (ordinary and generalized) triangular matrix rings are studied.

\section*{\textbf{CHAPTER 3. ISOMORPHISMS\\ OF INCIDENCE ALGEBRAS}}\label{chapter3}
\addtocontents{toc}{\textbf{CHAPTER 3. ISOMORPHISMS OF\\ \mbox{\,} $\quad$ INCIDENCE ALGEBRAS}\par}

\section[Some Definitions and Auxiliary Results]{Some Definitions and\\ Auxiliary Results}\label{section12}

The agreement from Section 3 that all rings are algebras over some commutative ring $T$ remains valid.

Let $Y$ and $X$ be arbitrary preordered sets and let $\overline{Y}$ and $\overline{X}$ be the corresponding partially ordered sets (see Section 2). Next, let $R$ be a ring and let $K'=I(Y,R)$, $K=I(X,R)$ be two incidence algebras. We formulate the following questions.

\textbf{(a)} When does an isomorphism $K'\cong K$ of algebras imply an isomorphism of partially ordered sets $\overline{Y}\cong\overline{X}$?

\textbf{(b)} When does an isomorphism of algebras $K'\cong K$ imply an isomorphism of preordered sets $Y\cong X$?

The questions \textbf{(a)} and \textbf{(b)} are some variations of the well-known isomorphism problem for incidence rings.

The third question is related to the description of isomorphisms between the algebras $K'$ and $K$.

\textbf{(c)} Under what conditions can any isomorphism $\varphi\colon K'\to K$ be <<diagonalized>>? This means that there are an inner automorphism $\nu$ of the algebra $K$ and a diagonal isomorphism $\gamma\colon K'\to K$ such that $\varphi=\nu\gamma$.

It should be noted that the concepts, statements and proofs from this section are technically and ideologically close to the concepts, statements and proofs of Section 5. For this reason, the presentation here will be brief and almost all the proofs will be omitted.

Similar to Section 3, we write down the decomposition $K=L\oplus M$ and a similar decomposition for the algebra $K'$: $K'=L_1\oplus M_1$, in which the subring $L_1$ and the ideal $M_1$ have a clear meaning.

Having an arbitrary algebra homomorphism $\varphi\colon K'\to K$, it is possible to make a $2\times 2$ matrix $\begin{pmatrix}\alpha&\gamma\\ \delta&\beta\end{pmatrix}$ in a standard way, where 
$$
\alpha\colon L_1\to L,\; \beta\colon M_1\to M,\; \gamma\colon M_1\to L,\; \delta\colon L_1\to M
$$
are some $T$-module homomorphisms.

We only deal with the <<triangular>> case, i.e. when $\gamma=0$. We consider only the isomorphisms of $\varphi$. Also, we do not distinguish between the isomorphism $\varphi$ and the corresponding $2\times 2$ matrix. Sometimes, we write a <<triangular isomorphism>>
$\varphi$, where 
$\varphi=\begin{pmatrix}\alpha&0\\ \delta&\beta\end{pmatrix}$, and a <<diagonal isomorphism>> $\varphi$, where 
$\varphi=\begin{pmatrix}\alpha&0\\ 0&\beta\end{pmatrix}$.

If $\begin{pmatrix}\alpha&0\\ \delta&\beta\end{pmatrix}$ is a triangular algebra isomorphism $K'\to K$, then $\alpha$ and $\beta$ are algebra isomorphisms $L_1\to L$ and $M_1\to M$, respectively (we assume $M_1$ and $M$ are non-unital algebras).

At the beginning of Section 3, the idempotents $e_x$ of the algebra $K$ are defined. We denote by $f_y$ similar idempotents of the algebra $K'$. 

Conditions \textbf{(I)} and \textbf{(II)} below can be considered as a transfer of Conditions \textbf{(I)} and \textbf{(II)} from Section 5 to the situation of two algebras $K'$ and $K$.

\textbf{(I)} Any isomorphism of $K'\to K$ is triangular.

\textbf{(II)} For any isomorphism $\varphi\colon K'\to K$ and every $x\in X$, the inclusion of $\varphi^{-1}(e_x)\in f_y+M_1$ is true for some $y\in Y$.

For both conditions, their symmetric counterparts are satisfied. Specifically, if \textbf{(I)} is true, then any isomorphism $K\to K'$ is also triangular. Before passing to \textbf{(II)}, we reveal the relationship between these conditions.

In fact, the proof of the following lemma repeats the proof of Lemma 5.6 and we omit the proof.

\textbf{Lemma 12.1.} For algebras $K'$ and $K$, the following assertions are true.

\textbf{1.} Condition \textbf{(II)} implies Condition \textbf{(I)}.

\textbf{2.} For an indecomposable ring $R$, Conditions \textbf{(II)} and \textbf{(I)} are equivalent.

\textbf{Lemma 12.1.} The following statements hold for the algebras $K'$ and $K$.

\textbf{1.} From Condition \textbf{(II)} Condition \textbf{(I)} follows.

\textbf{2.} For an indecomposable ring $R$, Conditions \textbf{(II)} and \textbf{(I)} are equivalent.

We recall the equality $L=\prod_{x\in X}R_x$ from Proposition 3.1. Now we return to the question of the symmetry of Condition \textbf{(II)}. We formalize the answer to it in the form of the following lemma.

\textbf{Lemma 12.2.} If Condition \textbf{(II)} is satisfied, then it also holds for the algebras $K$ and $K'$.

\textbf{Proof.} We have to check that for any isomorphism $\psi\colon K\to K'$ and any $y\in Y$, there is an inclusion of $\psi^{-1}(f_y)\in e_x+M$ for some $x\in X$.

Since Condition \textbf{(II)} is satisfied, it follows from Lemma 12.1 that all isomorphisms between $K'$ and $K$ in any direction are triangular. Therefore, we have
$$
\psi=\begin{pmatrix}\alpha&0\\ \delta&\beta\end{pmatrix}\; \text{and} \;
\psi^{-1}=\begin{pmatrix}\alpha^{-1}&0\\ \delta'&\beta^{-1}\end{pmatrix}.
$$
From application of Condition \textbf{(II)} to $\psi^{-1}$, it follows that for any $x\in X$, there exists $y\in Y$ such that $\psi(e_x)\in f_y+M_1$ for some $y\in Y$. Therefore, $\alpha(e_x)=f_y$, and $\psi^{-1}(f_y)=e_x+d$, where $d\in M$.

If $z\in X$, $z\ne x$, and $\alpha(e_z)=f_t$, then the equality $e_xe_z=0$ implies the equality $f_yf_t=0$, from which $y\ne t$. We conclude that the mapping $x\to y$, where $\alpha(e_x)=y$, sets the bijection of the set $X$ onto some subset $Y'$ of $Y$.

It is enough to check that $Y'=Y$. Otherwise, let $y\in Y$ and $y\notin Y'$. Then, let $\alpha(a)=f_y$, where $a\in L$. We write
$$
a=(a_x)=a_s+c,\;\text{where }\; a_s\ne 0,\; c\in \prod_{t\ne s}R_t.
$$
We get the relations $\psi(e_s)\in f_k+M_1$ and $\alpha(e_s)=f_k$ for some $k\in Y'$. Then we have $\alpha(e_sa)=f_kf_y=0$. Hence $e_sa=a_s=0$, which contradicts the choice of the $a_s$ element. We conclude that $Y'=Y$.~$\square$

In what follows, we use conditions $(1)$--$(5)$ given after Definition 5.1.

The following statement is some analogue of Lemma 5.7. Its proof also simulates the proof of this lemma.

\textbf{Lemma 12.3.} If all rings $R_x$ have Property $(5)$ from Section 5, then Condition \textbf{(I)} holds for the algebras $K'$ and $K$.

After proving Lemma 5.7, it is noted that the ring $M(n,R)$ satisfies Condition $(5)$ if the ring $R$ belongs to one of the following classes of rings:

\textbf{a)} local rings;

\textbf{b)} areas of the main left (or right) ideals;

\textbf{c)} commutative dedekind domains.

Taking into account Lemmas 12.1 and 12.3, we obtain the following. If a ring $R$ belongs to one of the classes of rings \textbf{a)}--\textbf{c)}, then the algebras $K'$ and $K$ satisfy the conditions \textbf{(I)} and \textbf{(II)}.

The statement below transfers Lemma 5.9 to the situation of two algebras $K'$ and $K$. At the same time, the proof of Lemma 5.9, after minor corrections, is suitable in this more general situation.

\textbf{Lemma 12.4.} If the quotient ring $R/J(R)$ is indecomposable, then Condition \textbf{(II)} holds for the algebras $K'$ and $K$.

Now we can reproduce the text at the end of Section 5. We do it very briefly.

The most common and convenient condition for answering  questions \textbf{(a)}--\textbf{(c)} can specify Condition \textbf{(II)}. The algebras $K'$ and $K$ satisfy Condition \textbf{(II)} if the ring $R$ belongs to one of the following classes of rings:

\textbf{a)} local rings;

\textbf{b)} principal left (or right) ideal domains;

\textbf{c)} commutative Dedekind domains;

\textbf{d)} $Y$ and $X$ are partially ordered sets and $R$ has no non-zero idempotents except 1;\\

\textbf{e)} the quotient ring $R/J(R)$ is indecomposable.

In the next section, we assume that the algebras $K'$ and $K$ satisfy Condition \textbf{(II)}. In particular, this is the case if $R$ is some ring from the list \textbf{a)}--\textbf{e)}.

\section{Questions $(a)$--$(c)$}\label{section13}

The symbols $Y$, $X$, $\overline{Y}$, $\overline{X}$, $K'$ and $K$ continue to have the meaning given to them in the previous section. There is also a recent agreement on the algebras $K'$ and $K$. A serious role is played by $(n,m)$-condition $(2)$ and Condition $(1)$, formulated before Theorem 9.1.

The arguments from the proof of the first of the two main results of this section are similar to the arguments at the beginning of Section 6. Nevertheless, for completeness and independence of presentation, we give detailed proofs.

\textbf{Theorem 13.1.} For algebras $K'$ and $K$ satisfying Condition \textbf{(II)}, the following assertions are true.

\textbf{1.} Any isomorphism of algebras $K'\to K$ induces the isomorphism of partially ordered sets $\overline{Y}\to \overline{X}$.

\textbf{2.} If the ring $R$ satisfies Condition $(1)$ or Condition $(2)$ before Theorem 9.1, then every isomorphism $K'\to K$ induces isomorphism of preordered sets $Y\to X$.

\textbf{Proof.}\\
\textbf{1.} We fix some isomorphism $\varphi\colon K'\to K$. It is triangular, since Condition \textbf{(II)} holds by Lemma 12.1. If $\varphi=\begin{pmatrix}\alpha&0\\ \delta&\beta\end{pmatrix}$, then $\alpha$ is an isomorphism of algebras $L_1\to L$. For every $y\in Y$, there is $x\in X$ such that $\alpha(f_y)=e_x$ (see the end of Section 5, Proposition 6.1 and the text before it). If $t\in Y$, $t\ne y$, and $\alpha(f_t)=e_s$, then it follows from $f_yf_t=0$ that $e_xe_s=0$; therefore, $x\ne s$.

Now let $z\in X$. Since the inverse isomorphism $\varphi^{-1}$ is of the form $\begin{pmatrix}\alpha^{-1}&0\\ *&*\end{pmatrix}$, we have that $\alpha^{-1}(e_z)=f_y$ for some $y\in Y$ and $\alpha(f_y)=e_z$. Therefore, that $\alpha$ induces the bijection of partially ordered sets $\overline{\tau}\colon \overline{Y}\to\overline{X}$. Namely, if $\alpha(f_y)=e_x$, then $\overline{\tau}(y)=x$ and, consequently, $\alpha(f_y)=e_{\overline{\tau}(y)}$. (Once again, we pay attention to the agreement in Section 3 about the subscripts $x,y,\ldots$ in the symbols $e_x,f_y,\ldots$~.)

We show that $\overline{\tau}$ is an automorphism of partially ordered sets. Let $y,z\in \overline{Y}$ and $y<z$. Then $(M_1)_{yz}=f_yM_1f_z\ne 0$. We have $\varphi(f_y)=\alpha(f_y)+\delta(f_y)=e_{\overline{\tau}(y)}+\delta(f_y)$; similarly, we have $\varphi(f_z)=e_{\overline{\tau}(z)}+\delta(f_z)$. There is an inner automorphism $\mu$ of the algebra $K$ such that 
$$
\mu(e_{\overline{\tau}(y)}+\delta(f_y))=e_{\overline{\tau}(y)},\quad
\mu(e_{\overline{\tau}(z)}+\delta(f_z))=e_{\overline{\tau}(z)}
$$
by Theorem 5.5. Now we have relations
$$
\mu\varphi(M_1)_{yz}=\mu\varphi(f_yM_1f_z)=e_{\overline{\tau}(y)}Me_{\overline{\tau}(z)}=M_{\overline{\tau}(y)\overline{\tau}(z)}
\ne 0.
$$
Therefore, $\overline{\tau}(y)<\overline{\tau}(z)$.

Conversely, if it is given that $\overline{\tau}(y)<\overline{\tau}(z)$, then by a similar method, using the isomorphism $\varphi^{-1}$, we can obtain the inequality $y<z$. All of the above leads to the idea that $\overline{\tau}$ is an isomorphism of partially ordered sets $\overline{Y}\to\overline{X}$.

\textbf{2.} We assume again that $\varphi=\begin{pmatrix}\alpha&0 \\ \delta&\beta\end{pmatrix}\colon K'\to K$ is an isomorphism. It follows from the proof of \textbf{1} that $\varphi$ induces an isomorphism $\alpha\colon L_1\to L$. The latter in turn induces an isomorphism between $R_y'$ and $R_{\overline{\tau}(y)}$ for every $y\in Y$ (it is assumed that we have the equality $L_1=\prod_{y\in Y}R_y'$ which is similar to the equality $L=\prod_{x\in X}R_x$ from Proposition 3.1). If $R_y'\cong M(n,R)$ and $R_{\overline{\tau}(y)}\cong M(m,R)$, then it follows from the $(n,m)$-condition $(2)$ or Condition $(1)$ that $n=m$. This implies the equivalence of the intervals $[y]$ in $Y$ and $[\overline{\tau}(y)]$ in $X$. Hence, the isomorphism $\overline{\tau}\colon\overline{Y}\to\overline{X}$ can be lifted to the isomorphism of preordered sets $\tau\colon Y\to X$.~$\square$

The second theorem gives a positive answer to the question \textbf{(c)} for our rings $K'$ and $K$. The proof of this theorem does not differ fundamentally from the proof of Theorem 6.2(1). For the reasons stated before Theorem 13.1, we prove this theorem as well.

\textbf{Theorem 13.2.} Let the algebras $K'$ and $K$ satisfy Condition \textbf{(II)}. Then every isomorphism between $K'$ and $K$ is equal to the composition of a diagonal isomorphism from $K'$ to $K$ and an inner automorphism of the algebra $K$.

\textbf{Proof.} We take an arbitrary isomorphism $\varphi\colon K'\to K$. Let $\overline{\tau}\colon \overline{Y}\to\overline{X}$ be the isomorphism from the proof of Theorem 13.1(1). We form an element $v=(v_{st})$ of $K$, where $v_{st}=\varphi(f_{\overline{\tau}{\,}^{-1}(t)})(s,t)$ for any $s,t\in X$. We remark that $ve_x=\varphi(f_{\overline{\tau}{\,}^{-1}(x)})e_x$ for every $x\in X$. We have $v_{tt}=1$, and therefore, the element $v$ is invertible in $K$ (Proposition 4.1).

So, for every $z\in Y$, we have relation $ve_{\overline{\tau}(z)}=\varphi(f_z)e_{\overline{\tau}(z)}$. In addition, the relation $\varphi(f_z)v=\varphi(f_z)e_{\overline{\tau}(z)}$ holds. Thus, we have relations
$$
\varphi(f_z)v=ve_{\overline{\tau}(z)},\qquad 
v^{-1}\varphi(f_z)v=e_{\overline{\tau}(z)}
$$
for every $z\in Y$. Let $\mu$ be the inner automorphism of the algebra $K$ defined by the element $v$. Then we have the relation $\mu\varphi(f_z)=e_{\overline{\tau}(z)}$, $z\in Y$. By setting $\gamma=\mu\varphi$, we obtain $\varphi=\mu^{-1}\gamma$, where $\gamma$ is a diagonal isomorphism $K'\to K$, and $\mu^{-1}$ is an inner automorphism of the algebra $K$.~$\square$

\textbf{Corollary 13.3.} Let the algebras $I(Y,R)$ and $I(X,R)$ satisfy Condition \textbf{(II)} and let $Y$, $X$ be partially ordered sets. Then it follows from isomorphism $I(Y,R)\cong I(X,R)$ that $Y\cong X$.

Some rings $R$ such that the algebras $I(Y,R)$ and $I(X,R)$ satisfy Condition \textbf{(II)}, are listed at the end of Section 12.

\textbf{Remarks.} If $X$ is a finite preordered set, then the incidence ring $I(X,R)$ is often called \textsf{structural matrix ring} (see the end of Section 3). It is well written about these rings in the introduction to \cite{DasW96}. In this paper, the following theorem is proved. Let $R$ be a semiprime Noetherian ring and let $B_1$, $B_2$ be two Boolean matrices of order $n$ such that $M(n,B_1,R)\cong M(n,B_2,R)$. Then there is a permutation $\tau$ with $B_2=\tau B_1$.

Isomorphisms of algebras $I(X,R)$ with finite set $X$ are studied in \cite{AbrHR99}, \cite{AbrHR02}, \cite{AnhW13}. The paper \cite{Vos80} contains a theorem which similar to Theorem 13.1, for the case when the factor ring $R/J(R)$ is indecomposable (see item \textbf{e} at the end of Section 12).
In \cite{Fro85}, the isomorphism problem is studied for direct sums of commutative indecomposable rings. The paper \cite{Ler82} is devoted to some analogue of the isomorphism problem for incidence algebras M\"obius categories.

The author of preprint \cite{Wak16} defined algebras $I^n(X,R)$ and called them <<$n$-th partial flag incidence algebras>>; these algebras generalize incidence algebras for $n=2$. In \cite{Khr22}, for finite partially ordered sets $Y$, $X$ and a commutative indecomposable ring $R$, it is proved that $Y\cong X$ provided $I^3(Y,R)\cong I^3(X,R)$.

Other generalizations of incidence rings are also known. For example, isomorphisms of <<formal matrix incidence rings>> are considered  in \cite{Tap17} and \cite{Tap18}. In \cite{For22}, several properties of <<skew incidence rings>> are found and isomorphism problem is solved for these rings.

\section*{\textbf{CHAPTER 4. DERIVATIONS OF \\INCIDENCE ALGEBRAS}}\label{chapter4}
\addtocontents{toc}{\textbf{CHAPTER 4. DERIVATIONS OF INCIDENCE ALGEBRAS}\par}

\section[Representation of Derivations by Matrices]{Representation of Derivations\\ by Matrices}\label{section14}

We use all the material in Chapters 1, 2 and 3. In particular, the incidence algebras $I(X,R)$ are usually denoted by the letter $K$. We do not impose any preconditions on the locally finite preordered set $X$, the ring $R$, which is a $T$-algebra, and the algebra $K$ itself.

A description of the derivation group of the algebra $K$ is obtained (sections 17--19). We note that in sections 14--16 we consider that $K$ is either an incidence algebra or an formal matrix algebra with zero trace ideals. Proofs are carried out for incidence algebras, but they allow transfer to formal matrix algebras (with appropriate corrections). The technique of working with formal matrix algebras with zero trace ideals is well reflected in \cite{KryT21}, \cite{KryT22}.

Here's what's interesting. Comparing the results and proofs in Chapters 2 and 4, it is not difficult to discover their certain similarity and parallelism. A lot of things in Chapter 4 can be obtained if everything <<multiplicative>> is converted to <<additive>>. Perhaps this indicates in favour of the fact that there is a class of certain linear maps of incidence rings, including both automorphisms and derivations.

Let $A$ be an algebra over some commutative ring $T$. A mapping $d\colon A\to A$ is called a \textsf{derivation} of the algebra $A$ if $d$ is a linear mapping, i.e. endomorphism of the $T$-module $A$, and $d(ab)=d(a)b+ad(b)$ for all $a,b\in A$. All derivations of the algebra $A$ form a $T$-module. We denote it by $\text{Der }A$. 

For an element $c\in A$, we define a mapping $d_c$ from $A$ to $A$ by assuming that $d_c(a)=ac-ca$, $a\in A$. Then $d_c$ is a derivation called \textsf{inner}. One says that $d_c$ \textsf{is defined by an element $c$}. Inner derivations of the algebra $A$ form a submodule of the $T$-module $\text{Der }A$. We denote it by $\text{In }(\text{Der }A)$. For designation the factor module $\text{Der }A/\text{In }(\text{Der }A)$, we use the symbol $\text{Out }A$. The elements of $\text{Out }A$ are called \textsf{outer derivations} of the algebra $A$, and $\text{Out }A$ is the \textsf{module of outer derivations} of the algebra $A$.

There is a notion of a derivation in more general form. Let $M$ be an $A$-$A$-bimodule. A \textsf{derivation} of the algebra $A$ with values in the bimodule $M$ is the homomorphism of $T$-modules $d\colon A\to M$ which satisfies the equality $d(ab)=d(a)b+ad(b)$ for all $a,b\in A$. Such a derivation $d$ is said to be \textsf{inner} if there is an element $c\in M$ such that $d(a)=ac-ca$, $a\in A$.

It is known that the $T$-module $\text{Der }A$ has the  structure of a Lie algebra. In this algebra, the multiplication $\circ$ is defined by the relation $d_1\circ d_2=d_1d_2 -d_2d_1$. Here we focus on additive structure of this algebra.

We assume that the algebra $A$ is a splitting extension of its ideal $M$ with the use of some of the algebra $L$, i.e. $A=L\oplus M$, where $\oplus$ is the sign of the group direct sum. The ideal $M$ is a natural $L$-$L$-bimodule. Additionally, $M$ is a non-unital algebra.

We know that every incidence algebra is such an extension. A formal matrix algebra with zero trace ideals also can be represented in the form of indicated splitting extension (see \cite{KryT21} and \cite{KryT22}).

We take an arbitrary derivation $d$ of the algebra $A=L\oplus M$. As every additive endomorphism, $d$ can be represented by the matrix $\begin{pmatrix}\alpha&\gamma\\ \delta&\beta\end{pmatrix}$ with respect to the direct decomposition $A=L\oplus M$ (about this see Sections 5 and 12). Here
$$
\alpha\colon L\to L,\;\beta\colon M\to M,\;\gamma\colon M\to L,\;\delta\colon L\to M
$$
is a $T$-module homomorphisms and
$$
d(a+b)=\begin{pmatrix}\alpha&\gamma\\ \delta&\beta\end{pmatrix}
\begin{pmatrix}a\\b\end{pmatrix}=
\begin{pmatrix}\alpha(a)+\gamma(b)\\ \delta(a)+\beta(b)\end{pmatrix}
$$
for all $a\in L$ and $b\in M$.

We do not distinguish a derivation $d$ and the matrix corresponding to it. Sometimes we write a <<triangular derivation $d$>> if $\gamma=0$, and a <<diagonal derivation $d$>> if $\gamma=0=\delta$.

In the case of triangular derivations, it is possible to obtain very meaningful information about the group $\text{Der}A$. We will soon see that for incidence algebras and formal matrix algebras with zero trace ideals, this is always the case.

So, let $x+y,\,s+t\in A=L\oplus M$. We write down relation
$$
d((x+y)(s+t))=d(x+y)\cdot(s+t)+(x+y)\cdot d(s+t).\eqno(*)
$$
Sequentially assigning values to arguments
$$
y=0=t,\; x=0=s,\;x=0=t,\; s=0=y
$$
and calculating the left and right parts in $(*)$, we get the equalities
$$
\alpha(xs)=\alpha(x)s+x\alpha(s),\quad \delta(xs)=\delta(x)s+x\delta(s),
$$
$$
\gamma(ys)=\gamma(y)s,\quad \gamma(xt)=x\gamma(t),
$$
$$
\beta(yt)=\beta(y)t+y\beta(t)+ \gamma(y)t+y\gamma(t),
$$
$$
\beta(xt)=\alpha(x)t+x\beta(t)+ \delta(x)t,
$$
$$
\beta(ys)=\beta(y)s+y\alpha(s)+y\delta(s).
$$
The first two relations mean that $\alpha$ is a derivation of the algebra $L$ and $\delta$ is a derivation of the algebra $L$ with values in the bimodule $M$. It follows from the third relation that $\gamma$ is an $L$-$L$-bimodule homomorphism.

The converse is also true. Let
$$
\alpha\colon L\to L,\; \beta\colon M\to M,\; \gamma\colon M\to L,\; \delta\colon L\to M
$$
be $T$-module homomorphisms such that the above relations hold. Then the transformation of the algebra $A$ defined by the matrix $\begin{pmatrix}\alpha&\gamma\\ \delta&\beta\end{pmatrix}$, i.e.
$$
\begin{pmatrix}\alpha&\gamma\\ \delta&\beta\end{pmatrix}
\begin{pmatrix}x\\y\end{pmatrix}=
\begin{pmatrix}\alpha(x)+\gamma(y)\\ \delta(x)+\beta(y)\end{pmatrix},\;
x\in L,\;y\in M,
$$
is its derivation.

Now we will see that any derivation of an incidence algebra or a formal matrix algebra with zero trace ideals is triangular. As established in Section 3, an incidence algebra is a splitting extension of $L\oplus M$. A formal matrix algebra with zero trace ideals has a similar property, see \cite{KryT21}, \cite{KryT22}. Next, the letter $K$ denotes one of these two algebras. Let $d=\begin{pmatrix}\alpha&\gamma\\ \delta&\beta\end{pmatrix}$ be some derivation of the algebra $K$.

\textbf{Lemma 14.1.} The relation $\gamma=0$ holds.

\textbf{Proof.} We assume that $K$ is an incidence algebra. We assume that $\gamma\ne 0$, where, as we know, $\gamma$ is an $L$-$L$-bimodule homomorphism $M\to L$. There is an idempotent $e_x$ such that $e_x(\gamma M)=\gamma(e_xM)\ne 0$. Next we have
$$
\gamma(e_xM)=e_x\gamma(e_xM)\subseteq e_xL=R_x\; \text{ and }\; \gamma(e_xM)\subseteq e_xL=R_x.
$$
The last inclusion leads to contradictory relations
$$
0=\gamma(e_xMe_x)=\gamma(e_xM)e_x\ne 0.
$$
Therefore, that $\gamma=0$.~$\square$

\textbf{Remark 14.2.} In the above proof, we implicitly use that the algebra $K$ has at least two distinct idempotents $e_x$ (they appeared in Section 3). Otherwise, the situation  is degenerating, i.e. $K$ is simply some matrix ring $M(n,R)$.

We summarize this section.

\textbf{Corollary 14.3.} Every derivation $d$ of an incidence algebra or a formal matrix algebra with zero trace ideals is triangular, $d=\begin{pmatrix}\alpha&0\\ \delta&\beta\end{pmatrix}$. Here $\alpha$ is a derivation of the algebra $L$, $\beta$ is a derivation of the algebra $M$, and $\delta$ is a derivation of the algebra $L$ with values in the bimodule $M$. 

\section{Some Properties of Derivations}\label{section15}

As it was agreed in the previous section, the letter $K$ denotes either an incidence algebra or a formal matrix algebra with zero trace ideals. However, the notation used in this case refers to incident algebras. By Corollary 14.3, every derivation of both one and the other algebra is triangular. We usually mean the equality $K=L\oplus M$.

We focus on inner derivations of the algebra $K$. We denote by $\text{In}_0(\text{Der }A)$ (resp., $\text{In}_1(\text{Der }A)$) the submodule of inner derivations defined by elements of $L$ (resp., $M$).

\textbf{Lemma 15.1.} There is a direct sum
$$
\text{In }(\text{Der }A)=\text{In}_0(\text{Der }A)\oplus
\text{In}_1(\text{Der }A).
$$

\textbf{Proof.} Not quite obvious that the defined submodules have the zero intersection. Let $b\in L$, $c\in M$ and $d_b=d_c$ (these symbols are defined in Section 14). For any $x\in X$, relations
$$
d_b(e_x)=e_xb-be_x=d_c(e_x)=e_xc-ce_x
$$
hold. Therefore,
$$
e_xc-ce_x=0,\quad e_xc=ce_x=e_xce_x=0.
$$
Therefore, $c=0$ and $d_b=d_c=0$.~$\square$

Now we define one homomorphism and several derivation modules. Namely, we denote by $f$ the module homomorphism $\text{Der }K\to\text{Der }L$ such that $f(d)=\alpha$ for every derivation $d=\begin{pmatrix}\alpha&0\\ \delta&\beta\end{pmatrix}$. The kernel $\text{Ker }f$ consists of derivations of the form $\begin{pmatrix}0&0\\ \delta&\beta\end{pmatrix}$. The image of this homomorphism is denoted by $\Omega$. Next, let $\Lambda$ be the submodule of diagonal derivations, i.e. derivations of the form $\begin{pmatrix}\alpha&0\\0&\beta\end{pmatrix}$. The derivations, which are representable by matrices $\begin{pmatrix}0&0\\0&\beta\end{pmatrix}$, are called \textsf{additive}. The letter $\Psi$ denotes the submodule of all additive derivations, and $\Psi_0$ is the submodule of inner derivations defined by central elements of the ring $L$ (these two submodules are considered in Sections 18 and 19). The last of the submodules is the submodule $\Phi$ equal to $\left\{\begin{pmatrix}\alpha&0\\ \delta&\beta\end{pmatrix}\,\Big|\,\alpha\in\text{In }(\text{Der }L) \right\}$.

The following inclusions hold:
$$
\text{Ker }f\subseteq \Phi,\quad \Psi\subseteq \Phi,\quad 
\text{In}_0(\text{Der }K)\subseteq \Lambda,
$$
$$
\text{In}_1(\text{Der }K)\subseteq\text{Ker }f,\quad 
\text{In }(\text{Der }K)\subseteq \Phi . 
$$
As will be seen from the following, information about the submodules introduced is extremely important for understanding the structure of the entire module $\text{Der}K$.

We highlight the following research directions on the problem of finding the structure of the derivation module $\text{Der}K$.

\textbf{1.} Calculation of the module $\Omega$.

\textbf{2.} Calculation of the submodules $\Psi$ and $\Psi_0$.

\textbf{3.} Calculation of the submodule $\Phi$.

If $\begin{pmatrix}\alpha&0\\ \delta&\beta\end{pmatrix}$ is a derivation of the algebra $K$, then $\alpha$ is a derivation of the algebra $L$. Since $L$ is the product of algebras $R_x$, $x\in X$ (Proposition 3.1), we need information about derivations of ring products. They are arranged quite simply.

\textbf{Proposition 15.2.} Let $S_i$, $i\in I$, be two algebras, $S=\prod_{i\in I}S_i$, and let $d$ be a derivation of the algebra $S$. Then $d(S_i)\subseteq S_i$ for any $i\in I$. In addition, $d$ acts on the product $\prod_{i\in I}S_i$ coordinate-wise and there is a module isomorphism $\text{Der }S=\prod_{i\in I}\text{Der }S_i$.

\textbf{Proof.} We denote by $e_i$ the identity element of the ring $S_i$. Then
$$
d(e_i)=d(e_i)e_i+e_id(e_i)\in S_i.
$$
For any element $a\in S_i$, we have
$$
d(a)=d(e_ia)=d(e_i)a+e_id(a)\in S_i.
$$
Therefore, $d(S_i)\subseteq S_i$, and we can set $d_i=d|_{S_i}$, where $d_i\in \text{Der }S_i$.

The coordinate-wise action of the derivation $d$ means that $d(a)=(d_i(a_i))$ for any element $a=(a_i)\in S$.

We take an arbitrary element $a=(a_i)\in S$. We fix a subscript $k\in I$ and write down $a=a_k+b$, where $b\in (1-e_k)S$. With the use of the idempotent $1-e_k$, it is easy to obtain that the derivation $d$ leaves the ring $\prod_{i\ne k}S_i$ in place.
We write down 
$$
d(a)=d(a_k)+d(b),\; d(a)=(c_i)=c_k+g.\; \text{where} \;g\in(1-e_k)S.
$$
Then $d(a_k)=c_k$, which confirms the fact of the coordinate action of $d$.

Now suppose that for each $i$, we have a derivation $d_i$ of the algebra $S_i$. Assuming $d(a)=(d_i(a_i))$ for the element $a=(a_i)\in S$, we obtain a derivation $d$ of the algebra $S$.

From all the above, we obtain that we have the canonical isomorphism $\text{Der}S\cong\prod_{i\in I}\text{Der}S_i$.~$\square$

For an algebra $L$ equal to $\prod_{x\in X}R_x$, we can write down the following useful fact.

\textbf{Corollary 15.3.} If $d=\begin{pmatrix}\alpha&0\\ \delta&\beta\end{pmatrix}$ is a derivation of the algebra $K$, then for the derivation $\alpha$ of the algebra $L$, all assertions of Proposition 15.2 are true. Therefore, relations $\alpha(e_x)=0$ and $d(e_x)=\delta(e_x)$ are true for every $x\in X$.

In the rest of the section, we pay attention to diagonal derivations. For them, Corollary 14.3 allows for amplification.

\textbf{Corollary 15.4.} Let $d=\begin{pmatrix}\alpha&0\\ 0&\beta\end{pmatrix}$ be a derivation. Then $\alpha$ is a derivation of the algebra $L$, $\beta$ is a derivation of the algebra $M$, and relations
$$
\beta(xt)=\alpha(x)t+x\beta(t),\quad \beta(ys)=\beta(y)s+y\alpha(s) 
$$
hold for all $x,s\in L$ and $y,t\in M$. The converse is also true. If some endomorphisms $\alpha$ and $\beta$ of $T$-modules $L$ and $M$, respectively, satisfy the above properties, then $\begin{pmatrix}\alpha&0\\ 0&\beta\end{pmatrix}$ is a derivation of the algebra $K$.

\textbf{Proposition 15.5.} A derivation $d=\begin{pmatrix}\alpha&0\\ \delta&\beta\end{pmatrix}$ is a diagonal derivation if and only if $d(e_x)=0$ for every $x\in X$.

\textbf{Proof.} The derivation $d$ is diagonal if only if $e_xd(a)e_y=0$ for all $a\in L$ and all non-equivalent $x,y\in X$.

We assume that $d(e_x)=0$ for every $x\in X$. Then for any $a,b\in L$, we have
$$
d(e_xa)=d(e_x)a+e_xd(a)=e_xd(a)\; \text{and similarly }\; d(be_y)=d(b)e_y.
$$
We obtain
$$
e_xd(a)e_y=d(e_xa)e_y= d(e_xae_y)=d(0)=0.
$$

Now let $d$ be a diagonal derivation. Then
$$
d(e_x)=\alpha(e_x)+\delta(e_x)=\alpha(e_x)=0.\quad\square
$$

\textbf{Corollary 15.6.} If $d=\begin{pmatrix}\alpha&0\\ 0&\beta\end{pmatrix}$, then the relation
$$
\beta(e_xce_y)=e_x\beta(c)e_y,\quad \text{where}\; c\in M
$$
holds.

\textbf{Proof.} Taking into account Corollary 15.4 and Proposition 15.5, we write down relations
$$
\beta(e_xce_y)=\alpha(e_x)ce_y+ e_x\beta(c)e_y+e_xc\alpha(e_y)=
e_x\beta(c)e_y.\;\square
$$ 

The definitions of the bimodules $M_{xy}$, $W_k$ and the ideals $V_k$ are given in Section 3.

\textbf{Proposition 15.7.} We take an arbitrary diagonal derivation $d=\begin{pmatrix}\alpha&0\\ 0&\beta\end{pmatrix}$ of the algebra $K$. For the derivation $\beta$ of the algebra $M$ we have an inclusion $\beta M_{xy}\subseteq M_{xy}$ for all $x,y\in X$. In addition, $\beta$ acts  coordinate-wise on the product $\prod_{x,y\in X}M_{xy}$. For $k\ge 1$, inclusions $dV_k\subseteq V_k$ and $dW_k\subseteq W_k$ are also valid.

\textbf{Proof.} The inclusion $\beta M_{xy}\subseteq M_{xy}$ follows from the relation $M_{xy}=e_xMe_y$ and Corollary 15.6. 
We can map a vector $(\beta_{xy})$ to the derivation
$\beta$, where $\beta_{xy}=\beta|_{M_{xy}}$. We show that for $b=(b_{xy})\in M=\prod_{x,y\in X}M_{xy}$, we have $\beta(b)=(\beta_{xy}(b_{xy}))$ which means that $\beta$ acts coordinate-wise. Let $\beta(b)=(c_{xy})$ for all $x,y$. We verify that $\beta_{xy}(b_{xy})=c_{xy}$ for all $x,y$. We fix elements $x,y$ and write down $M=M_{xy}\oplus N$, where $N$ is the product of all bimodules $M_{st}$ besides $M_{xy}$. Let $b=b_{xy}+f$, where $f\in N$. Then
$$
\beta(b)=\beta(b_{xy})+\beta(f),\; \text{where}\; \beta(b_{xy})\in M_{xy}.
$$ 
In addition, we write down
$$
\beta(b)=c_{xy}+g,\; \text{where}\; c_{xy}\in M_{xy},\; g\in N.
$$ 
By Corollary 15.6, we have $e_x\beta(f)e_y=\beta(e_xfe_y)=0$. Therefore, $\beta(f)$ has the zero projection in $M_{xy}$. We obtain that $\beta_{xy}(b_{xy})=c_{xy}$, which is required.

The rest of the inclusion follows from what has already been said.~$\square$

\section{Main Decompositions and Isomorphisms for Module $\text{Der }K$}\label{section16}

As before, $K$ denotes an incidence algebra $I(X,R)$ or a formal matrix algebra with zero trace ideals represented in the form $K=L\oplus M$, as in Section 3. We preserve previously accepted notation. It is important that derivations of the algebra $K$ are triangular (Corollary 14.3).

We recall one more very useful equality from Lemma 15.1:
$$
\text{In }(\text{Der }A)=\text{In}_0(\text{Der }A)\oplus
\text{In}_1(\text{Der }A).
$$
We also recall that the modules listed below appeared at the beginning of section 5.

It is great that derivations of the algebra $K$ can be diagonalized in a certain sense (cf. Theorem 6.2).

\textbf{Theorem 16.1.} There are the following module relations.

\textbf{1.} $\text{Der }K=\text{In}_1(\text{Der }K)\oplus \Lambda$.

\textbf{2.} $\text{Ker }f=\text{In}_1(\text{Der }K)\oplus \Psi$.

\textbf{3.} $\Phi=\text{In }(\text{Der }K)+\Psi=\text{In}_1(\text{Der }K)\oplus (\text{In}_0(\text{Der }K)+\Psi)$.

\textbf{Proof.}

\textbf{1.} For a derivation $d=\begin{pmatrix}\alpha&0\\ \delta&\beta\end{pmatrix}$ of the algebra $K$, we define a function $g\in I(X,R)$, by setting
$$
g(x,y)=\begin{cases}
d(e_y)(x,y),\quad\text{if }\;x\le y,\\
0,\quad\text{if }\;x\not\le y.
\end{cases}
$$
We also can write $g=(g_{xy})=(d(e_y)_{xy})$. In addition, $g\in M$, since $d(e_y)\in M$ (Corollary 15.3). So, $g_{xy}=d(e_y)_{xy}$.

For every $x\in X$, the relation $ge_x=d(e_x)e_x$ holds. We write down the element $g$ in the form $(d(e_y)e_y)_{y\in X}$. We also use a similar form for other elements. If $x\ne y$, then 
$$
e_xe_y=0,\; 0=d(e_xe_y)=d(e_x)e_y+e_xd(e_y),\; e_xd(e_y)=-d(e_x)e_y.
$$
We set $d'=d+d_g$, where $d_g$ is the inner derivation of the algebra $K$ defined by the element $g$ (the designation $d_g$ was given at the beginning of Section 4). To verify that $d'$ is a diagonal derivation, we make the following transformations:
$$
e_xg=e_x(d(e_y)e_y)_{y\ne x}=(e_xd(e_y)e_y)_{y\ne x}=
$$
$$
=-((d(e_x)e_y)e_y)_{y\ne x}=-(d(e_x)e_y)_{y\ne x}.
$$
Next we have the relation
$$
d_g(e_x)=e_xg-ge_x=-(d(e_x)e_y)_{y\ne x}-d(e_x)e_x=
$$
$$
=-(d(e_x)e_y)_{y\in X}=-d(e_x).
$$
We obtain that $d'(e_x)=0$, and $d'$ is a diagonal derivation by Proposition 15.5. Thus, $d=-d_g+d'$, where $d_g\in 
\text{In}_1(\text{Der }K)$, $d'\in\Lambda$. We take an arbitrary derivation $d$ from the intersection $\text{In}_1(\text{Der }K)\cap \Lambda$. Let $d=d_g$, where $g\in M$. It follows from Proposition 15.5 that $d(e_x)=0$ for all $x$. Therefore, we obtain
$$
e_xg-ge_x=0 \text{ and } e_xg=e_xge_x=0 \text{ for any }x.
$$
Therefore, $g=0$, $d=0$ and we have proved \textbf{1}.

\textbf{2.} The assertion follows from \textbf{1} and relations
$$
\text{In}_1(\text{Der }K)\subseteq \text{Ker }f \text{ and } 
\text{Ker }f\cap\Lambda=\Psi .
$$

\textbf{3.} We take an arbitrary derivation $d=\begin{pmatrix}\alpha&0\\ \delta&\beta\end{pmatrix}$ from $\Phi$. Let an inner derivation $\alpha$ of the algebra $L$ be defined by an element $c\in L$. Then we have $d-d_c\in \text{Ker }f$, where $d_c$ is the inner derivation of the algebra $K$ defined by the element $c$. Thus,
$$
d\in\text{In}_0(\text{Der }K)+\text{Ker }f,\quad \Phi=\text{In}_0(\text{Der }K)+\text{Ker }f=
$$
$$
=\text{In}_0(\text{Der }K)+(\text{In}_1(\text{Der }K)\oplus\Psi)=
\text{In }(\text{Der }K)+\Psi.
$$

It remains to verify that $\text{In}_1(\text{Der }K)\cap(\text{In}_0(\text{Der }K)+\Psi)=0$. Let $d_b=d_a+\gamma$, where
$$
d_b\in\text{In}_1(\text{Der }K),\; d_a\in\text{In}_0(\text{Der }K),\; \gamma\in\Psi.
$$

Since $d_a,\gamma\in \Lambda$, we have $d_b\in \Lambda$. Therefore, if $d_b=\begin{pmatrix}0&0\\ \delta&\beta\end{pmatrix}$, then $\delta=0$. Thus, $fb-bf=0$ for all elements $f\in L$. In particular, $e_xb-be_x=0$. Therefore, $e_xb=e_xbe_x=0$ for any $x$. Consequently, $b=0$ and $d_b=0$. The equality $\Phi=\text{In}_1(\text{Der }K)\oplus (\text{In}_0(\text{Der }K)+\Psi)$ is proved.~$\square$

We gather several useful equalities and isomorphisms.

\textbf{Proposition 16.2.} The following relations and isomorphisms hold.

\textbf{1.} $\Psi\cap\text{In }(\text{Der }K)=\Psi\cap\text{In}_0(\text{Der }K)=\Psi_0$.

\textbf{2.} $\Lambda/(\text{In}_0(\text{Der }K)+\Psi)\cong
\Omega/\text{In }(\text{Der }L)$.

\textbf{3.} $\Phi/\text{Ker }f\cong\text{In}_0(\text{Der }K)/\Psi_0
\cong\text{In }(\text{Der }L)$.

\textbf{4.} $\Phi/\text{In }(\text{Der }K)\cong\Psi/\Psi_0$.

\textbf{Proof.}

\textbf{1.} We only verify the inclusion $\Psi\cap\text{In }(\text{Der }K)\subseteq\Psi_0$. We take an arbitrary derivation $d$, equal to $d_1+d_0$, where $d\in \Psi$, $d_1\in \text{In}_1(\text{Der }K)$, $d_0\in \text{In}_0(\text{Der }K)$. Let $d_1$ be defined by an element $c\in M$ and let $d_0$ be defined by an element $a\in L$. Then $d$ is defined by an element $a+c$. 

For any element $b\in L$, we have $d(b)=(bc-cb)+(ba-ab)$. Since $d\in\Psi$, we have $ba-ab=0$. Consequently, $a\in C(L)$, i.e. $a$ is a central element of $L$. We also have $bc-cb=0$. In particular, $e_xc-ce_x=0$ for all $x$. Similar to the end of the proof of Theorem 16.1(3), we obtain $c=0$. So, $d_1=0$ and $d=d_0$. Since $a\in C(L)$, we have $d\in\Psi_0$.

\textbf{2.} First, we remark that the submodule $\text{In }(\text{Der }L)$ is contained in $\Omega$. We denote by $\pi$ the canonical epimorphism $\Omega\to\Omega/\text{In }(\text{Der }L)$. The kernel of the homomorphism $\pi f|_{\Lambda}\colon \Lambda\to\Omega/\text{In }(\text{Der }L)$ is equal to $\Phi\cap\Lambda=\text{In}_0(\text{Der }K)+\Psi$ (see the proof of Theorem 16.1(3)).

\textbf{3.} By considering Theorem 16.1 and \textbf{1}, we write down relations
$$
\Psi/\text{Ker }f\cong(\text{In}_0(\text{Der }K)+\Psi)/\Psi\cong 
\text{In}_0(\text{Der }K)/(\text{In}_0(\text{Der }K)\cap \Psi)=
$$
$$
=\text{In}_0(\text{Der }K/\Psi_0\cong\text{In }(\text{Der }L).
$$

\textbf{4.} Again, taking into account Theorem 16.1 and  \textbf{1}, we have relations
$$
\Phi/\text{In }(\text{Der }K)=(\text{In }(\text{Der }K)+\Psi)/\text{In }(\text{Der }K)\cong\Psi/(\text{In }(\text{Der }K)\cap\Psi)\cong\Psi/\Psi_0.\;\square
$$

In connection with Proposition 16.2(4), we remark that some criteria for coinciding $\Psi$ with $\Psi_0$ are found in Sections 18 and 19.

The following statement is directly derived from Theorem 16.1 and Proposition 16.2.

\textbf{Corollary 16.3.} The following isomorphisms are true.

\textbf{1.} $\text{Der }K/\text{Ker }f\cong\Omega\cong\Lambda/\Psi$.

\textbf{2.} $\text{Der }K/\text{Ker }\Phi\cong\Omega/\text{Der }L$.

In conclusion of the section, we note that the modules $\Omega$ and $\Psi$ are studied in the three remaining sections.
 
\section{Structure of Modules $\text{Der }K$ and $\text{Out }K$}\label{section17}

\textbf{In the remaining Sections 17--19, the letter $K$ denotes only some incidence algebra (and not the incidence algebra or a formal matrix algebra, as in Sections 14--16).}\\
It is essentially used here that all $R_x$ are ordinary matrix rings and all $M_{xy}$ are ordinary matrix groups. Therefore, it is impossible to directly transfer the results of sections 17-19 to formal matrix algebras.

We formulate a number of questions concerning the structure of the modules $\Omega$, $\text{Der}K$ and $\text{Out}K$.

\textbf{1.} Which derivations from $\text{Der}L$ belong to $\Omega$?

\textbf{2.} What is the structure of the modules $\Omega$ and $\Omega/\text{In}(\text{Der}L)$?

\textbf{3.} What is the structure of the modules $\text{Der}K$ and $\text{Out}K$?

Satisfactory answers to all three these questions will be given.

\textbf{Proposition 17.1.} Let $H$ be some $T$-algebra and let $\alpha$ and $\gamma$ be derivations of $H$. An endomorphism $\beta$ of the $T$-module $H$ such that 
$$
\beta(ab)=\alpha(a)b+a\beta(b)=\beta(a)b+a\gamma(b)\; \text{ for all}\, a,b\in H 
$$ 
exists if and only if $\gamma-\alpha$ is an inner derivation.

\textbf{Proof.} We assume that the indicated endomorphism $\beta$ exists. By giving the elements $a$ and $b$ in turn the value 1, we obtain relations
$$
\beta(c)=\beta(1)c+\gamma(c)=\alpha(c)+c\beta(1)
$$
for every $c\in H$. Therefore, $(\gamma-\alpha)(c)= c \beta(1)-\beta(1) c$. In other words, $\gamma-\alpha$ is the inner derivation defined by the element $\beta(1)$. 

Now let $\gamma-\alpha$ be the inner derivation defined by the element $d\in H$. Then
$$
(\gamma-\alpha)(c)=cd-dc\; 
\text{ and }\; dc+\gamma(c)=\alpha(c)+cd
$$ 
for any $c\in H$. We define an endomorphism $\beta$ of the $T$-module $H$ by setting
$$
\beta(c)=\alpha(c)+cd=dc+\gamma(c), \; c\in H.
$$
For arbitrary elements $a,b\in H$, we have
$$
\beta(ab)=\alpha(ab)+abd=dab+\gamma(ab)=
$$ 
$$
=\alpha(a)b+a\alpha(b)+abd=dab+\gamma(a)b+a\gamma(b).
$$ 
We also have
$$
\alpha(a)b+a\beta(b)=\alpha(a)b+a\alpha(b)+abd,
$$
$$
\beta(a)b+a\gamma(b)=dab+\gamma(a)b+a\gamma(b).
$$
These equalities imply the required result.~$\square$

At this point, we can repeat the text of the two paragraphs at the beginning of Section 8. There, for the positive integers $k$, $\ell$ and $c$, where $c=\text{the least common multiple}(k,\ell)$, we introduced matrix algebras $P$, $Q$, $H$ and the $P$-$Q$-bimodule $V$ of matrices. In addition, the numbers $\ell'=c/k$ and $k'=c/\ell$ are determined.

Let $\alpha$ and $\gamma$ be derivations of the algebras $P$ and $Q$, respectively. They induce derivations $\overline{\alpha}$ and $\overline{\gamma}$ of the algebra $H$, respectively. Namely, $\overline{\alpha}(A)=(\alpha(a_{ij}))$ for any matrix $A=(a_{ij})\in H$, $a_{ij}\in P$, and similarly for $\overline{\gamma}$. The derivations $\overline{\alpha}$ and $\overline{\gamma}$ are called \textsf{ring block derivations} (they are also called \textsf{induced derivations}). We keep the notation $\alpha$ and $\gamma$ for them. This agreement is already in effect in the following proposition, which extends Proposition 17.1 for matrix rings.

\textbf{Proposition 17.2.} If $\alpha\in\text{Der }P$ and $\gamma\in\text{Der }Q$, then the existence of an endomorphism $\beta$ of the $T$-module $V$ such that relations 
$$
\beta(pa)=\alpha(p)a+p\beta(a)\, \text{ and }\,
\beta(aq)=\beta(a)q+a\gamma(q)
$$ 
hold for all $p\in P$, $q\in Q$ and $a\in V$ is equivalent to the property that $\gamma-\alpha$ is an inner derivation of the algebra $H$.

\textbf{Proof.} Let $\beta$ be a $T$-endomorphism of the module for which the equalities written in the proposition are fulfilled. It is clear that $\beta$ induces a (block) endomorphism $\overline{\beta}$ of a $T$-module $H$, where $\overline{\beta}(A)=(\beta(A_{ij}))$ for each matrix $A=(A_{ij})$ from $H$. It is assumed that the matrix $A$ is represented in block form, i.e. $A_{ij}$ are $\ell'\times k'$ blocks. For $\overline{\beta}$, the equalities from Proposition 17.1 hold. So $\gamma-\alpha$ is the inner derivation of the algebra $H$.

Now we assume that $\gamma-\alpha$ is an inner derivation of the algebra $H$. We denote by $\beta$ the endomorphism of the module $H$ which exists by Proposition 17.1.

Let $e_1,\ldots,e_{\ell'}$ and $f_1,\ldots,e_{k'}$ be diagonal matrix units corresponding to two block partitions of matrices from $H$. We have $\alpha(e_i)=0$ for $i=1,\ldots,\ell'$, $\gamma(f_j)=0$ for $j=1,\ldots,k'$, and $\beta$ satisfies the equalities from Proposition 17.1. Therefore, the endomorphism $\beta$ induces an endomorphism on the $T$-module $e_iHf_j$ for any unequal $i$ and $j$. Multiplications in $H$ are made in block form (this applies to all three block decompositions). Therefore, the restriction of $\beta$ to $e_iHf_j$, i.e. in fact, on $V$, satisfies the equalities from the proposition.~$\square$

The derivation $\alpha$ of the algebra $L$ is contained in $\Omega$ exactly when there is a derivation of $\beta$ of the algebra $M$ satisfying the equalities from Corollary 15.4. Before Proposition 3.1, the formulas of the bimodule multiplication in $M$ and the multiplication in $M$ are given. Given these formulas and the fact of the coordinate action of $\beta$ (Proposition 15.7), it is not difficult to make sure that $\beta$ satisfies equality
$$
\beta(cd)=\beta(c)d+c\beta(d),\eqno(1)
$$ 
where $c\in M_{xz}$, $d\in M_{zy}$ and $x<z<y$, and relations
$$
\beta(ac)=\alpha(a)c+a\beta(c),\quad
\beta(db)=\beta(d)b+d\alpha(b),\eqno(2)
$$
where $a\in R_x$, $c,d\in M_{xy}$, $b\in R_{y}$.

Let $n_x$ be the order of matrices from the rings $R_x$. We set $c_{xy}=\text{the least common multiple}(n_x,n_y)$ for all $x,y\in X$ such that $x<y$. We denote by $H_{xy}$ the matrix ring $M(c_{xy},R)$. It can be represented as a block matrix ring over the rings $R_x$, $R_y$ and as $R_x$-$R_y$-bimodule of block matrices over $M_{xy}$. The derivations rings $R_x$ and $R_y$ we assume ring (block) derivations the rings $H_{xy}$.

\textbf{Proposition 17.3.} A derivation $\alpha=(\alpha_x)$ of the algebra $L=\prod_{x\in X}R_x$ is contained in the module $\Omega$ if and only if $\alpha_y-\alpha_x$ is an inner derivation of the algebra $H_{xy}$ for all $x,y\in X$ such that $x<y$. 

\textbf{Proof.}\\
\textbf{Necessity.} If $\alpha\in \Omega$, then there is derivation $\begin{pmatrix}\alpha&0\\ 0&\beta\end{pmatrix}$ such that $\beta$ satisfies relations $(2)$ by Corollary 15.4. By Proposition 17.2, $\alpha_y-\alpha_x\in\text{In }(\text{Der }H_{xy})$.

\textbf{Sufficiency.} By Proposition 17.2, there exists an endomorphism $\beta_{xy}$ of the $T$-module $M_{xy}$ for any $x$ and $y$ such that relations $(2)$ hold.

Let $\beta$ be an endomorphism of the $T$-module $M$ which maps an element $(d_{xy})$ to $(\beta_{xy}(d_{xy}))$ for every $(d_{xy})\in M=\prod_{x,y\in X}M_{xy}$. This $\beta$ satisfies two relations from Corollary 15.4. In order for the transformation $\begin{pmatrix}\alpha&0\\ 0&\beta\end{pmatrix}$ of the algebra $K=L\oplus M$ to be its derivation, it remains to verify that $\beta$ is a derivation of the algebra $M$. Thus, it came down to checking equality $(1)$.

We fix three elements $x,z,y$ such that $x<z<y$. We introduce into consideration one more matrix ring. We set $H=M(d,R)$, where $d=\text{the least common multiple}(n_x,n_z,n_y)$. The ring $H$ is block matrix ring over each of the rings $H_{xz}$, $H_{zy}$, $H_{xy}$, $R_x$, $R_z$, $R_y$. It can also be considered as a bimodule of block matrices over each of the bimodules $M_{xz}$, $M_{zy}$, $M_{xy}$. The derivations $\alpha_x$, $\alpha_z$, $\alpha_y$ are considered a ring (block) derivations of the algebra $H$. The $T$-endomorphisms $\beta_{xz}$, $\beta_{zy}$, $\beta_{xy}$ are considered as (block) endomorphisms of the bimodule $H$. The last ones satisfy relations with respect to derivations $\alpha_x$ and $\alpha_z$, $\alpha_z$ and $\alpha_y$, $\alpha_x$ and $\alpha_y$ (of the type recorded in Proposition 17.1). 

The differences $\alpha_z-\alpha_x$, $\alpha_y-\alpha_z$ and $\alpha_y-\alpha_x$ are inner derivations of the algebra $H$. 
Naturally, they induce (in the sense disclosed in the proof of Proposition 17.1) the same mappings $\beta_{xz}$, $\beta_{zy}$, $\beta_{xy}$ as appeared above.

Let inner derivations $\alpha_z-\alpha_x$, $\alpha_y-\alpha_z$ and $\alpha_y-\alpha_x$ be defined by elements $d_{xz}$, $d_{zy}$ and $d_{xy}$ of the algebra $H$, respectively. Then $d_{xy}=d_{xz}+d_{zy}$. We write out from the proof of Proposition 17.1 how endomorphisms act endomorphisms $\beta_{xz}$, $\beta_{zy}$ and $\beta_{xy}$:
$$
\beta_{xz}(a)=\alpha_x(a)+ad_{xz}=d_{xz}a+\alpha_z(a),
$$
$$
\beta_{zy}(b)=\alpha_z(b)+bd_{zy}=d_{zy}b+\alpha_y(b),
$$
$$
\beta_{xy}(c)=\alpha_x(c)+cd_{xy}=d_{xy}c+\alpha_y(c),
$$
where $a,b,c\in H$.

We verify the following relation:
$$
\beta_{xy}(b)=\beta_{xz}(a)b+a\beta_{zy}(b)\eqno(*)
$$
in $H$ for arbitrary $a,b\in H$. 
To do this, we transform the right part of the equality $(*)$ to the left part:
$$
\beta_{xz}(a)b+a\beta_{zy}(b)=\alpha_x(a)b+ad_{xz}b+a\alpha_z(b)+abd_{zy}=
$$
$$
=\alpha_x(a)b+a(d_{xz}b+\alpha_z(b))+abd_{zy}=
$$
$$
=\alpha_x(a)b+a(\alpha_x(b)+bd_{xz})+abd_{zy}=
$$
$$
=\alpha_x(a)b+a\alpha_x(b)+abd_{xz}+abd_{zy}=
$$
$$
=\alpha_x(a)b+a\alpha_x(b)+abd_{xy}=
$$
$$
=\alpha_x(ab)+abd_{xy}=\beta_{xy}(ab).
$$
The equality $(*)$ is proved. Multiplications occurring in it are performed over blocks. This implies that relation $(1)$ holds.~$\square$

We pass to the question of the structure of the modules $\Omega$ and $\text{Der}K$. The concept of a ring derivation of the algebra $M(n,R)$ is actually already given above. We clarify that if $\alpha\in\text{Der}R$, then the derivation of the algebra $M(n,R)$, which maps the matrix $(a_{ij})$ to the matrix $(\alpha(a_{ij}))$, is called a \textsf{ring} (or \textsf{induced}) derivation. The following result is known.

\textbf{Theorem 17.4 \cite{Jon95}.} Every derivation of the matrix algebra $M(n,R)$ is the sum of a ring derivation and an inner derivation.

Let $\varepsilon$ be some derivation of the algebra $R$. Then $\varepsilon$ provides a ring derivation $\varepsilon_x$ of the algebra $R_x$ for every $x\in X$. We denote by $\varepsilon_L$ derivation $(\varepsilon_x)$ of the algebra $L$. For any $x,y\in X$ with $x<y$, the derivation $\varepsilon$ similarly induces a mapping $\varepsilon_{xy}$ on the bimodule $M_{xy}$. Let $\varepsilon_M =(\varepsilon_{xy})$ be an endomorphism of the $T$-module $M$ which acts coordinate-wise. We set $\overline{\varepsilon}=\begin{pmatrix}\varepsilon_L&0\\ 0&\varepsilon_M\end{pmatrix}$. It follows from Corollary 15.4 or Proposition 17.3 that $\overline{\varepsilon}$ is a derivation of the algebra $K$. 

The derivations $\varepsilon_L$, $\varepsilon_M$ and $\overline{\varepsilon}$ are also called \textsf{ring derivations}. We denote by $D$ submodule of all ring derivations of the algebra $L$. We denote the submodule of all ring derivations of the algebra $K$ by the same letter $D$.
It is clear that $D\subseteq \Lambda$ and $D\cap\Psi=0$ in $\text{Der }K$, and the module $D$ is canonically isomorphic to $\text{Der }R$. Next, let $D_0$ be the submodule in $\text{Der }L$ consisting of ring derivations $\varepsilon_L$ of the algebra $L$ such that $\varepsilon\in\text{In }(\text{Der }R)$. Let $D_0$ denote a similar submodule in $\text{Der }K$. These submodules are isomorphic to $\text{In }(\text{Der }R)$ .

One simple result will be useful to us.

\textbf{Lemma 17.5.} If $d$ is a ring derivation and an inner derivation of the algebra $M(n,R)$, then it is defined by a scalar matrix. In other words, $d$ is a ring derivation defined by an inner derivation of the algebra $R$.

We can use Lemma 17.5 to directly verify the following lemma.

\textbf{Lemma 17.6.} 

\textbf{1.} There is an isomorphism $D/D_0\cong \text{Out }R$ in $\text{Der }L$.
 
\textbf{2.} $D_0=D\cap \text{In }(\text{Der }L)$ is true in $\text{Der }L$. 

\textbf{3.} In $\text{Der }K$, we have 
$D_0=D\cap\text{In}_0(\text{Der }K)=D\cap(\text{In}_0(\text{Der }K)+\Psi)$, 
$$
(\Psi\oplus D)\cap\text{In}_0(\text{Der }K)=\Psi_0\oplus D_0.
$$ 

We recall that in Section 2 we agreed to consider a preordered set $X$ as a directed graph.

We formulate the main result of the chapter. It is useful to compare it with Theorem 9.1

\textbf{Theorem 17.7.} Let $X$ be a connected set. For the algebra $K$, where $K=I(X,R)$, we have the following relations and isomorphisms.

\textbf{1(a).} $\Omega=D+\text{In }(\text{Der }L)$, $\Omega/\text{In }(\text{Der }L)\cong \text{Out }R$.

\textbf{1(b).} $\Lambda=\text{In}_0(\text{Der }K)+(\Psi\oplus D)$,
$\text{Der }K=\text{In }(\text{Der }K)+(\Psi\oplus D)=$
$$
=\text{In}_1(\text{Der }K)\oplus (\text{In}_0(\text{Der }K)+(\Psi\oplus D)).
$$

\textbf{2.} $\text{Out }K\cong \Psi/\Psi_0\oplus \text{Out }R$.

\textbf{Proof.}

\textbf{1(a).} We take an arbitrary derivation $\alpha=(\alpha_x)$ in $\Omega$. It follows from Theorem 17.4 that for every $x\in X$, we have $\alpha_x=\rho_x+\mu_x$, where $\rho_x$ is a ring derivation and $\mu_x$ is an inner derivation of the algebra $R_x$. We form derivations $\rho=(\rho_x)$ and $\mu= (\mu_x)$ of the ring $L$ and obtain the relation $\alpha=\rho+\mu$ in $\text{Der }L$. Since $\alpha,\mu\in \Omega$, we have $\rho\in \Omega$.

We will prove the following property: up to an inner derivation, it is possible to obtain that $\rho_x=\rho_y$ for all $x,y\in X$.

We fix some element $t\in X$. We take an arbitrary element $x\in X$ and choose a semipath from $t$ to $x$ in $X$. 
By applying Proposition 17.3 several times or by applying the induction with respect to the length of the selected semipath, we can obtain that $\rho_x=\rho_t+\nu_x$ for some inner ring derivation $\nu_x$. We limit ourselves to a short comment. Equalities of the form $\rho_s=\rho_z+\nu_s$ arise in the ring $H_{sz}$ (see Proposition 17.3). We set $c=\text{LCM}(n_t,n_{z_1},\ldots,n_{z_m},n_x)$, where $z_1,\ldots,z_m$ are the vertices of the selected semipath from $t$ to $x$. Then all ring derivations that appear can be considered as derivations of any ring $R_x$.

Now we obtain $\alpha_x=\rho_x+\mu_x=\rho_t+(\nu_x+\mu_x)$ for every $x\in X$. By setting $\gamma=(\nu_x+\mu_x)$ and $\rho=(\rho_x)$, we obtain $\alpha=\rho+\gamma$, where $\rho\in D$, $\gamma\in \text{In }(\text{Der }L)$. Therefore, the relation $\Omega=D+\text{In }(\text{Der }L)$ is proved. By using Lemma 17.6, we obtain $\Omega/\text{In }(\text{Der }L)\cong 
D/D_0\cong \text{Out }R$.

\textbf{1(b).} Let $d=\begin{pmatrix}\alpha&0\\ 0&\beta\end{pmatrix}\in\Lambda$. Then $\alpha\in\Omega$ and we have $\alpha=\rho+\gamma$, where $\rho\in D$ and $\gamma\in\text{In }(\text{Der }L$ (see \textbf{1(a)}). Derivations $\rho$ and $\gamma$ induce a ring derivation and an inner derivation, respectively (it is contained in $\text{In}_0(\text{Der }K)$). We keep the designations $\rho$ and $\gamma$ for them. We denote by $\psi$ the difference $d-\rho-\gamma$. Then $d=\gamma+\psi+\rho$, where $\gamma\in\text{In}_0(\text{Der }K)$, $\psi\in\Psi$, $\rho\in D$. This leads to the first equality and also leads to the second equality if we take into account Theorem 16.1.

\textbf{2.} With the use of Theorem 16,1, Lemma 15.1 and Lemma 17.6, we obtain relations
$$
\text{Out }K=\dfrac{\text{In}_1(\text{Der }K)\oplus\Lambda}{\text{In}_1(\text{Der }K)\oplus\text{In}_0(\text{Der }K)}\cong
$$
$$
\dfrac{\text{In}_0(\text{Der }K)+(\Psi\oplus D)}{\text{In}_0(\text{Der }K)}\cong\dfrac{\Psi\oplus D}{(\Psi\oplus D)\cap\text{In}_0(\text{Der }K)} =
$$
$$
=\frac{\Psi\oplus D}{\Psi_0\oplus D_0}\cong \Psi/\Psi_0\oplus D/D_0\cong \Psi/\Psi_0\oplus\text{Out }R.\quad\square
$$

We have reduced the problem of finding the structure of modules $\text{Der}K$ and $\text{Out}K$ to a similar problem for the modules $\Psi$ and $\Psi/\Psi_0$. The last sections are devoted to these modules.

We extend Theorem 17.7 to incidence algebras $I(X,R)$, where $X$ is an arbitrary preordered set. As before, let $K=I(X,R)$, $X_i$, $i\in I$, be all connected components of the set $X$. We have $K=\prod_{i\in I}K_i$. Next, by considering Proposition 15.2, we obtain relations
$$
\text{Der }K=\prod_{i\in I}\text{Der }K_i,\quad
\text{In }(\text{Der }K)=\prod_{i\in I}\text{In }(\text{Der }K_i),
$$
$$
\text{Out }K=\prod_{i\in I}\text{Out }K_i
$$
and similar relations for $\text{In}_0(\text{Der }K)$ and $\text{In}_1(\text{Der }K)$.

The symbols $L_i$, $M_i$, $\Lambda_i$, $\Omega_i$, $\Psi_i$ and $D_i$ have a clear meaning in relation to the algebra $K_i$. Again, the equalities $L=\prod_{i\in I}L_i$, $M=\prod_{i\in I}M_i$ and similar equalities are valid for the remaining modules $\Lambda$, $\Omega$, $\Psi$ and $D$.

With the use of the above relations it is not difficult to prove the following theorem.

\textbf{Theorem 17.8.} For an arbitrary incidence algebra $K$, there are the following relations.

\textbf{1(a).} $\Omega=D+\text{In }(\text{Der }L)$,
$\Omega/\text{In }(\text{Der }L)\cong \prod_{|I|}\text{Out }R$.

\textbf{1(b).} $\Lambda=\text{In}_0(\text{Der }K)+(\Psi\oplus D)$,\\
$\text{Der }K=\text{In }(\text{Der }K)+(\Psi\oplus D)=
\text{In}_1(\text{Der }K)\oplus(\text{In}_0(\text{Der }K)+(\Psi\oplus D))$.

\textbf{2.} $\text{Out }K\cong\Psi/\Psi_0\oplus\prod_{|I|}\text{Out }R$.

\textbf{Corollary 17.9 \cite{SpiO97}.} Let $X$ be a partially ordered set and let $R$ be a commutative ring. We consider only derivations of the $R$-algebra $K$. Then we have the following relations and isomorphism
$$
\Omega=0,\; D=0,\; \Psi_0=\text{In}_0(\text{Der }K),\; \text{Der }K=\text{In}_1(\text{Der }K)\oplus\Psi,\; \text{Out }K\cong\Psi/\Psi_0. 
$$ 

\section{Additive Derivations}\label{section18}

This section is an illustration of the observation at the beginning of Section 14 about a certain similarity between Chapters 2 and 4. We can say that this section is some additive version of section 10. The submodules $\Psi$ and $\Psi_0$, defined in Section 15, are considered here. They are related by equalities
$$
\Psi\cap\text{In }(\text{Der }K)=\Psi\cap\text{In}_0(\text{Der }K)=\Psi_0
$$
(Proposition 16.2). The derivations from $\Psi$ are said to be \textsf{additive} and derivations from $\Psi_0$ are said to be \textsf{potential}.

We will get some general information about the modules $\Psi$ and $\Psi_0$. In the case of a finite connected set $X$, the module $\Psi_0$ is a direct summand in $\Psi$ and it is possible to calculate it. In addition, the conditions for matching $\Psi$ with $\Psi_0$ are proposed.

We take an arbitrary additive derivation 
$\begin{pmatrix}0&0\\ 0&\beta\end{pmatrix}$. Here $\beta$ is a derivation of the algebra $M$ and relations
$$
\beta(xt)=x\beta(t),\quad \beta(ys)=\beta(y)s\; \text{ for all }\; x,s\in L,\, y,t\in M
$$ 
hold by Corollary 15.4. Therefore, $\beta$ is an endomorphism of the $L$-$L$-bimodule $M$. Conversely, if a mapping $\beta\colon M\to M$ satisfies indicated properties, then $\begin{pmatrix}0&0\\ 0&\beta\end{pmatrix}$ is an additive derivation. 

We recall again that the bimodules $M_{xy}$ are defined in Section 3. They also appear in Section 17. Let we have an additive derivation 
$\begin{pmatrix}0&0\\ 0&\beta\end{pmatrix}$. As in Section 17, the symbol $\beta_{xy}$ denotes the restriction of $\beta$ to $M_{xy}$. Since $\beta$ is a derivation of the algebra $M$, we have the relation
$$
\beta_{xy}(ab)=\beta_{xz}(a)b+a\beta_{zy}(b)\eqno(1)
$$ for any $x,z,y\in X$ with $x<z<y$ and any $a\in M_{xy}$, $b\in M_{zy}$ 
(see the text after of the proof of Proposition 17.2). Since $\beta$ also is an endomorphism of the $L$-$L$-bimodule $M$, we have that $\beta_{xy}$ is an endomorphism of the $R_x$-$R_y$-bimodule $M_{xy}$ for any $x,y\in X$ with $x<y$.

Next, it follows from Proposition 10.1 that for any such $x$ and $y$, there is an element $c_{xy}\in C(R)$ such that $\beta_{xy}(a)=c_{xy}a$, $a\in M_{xy}$. By considering the relation $M_{xz}\cdot M_{zy}=M_{xy}$, we use $(1)$ to obtain the relation
$$
c_{xy}=c_{xz}+c_{zy},\eqno(2)
$$
where $x<z<y$. We obtain that there is a correspondence between a given derivation $\psi\in\Psi$ and a system of central elements $c_{xy}$ $(x,y\in X$, $x<y)$ of the ring $R$. 

Conversely, every system of elements $\{c_{xy}\in C(R)\,|\,x<y\}$, such that relation $(2)$ holds, leads to an additive derivation. 
Specifically, if $b=(b_{xy})\in M$, then we set $\psi(b)=(c_{xy}b_{xy})$ and $\pi(a)=0$ for $a\in L$.

Based on the above, we write down the following assertion.

\textbf{Proposition 18.1.}

\textbf{1.} The is one-to-one correspondence between additive derivations and systems of the elements $\{c_{xy}\in C(R)\,|\,x<y\}$ such that $c_{xy}=c_{xz}+c_{zy}$ for all $x,z,y\in X$, where $x<z<y$.

\textbf{2.} The module $\Psi$ can be embedded in the direct sum $\mathcal{M}$ of copies of the $T$-module $C(R)$, where $\mathcal{M}=|\{(x,y)\,|\,x<y\}|$.

In general, it seems difficult to find the image of the embedding from Proposition 18.1 (see Section 19).

For potential derivations, i.e. for elements of the submodule $\Psi_0$, Proposition 18.1(1) allows for amplification. We take the potential derivation $\psi=\begin{pmatrix}0&0\\ 0&\beta\end{pmatrix}$; let it be defined by an element $v=(v_x)\in L$, where $v_x\in C(R)$ (we identify scalar matrices from $R_x$ with corresponding elements of $R$). For any $x,y$ $(x<y)$ and every $a\in M_{xy}$, we have relations
$$
\beta(a)=av-va=av_y-v_xa=(v_y-v_x)a.
$$

We conclude that for the derivation $\psi$, the system of elements $\{c_{xy}\,|\,x<y\}$ from Proposition 18.1 consists of elements $v_y-v_x$.

\textbf{Proposition 18.2.} Let $\psi$ be an additive derivation. By Proposition 18.1, $\psi$ corresponds to the system of the central elements $c_{xy}$ which satisfy equality $(2)$. The following assertions are equivalent:

\textbf{1)} $\psi$ is a potential derivation;

\textbf{2)} $\psi$ is an inner derivation;

\textbf{3)} there is an element $v=(v_x)$ of $L$ such that $v_x\in C(R)$ and $c_{xy}=v_y-v_x$ for all $x,y$ ($x<y$).

\textbf{Proof.} The implication \textbf{2)}\,$\Rightarrow$\,\textbf{1)} follows from definitions and Proposition 16.2. 

If we have \textbf{3)}, then $\psi$ coincides with the inner derivation defined by the element $v$.~$\square$

As after Proposition 10.2, we further assume that $X$ is finite, connected, partially ordered set.

Based on an arbitrary set $\{c_{xy}\,|\,x,y\in X,\,x<y\}$ of  central elements of the ring $R$, we assign  weights to edges and semipaths of the graph $X$. If $(x,y)$ is an edge, then we set
$$
w(x,y)=c_{xy}\, \text{ for }\, x<y\, \text{ and }\, w(x,y)=-c_{yx}, \text{ for }\, x>y.
$$
Next for a semipath $P=z_1z_2\ldots z_{k+1}$ in $X$, we set
$$
w(P)=w(z_1,z_2)+\ldots +w(z_k,z_{k+1}).
$$
Let we have some additive derivation $\psi$ of the algebra $K$ and let $\{c_{xy}\,|\,x<y\}$ be a system of central elements of the ring $R$ corresponding to $\psi$ by Proposition 18.1. With the use of this system, we obtain weithts of edges and semipaths of the graph $X$ as it is indicated in the previous paragraph. Now we assume that $\psi$ is a potential derivation, i.e. $\psi$ is an inner additive derivation (Proposition 18.2) and $\psi$ is defined by a central element $v=(v_x)$ the rings $L$. In that case for any two elements $x,y\in X$, the weight of any semipath from $x$ to $y$ is equal to $v_y-v_x$. Consequently, this weight does not depend on  the specific semipath.

We choose some spanning tree $T$ of the graph $X$. To each edge $(x,y)$ in $T$, where $x<y$, we assign some element $c_{xy}\in C(R)$. Using these elements, we set the weight $w$ of these edges and semipaths in $T$ in a way similar to the one used above. If the edge $(x,y)$ is not contained in $T$ and $x<y$, then we assume $c_{xy}=w(P)$, where $P$ is the only semipath in $T$ from $x$ to $y$. As a result, we have a system of elements $\{c_{xy}\in(R)\,|\,x,y\in X,\,x<y\}$ and the weights of the semipaths in $X$.

\textbf{Lemma 18.3.} The system $\{c_{xy}\,|\,x<y\}$ of the elements defined above satisfies the equalities $(2)$. Therefore, it gives an additive derivation according to Proposition 18.1(1).

\textbf{Proof.} We take arbitrary elements $x,z,y\in X$ with $x<z<y$. Let $P_{xz}$, $P_{zy}$ and $P_{xy}$ be semipaths in $T$ from $x$ to $z$, from $z$ to $y$ and from $x$ to $y$, respectively. Based on the uniqueness of semipaths in $T$, we can write $P_{xy}=P_{xz}\,P_{zy}$. Then we get
$$
w(P_{xy})=w(P_{xz})+w(P_{zy}),\quad c_{xy}=c_{xz}+c_{zy}.
$$
It remains to refer to proposition 18.1.~$\square$

In this section, the bimodules $M_{xy}$ are nothing but rings  $R$. Using this circumstance, we define analogues of matrix units. Let $x,y\in X$ and $x<y$. We denote by $e_{xy}$ such a function $X\times X\to R$ that $e_{xy}(s,t)=1$ for $s=x$, $t=y$ and $e_{xy}(s,t)=0$ for all remaining pairs $(s,t)$. The functions $e_{xy}$ have the property that $e_{xz}e_{zy}=e_{xy}$ for $x<z<y$.

Recall that we have fixed some spanning tree $T$ of the graph $X$. Let $\Psi_1$ be the submodule of those additive derivations that cancel elements $e_{xy}$ for all edges $(x,y)$ of $T$ such that $x<y$. Or, equivalently, $c_{xy}=0$ for the specified edges $(x,y)$. It means that $\{c_{xy}\in C(R)\,|\,x<y\}$ is a system of elements from Proposition 18.1, corresponding to some derivation from $\Psi_1$.

\textbf{Theorem 18.4.} The submodule $\Psi_0$ is a direct summand of $\Psi$. Namely, there is a direct module decomposition $\Psi=\Psi_1\oplus\Psi_0$.

\textbf{Proof.} We take an arbitrary derivation $\psi\in \Psi$.  Let $\{c_{xy}\in C(R)\,|\,x<y\}$ be the system of elements corresponding to $\psi$ (see Proposition 18.1). For every vertex in $X$, we define some central element of $R$.

In the tree $T$, we fix some root $x_0$. Let $y$ be some adjacent to $x_0$ vertex in $T$. If $x_0<y$, then we set $v_{x_0}=0$ and $v_y=v_{x_0}+c_{x_0y}$. For $x_0>y$, we set $v_y=v_{x_0}-c_{yx_0}$. Next we choose some vertex $z\in T$ adjacent to $y$. We set
$$
v_z=v_y+c_{yz}\, \text{ for }\, y<z\, \text{ and }\, v_z=v_y-c_{zy}\, \text{ for }\, y>z.
$$ 

We proceed in this way until we stop at some hanging vertex in $T$. In particular, we get a semipath in $T$ from $x_0$ to this hanging vertex. We assume that there are vertices in $T$ that are not included in this semipath. Let $s$ be the first vertex in the semipath with a degree greater than two. In this case, starting from the top of $s$, we form another semipath in a similar way. We continue this process until the moment when an element $v_x\in C(R)$ is mapped to each vertex $x$ in $T$ (hence, in $X$).

We form a central element $v=(v_x)$ of $L$. Let $\psi_0$ be the inner derivation of the algebra $K$ induced by the element $v$. Then $\psi_0\in\Psi_0$. We set $\psi_1=\psi-\psi_0$. If $(x,y)\in T$ and $x<y$, then $\psi(e_{xy})=c_{xy}e_{xy}$ and $\psi_0(e_{xy})=(v_y-v_x)e_{xy}$. It follows the from construction procedure of the elements $v_x$ that $c_{xy}=v_y-v_x$. We obtain that the derivations $\psi$ and $\psi_0$ have equal values on the elements $e_{xy}$, where $(x,y)\in T$. Therefore, we obtain $\psi_1\in \Psi_1$. We obtain
$$
\psi=\psi_1+\psi_0,\quad \psi_1\in \Psi_1,\quad \psi_0\in \Psi_0.
$$

Still need to check that $\Psi_1\cap\Psi_0=0$. Let $\xi\in\Psi_1\cap\Psi_0$. Then $\xi(e_{xy})=0$ for every edge $(x,y)\in T$, where $x<y$.
Since $\xi$ is an inner derivation, $\xi$ is defined by some central element $v=(v_x)$ of the ring $L$. We have relations $0=\xi(e_{xy})=(v_y-v_x)e_{xy}$. Therefore, $v_x=v_y$. If the edge $(x,y)$ is not contained in $T$, then with the use of semipath in $T$ from $x$ to $y$, we obtain $v_x=v_y$. Thus, $\xi=0$.~$\square$
 
Designations $m(X)$ and $\lambda(X)$ were given before Corollary 10.6.

\textbf{Corollary 18.5.} The module $\Psi_0$ of potential derivations is isomorphic to the direct sum of $m(X)-\lambda(X)$ of copies of the module $C(R)$.

\textbf{Proof.} We set $\displaystyle{G=\bigoplus_{m(X)-\lambda(X)}C(R)}$. We take an arbitrary derivation $\psi\in\Psi$. Let $\{c_{xy}\,|\,x<y\}$ be a system of the elements corresponding to the derivation $\psi$, by Proposition 18.1. We define a group homomorphism $g\colon \Psi\to G$ by setting
$$
g(\psi)=(c_{xy}),\quad \text{ where }\; (x,y)\in T,\; x<y .
$$
Lemma 18.3 guarantees that the homomorphism $g$ is surjective. In addition, we have $\text{Ker }g=\Psi_1$. Thus, $\Psi/\Psi_1\cong G$, and Theorem 18.4 implies an isomorphism $\Psi/\Psi_1\cong \Psi_0$. Consequently, $\Psi_0\cong G$.~$\square$

\textbf{Corollary 18.6.} The modules $\Psi$ and $\Psi_0$ coincide if and only if for any $\psi\in \Psi$ such that $\psi(e_{xy})=0$ for every edge $(x,y)\in T$, we have $\psi=0$.

\textbf{Corollary 18.7 \cite{ForP21}.} An additive derivation $\psi$ is inner if and only if for any elements $x,y\in X$, weithts of any two semipaths from $x$ to $y$ are equal.

\textbf{Proof.} Let the additive inner derivation $\psi$ be defined by an element $v=(v_x)$. Then the weight of every semipath from $x$ to $y$ is equal to $v_y-v_x$ (see the text before Lemma 18.3).

Now we assume that it is fulfilled the condition for weights of semipaths for the derivation $\psi$. We write down $\psi=\psi_1+\psi_0$, where $\psi_1\in\Psi_1$, $\psi_0\in\Psi_0$ (Theorem 18.4). For every edge $(x,y)$ in $T$, we have $\psi(e_{xy})=\psi_0(e_{xy})$. If $(x,y)$ is an edge not contained in $T$, then $\psi(e_{xy})=c_{xy}e_{xy}$ and $w(x,y)=c_{xy}$. On the other hand, $w(P)=v_y-v_x$, where $P$ is a semipath in $T$ from $x$ to $y$, and the element $v=(v_x)$ defines an inner derivation $\psi_0$. Consequently, $c_{xy}=v_y-v_x$. Thus, the values of $\psi$ and $\psi_0$ coincide on all elements $e_{xy}$. Therefore, $\psi=\psi_0$.~$\square$

We fix some additive derivation $\psi$. Let $C$ be a simple cycle of the graph $X$. We can define its weight $w(C)$ (relative to $\psi$) as the weight of the semipath from any vertex of the cycle to it. When changing the direction of the cycle traversal, the weight changes the sign. That is, the entire cycle is defined up to a multiplier $\pm 1$.

Let $C_1,\ldots,C_k$ be a fundamental system of cycles, associated with the given the given spanning tree $T$ of the graph $X$ (it is known that $k=\lambda(X)$). The following result clarifies \cite[Corollary 3.10]{ForP21}.

\textbf{Corollary 18.8.} An additive derivation $\psi$ is inner if and only if $w(C_1)=\ldots=w(C_k)=0$.

\textbf{Proof.} Necessity of the condition follows from the first paragraph of the proof of Corollary 18.7.

We assume that $w(C_1)=\ldots=w(C_k)=0$. Let $\psi\in \Psi$. As in the proof of Corollary 18.7, let $\psi=\psi_1+\psi_0$. We take an arbitrary edge $(x,y)$, $x<y$, not contained in $T$. It belongs to some cycle $C_i$ (and only to $C_i$). Let $P$ be a semipath in the cycle $C_i$ from $x$ to $y$. We have relations
$$
0=w(C_i)=w(P)-w(x,y),\quad w(x,y)=w(P).
$$
Similar to the proof of Corollary 18.7, we obtain $c_{xy}=v_y-v_x$ and $\psi=\psi_0$.~$\square$

At the end of the section, we give a brief comment about additive or potential derivations. In \cite{SpiO97}, they are introduced as part of a functional approach to incident rings. A function $\sigma\colon X\times X\to R$ is said to be \textsf{additive} if $\sigma(x,y)=\sigma(x,z)+\sigma(z,y)$ for all $x,z,y$ with $x<z<y$ and $\sigma(x,y)\in C(R)$ for all $x,y$, where $x<y$.

If $\sigma$ is some additive function, then the mapping 
$$
d_{\sigma}\colon K\to K,\quad d_{\sigma}(f)=\sigma*f,\; f\in K,
$$
is an additive derivation of the algebra $K$ ($*$ is the Hadamard product and $K$ is $I(X,R)$, as before). Conversely, every additive derivation can be obtained based on some additive function $\sigma$.

Now let $q\colon X\to C(R)$ be an arbitrary mapping. We define a function $\tau_q$, contained in $I(X,R)$, by setting
$$
\tau_q(x,y)=\begin{cases}
q(y)-q(x)\, \text{ for }\, x\le y,\\
0\, \text{ for }\, x\not\le y.
\end{cases}
$$
The function $\tau_q$ is additive; it is called a <<potential function>>. In addition, the additive derivation $d_{\tau_q}$ is potential derivation. Any potential derivation can be obtained in such a manner.

\section{Calculation of Modules $\Psi$ and $\Psi/\Psi_0$}\label{section19}

We continue the study of additive derivations of the algebra $I(X,R)$, started in the previous section. We assume that $X$ is a finite, connected, partially ordered set.

Corollary 18.5 asserts the presence of the isomorphism $\displaystyle{\Psi_0\cong\bigoplus_{m(X)-\lambda(X)}C(R)}$ for the module $\Psi_0$ of potential derivations (i.e. additive inner derivations). In this section, in the case when the center $C(R)$ is a field, we obtain simple formulas for finding the dimensions of the space $\Psi$ of all additive derivations and the space $\Psi/\Psi_0$ of additive outer derivations. We recall that $m(X)$ (or simply $m$) is the number of edges and $\lambda(X)$ is the cyclomatic number of the graph $X$.

Let $e_1,\ldots,e_m$ be the edges of the graph $X$ and let $\Delta_1,\ldots,\Delta_s$ be all cycles of length 3, i.e. triangles. Thus we fix some numbering of edges and triangles.

Now we make one $s\times m$ matrix $P=(p_{ij})$. Its construction is similar to the construction of the well-known matrix of cycles. With the difference that we take into account the orientation of the edges in a certain way. The columns of the matrix $P$ correspond to the edges $e_1,\ldots,e_m$ and the rows correspond to the triangles $\Delta_1,\ldots,\Delta_s$.

We take some triangle $\Delta_t$. Then it is necessarily possible to denote its vertices by $x,z,y$ in such a way that $x<z<y$. Let $(x,y)=e_i$, $(x,z)=e_j$ and $(z,y)=e_k$. We assume that $p_{ti}=1$, $p_{tj}=p_{tk}=-1$. In the remaining positions of the row with integer $t$, we put $0$. We call the resulting matrix $P=(p_{ij})$ the \textsf{matrix of triangular cycles} of the set (or graph) $X$.

It follows from Proposition 18.1 that to obtain a specific additive derivatio of $\psi$, one must have a system of elements $\{c_{xy}\in C(R)\,|\,x<y\}$ such that $c_{xy}=c_{xz}+c_{zy}$ as soon as $x<z<y$. Therefore, these elements $c_{xy}$ are solutions over $C(R)$ of a system of $s$ linear homogeneous equations with the matrix $P$. Conversely, any solution of this system gives the elements $c_{xy}$ and the corresponding additive derivation. Let $y_1,\ldots,y_m$ be some variables and $Y=(y_1,\ldots,y_m)^t$. We write down the system of equations discussed in the form of a matrix equation $PY=0.\quad(1)$

\textbf{Corollary 19.1.} The module $\Psi$ of additive derivations is isomorphic to the submodule in $\oplus_mC(R)$ consisting of all solutions of the matrix equation $PY=0$.

Next, we assume that the center of $C(R)$ is a field (this is so, for example, if $R$ is a field). We denote this field by $F$. Let $r(P)$ be the rank of the matrix $P$. It can be called the \textsf{triangular rank} of the graph $X$. It does not depend on the choice of numbering of edges and triangles, i.e. it is an invariant of the graph $X$.

We write down the main result of the section.

\textbf{Theorem 19.2.} Let $X$ be a finite connected partially ordered set and let $C(R)$ be a field. For the algebra $I(X,R)$, the following assertions are true.

\textbf{1.} The following relations for space dimensions hold.

$\text{dim }_F\Psi=m(X)-r(P)$, $\text{dim }_F\Psi_0=m(X)-\lambda(X)$ (Corollary 18.5) and $\text{dim }_F\Psi/\Psi_0=\lambda(X)-r(P)$.

\textbf{2.} Every additive derivation is inner if and only if $r(P)=\lambda(X)$.

The brief considerations in this section can be clarified and supplemented by involving some vector space. We limit ourselves to a comment.

Let $F$ be a field and let $V$ be a vector space with the basis consisting of edges $e_1,\ldots,e_m$. Ones often take $\mathbb{Z}_2$ as $F$. In this case, all Eulerian subgraphs in $X$ form a subspace in $V$ of dimension $\lambda(X)$. (This subspace of $C(X)$ is called the \textsf{cycle space}.) Every fundamental system of cycles is the basis of the space $C(X)$. However, in our situation, the field $\mathbb{Z}_2$, apparently, is not suitable. For example, you can take the field $\mathbb{Z}_3$. Or, if $C(R)$ is a field of characteristic not equal to $2$, you can put $F=C(R)$.

\textbf{Remarks.} In \cite{Khr12}, a description of the Lie algebra of outer derivations of a finitary incidence algebra is obtained in terms of a certain cohomology group. The content of the paper \cite{KayKW19} is generally clear from its title. The same can be said for \cite{GhaG23}.

\addtocontents{toc}{\textbf{Bibliography$\qquad\qquad\qquad\qquad\qquad\qquad\qquad\qquad\qquad\qquad\;\;\,$ \pageref{biblio}}\par}

\label{biblio}

\end{document}